\documentclass[12pt,twoside]{amsart}
\usepackage{amssymb}
\usepackage{enumerate}

\newtheorem{theorem}{Theorem}[section]
\newtheorem{thm}[theorem]{Theorem}
\newtheorem{lem}[theorem]{Lemma}
\newtheorem{cor}[theorem]{Corollary}
\newtheorem{prop}[theorem]{Proposition}
\newtheorem{defi}[theorem]{Definition}
\newtheorem{rem}[theorem]{Remark}
\theoremstyle{definition}

\theoremstyle{remark}

\numberwithin{equation}{section}
\def \Bbb{\mathbb}

\def\onto{{\kern3pt\to\kern-8pt\to\kern3pt}}

\def\<{\langle}
\def\>{\rangle}
\def\|{{\ |\ }}

\def\onto{\twoheadrightarrow}

\def\-{\underline}

\def\N{\Bbb N}

\def\R{\Bbb R}
\def\C{\Bbb C}
\def\B{\Bbb B}
\def\Cn{\Bbb C^n}

\def\M {\Bbb M}
\def\H{\Bbb H}
\def\E{\Bbb E}

\def\<{\langle}
\def\>{\rangle}

\catcode`\@=11
\def\serieslogo@{\relax}
\def\@setcopyright{\relax}
\catcode`\@=12

\title[The $\overline{\partial}$-equation
 on a class of convex domains]
{  Optimal  Lipschitz Estimates for the $\overline\partial$ equation \\
on a class of convex domains}

\begin{document}

\author{Vi\d{\^e}t Anh  Nguy\^en}
\address{Vi\d{\^e}t Anh Nguy\^en\\
Carl von  Ossietzky Universit\"{a}t Oldenburg \\
Fachbereich  Mathematik\\
Postfach 2503, D--26111\\
 Oldenburg, Germany}
\email{nguyen@mathematik.uni-oldenburg.de}

\author{El Hassan Youssfi}
\address{El Hassan Youssfi\\
L.A.T.P. U.M.R. C.N.R.S 6632 \\
Centre de Math\'ematiques et d'Informatique\\
Rue Joliot--Curie\\
Universit\'e de Provence, F--13453\\
Marseille cedex 13 \\
France}
\email{youssfi@gyptis.univ-mrs.fr}

\subjclass[2000]{Primary 32A22, 32A25}
\date{}

\keywords{ The $\overline{\partial}$-equation, optimal  Lipschitz
estimate,
 Berndtsson formula,   Martinelli-Bochner formula,  Cauchy formula.  }

\begin{abstract}
 In this paper, we consider the Cauchy-Riemann equation $\overline\partial u= f$ in
 a new class of convex domains  in $\C^n.$  We prove that under  $L^p$   data, we can choose
 a solution in the Lipschitz space $\Lambda_{\alpha} ,$  where  $\alpha$
  is an  optimal positive number given  explicitly  in terms of    $p.$
  \end{abstract}


\maketitle



{\indent 1. \ Introduction and statement of the main results \dotfill 1}

{\indent 2. \  The complex manifolds $\H_N$ and $\M_N$ \dotfill 3}

{\indent 3. \ Integral formulas on  the complex manifolds $\M_N$
  \dotfill 10}

{\indent 4. \  Local coordinate systems on  the complex manifolds $\H_n$ and $\H_m$ \dotfill  20}

{\indent 5. \  Reduction of estimates from  $\M_N$  to $\B_N$ \dotfill 22}

{\indent 6. \  Integral kernels \dotfill 28}

{\indent 7. \ Integral estimates \dotfill 30}

{\indent 8. \  Lipschitz estimates  on the complex manifold $\M_N$  \dotfill 35}

{\indent 9. \  Stokes type theorem on the complex manifold $\M_N$ and applications  \dotfill 40}

{\indent 10. \ Proof of the main results \dotfill 45}


\bigskip

\section{ Introduction and   statement of the main results}

 For  every $m$-uplet  of   positive integers    $N:=(n_1,\ldots,n_m),$ we consider  the following domain~:
  \begin{equation}
  \Omega_N:= \left\{  Z=(Z_1,\ldots,Z_m)     \in  \C^{n_1}\times \cdots \times \C^{n_m}:  \sum\limits_{j=1}^{m}
  \left (|Z_j|^{2} + |Z_j\bullet Z_j| \right ) < 1 \right\},\end{equation}
  where $z\bullet w:= \sum\limits_{j=1}^{k}  z_j w_j$
  and  $\vert z\vert := \sqrt{z\bullet \overline{z}},$  for  all elements $z:=(z_1,\ldots,z_k)$ and $w:=(w_1,\ldots,w_k)$ of  $ \C^k.$

  The euclidean   ball of radius $\frac{\sqrt{2}}{2}$  in  $\C^m$
   and the   minimal  ball in $\C^{n_1}$ correspond respectively to
    the cases $n_1=\cdots=n_m=1$  and $m=1.$ The domains $\Omega_N$ were introduced by
    the second author in \cite{Y} where he computed their Bergman and Szeg\"o kernels.
      We should point out
      that  these domains  are convex  but  they are neither
       strictly pseudoconvex nor piecewisely smooth except for the case of the
        euclidean  balls.

    Optimal   estimates  for  the $\overline\partial$-equation    were  considered for the  category of
smooth  domains   by several  authors.
 In \cite{Kr1},  Krantz  obtained  the optimal Lipschitz and $L^p$ estimates for   smooth strongly  pseudoconvex domains.   Later in  \cite{CKM},   Chen, Krantz  and Ma   established that   this  kind of   regularity  holds  for
 smooth complex  ellipsoids.   The  general  case  of   smooth convex domains of finite
  type  was  considered   only  recently
 in the  works of   Cumenge (\cite{Cu1},\cite{Cu2}), Diederich-Fischer-Forn{\ae}ss \cite{DFF},
   Fischer \cite{Fi}  and  Hefer \cite{He}.  The aim of the  present  paper is
 to  study  the  optimal  Lipschitz regularity  for the  $\overline\partial$-equation
  in the  class of  convex domains
 $\Omega_N.$

  To state  the  main results,
we fix some  notations and suppose without loss of  generality
that $n_1\leq \cdots\leq  n_m.$ Since the case of the   euclidean
balls  is well-known, we shall assume  that $ N\not = (1,\ldots,1)$
and   let   $l$ denote the smallest nonnegative integer such that
  $n_{l+1} >1.$
  We set $\vert N\vert:= \sum\limits_{j=1}^{m} n_j.$

    The Lipschitz spaces we use herein are  the classical ones and   those given for   $0< \alpha \leq
    1,$ by
\begin{equation*}
      \Lambda_{\widetilde{\alpha}}(\Omega_N):=\left\{ f\,\, : \Vert f
\Vert_{L^{\infty}(\Omega_N)    } +
\underset{\overset{z,z+h \in \Omega_N}{ 0 < \vert h \vert <
\frac{1}{2}}}{\sup}\frac{\vert f(z+h)-f(z) \vert} { \vert h
\vert^{\alpha} \vert\log{\vert h \vert}  \vert } \equiv  \Vert f
\Vert_{\Lambda_{\widetilde{\alpha} }  (\Omega_N)         }
 <\infty        \right\}.
    \end{equation*}


 The first main result is  the following. It
    generalizes our previous  result  \cite{VY1}:
 \begin{thm}
  Suppose that  $N:=(n_1,\ldots,n_m)$  is as above
   and the domain  $\Omega_N$ is given  by
 (1.1). Let
 $$\alpha=\alpha (N,p):= \left\lbrace
\begin{array}{l}
 \frac{1}{2}-\frac{\vert N\vert +m-l+1}{p},\ \ \text{if}\  N \not =  (2,\ldots,2)\  \text{and} \ p>2(\vert N\vert +m-l+1);\\
 \frac{1}{2}-\frac{3m}{p},\qquad \ \  \    \ \ \text{if}\  N =  (2,\ldots,2)\  \text{and} \ p> 6m.
  \end{array}  \right.$$
  Then
  for every $\overline\partial$-closed $(0,1)$-form  $f$ with coefficients in $L^p(\Omega_N),$ there exists a
function $u$ defined on  $\Omega_N$   that satisfies
$\overline\partial u=f$
 (in the distribution sense) and the   estimate
$$\left\lbrace
\begin{array}{l}
\Vert u\Vert_{\Lambda_{\alpha}(\Omega_N)}           \leq C_p   \Vert f\Vert_{L^p(\Omega_N)},\qquad\ \  \text{if}\  p < \infty;\\
\Vert u\Vert_{\Lambda_{\widetilde{\frac{1}{2}}}(\Omega_N)}  \leq C_{\infty}   \Vert f\Vert_{L^{\infty}(\Omega_N)}
,\qquad  \text{if}\  p = \infty
. \end{array}  \right.$$
\end{thm}

 The following result asserts that   the regularity   in Theorem 1.1 is sharp.

\begin{thm}
Let $N,$ $\Omega_N,$ $p,$ and $\alpha:=\alpha(N,p)$ be as in the  statement of Theorem 1.1.
Then  there exists a $\overline\partial$-closed $(0,1)$-form  $f$ with coefficients in $\mathcal{C}^{\infty}(\Omega_N)$   that  satisfies
   $$\left\lbrace
\begin{array}{l} f\in L^s(\Omega_N),\ \forall s <p, \qquad \text{if}\  p<\infty;\\
   f\in L^{\infty}(\Omega_N),\qquad \qquad \ \ \ \ \text{if}\  p=\infty;
    \end{array}  \right. $$
and  if $u$ is a function  satisfying  $\overline\partial u=f,$  then $u \not\in \Lambda_{\alpha+\epsilon}(\Omega_N),\  \forall
   \epsilon >0.$
\end{thm}

These results have been announced in \cite{VY2}.

   Theorem 1.2 implies that if  $N \not =  (2,\ldots,2)$  and   $ p\leq 2(\vert N\vert +m-l+1)$
or  if  $N =  (2,\ldots,2)$  and  $  p\leq 6m,$
   then  we can not solve  the  $\overline\partial$-equation  on  $\Omega_N$ under
    $L^p$ data with the Lipschitz regularity given above.

We observe that for  $N=(2),$  the domain  $\Omega_{(2)}$ is
linearly biholomorphic to the Reinhardt triangle $\{
(z_1,z_2)\in\C^2:\vert z_1\vert +\vert z_2\vert <1\}.$ The
reduction of our  Theorems 1.1 and 1.2 to this case, compared with
the results obtained for domains of finite type
   (\cite{Kr1},\cite{Cu1},\cite{Cu2},\cite{DFF},\cite{Fi},\cite{He}), shows that  domain $\Omega_{(2)}$  has
     the same gain of smoothness    for the
$\overline\partial$-equation   as
   strictly pseudoconvex smooth domains in $\C^2.$
   Our results show also that  there exist smooth domains of finite type
 for which  the  gain of smoothness for the $\overline\partial$-equation is  worse than  that of the  singular domains    $\Omega_N.$

\medskip

The paper is  organised as follows.

 In Section 2 we  introduce the main tools and prove preliminary
results. The objects used are a complex manifold $\H_N$, its
intersection  $\M_N$ with the euclidean unit ball and a proper holomorphic
mapping $F_N$ relating the $\overline\partial$-equation on $\M_N$
to that on $\Omega_N.$ We establish in this section
 Proposition 2.5    which   gives an  integral  representation  formula of  Berndtsson type
for the  complex manifold  $\M_N.$ From this result we derive in
Section 3  a formula of Martinelli-Bochner type (Theorem 3.1) and
two formulas of  Cauchy type (Theorems 3.3 and 3.6)  for the
complex manifold $\M_N.$ These  integral representations play  a
peculiar role in the construction of the
$\overline\partial$-solving operators on $\M_N$  and   $\Omega_N.$

 In Section 4  we give appropriate local coordinates  on the
  complex  manifold   $\H_N$   which   permit us to
    prove Theorem  5.6 in Section 5. The latter result will be  called   Theorem of reduction
of estimates since from broad outlines,
  it  reduces  certain  integral estimates on
$\M_N$  to   analogous integrals, but  simpler, which are taken  on
some balls of   $\C^{\vert N\vert}.$ This result, combined  with Section 6, allows us
to establish integral estimates in Section 7.

  An operator  solution $T_1$ of the
 $\overline\partial$-equation  on  $\M_N$ is constructed in Section 8 and
 related      Lipschitz  estimates are established there. The  formula for  $T_1$ is explicit
 and contains an integral term  taken  over the boundary $\partial\M_N$  of  $\M_N.$
 In order to  handle  this  term,  we  prove a sort of Stokes theorem  in
Section  9 which allows us to transform these integral estimates
into analogous ones
   taken
over  $\M_N$ and then apply  the Theorem of reduction of
estimates.

  Theorem 1.1 is proved in Section 10.   By means of the  operator
  $T_1$ and
    the proper holomorphic mapping  $F_N,$
 we define  an  operator  $T,
 $ solution of the $\overline\partial$-equation  on the domain
 $\Omega_N$ and
 transfer  the Lipschitz regularity  for  $T_1$
  to that of the operator  $T.$ Finally,   we prove
    Theorem  1.2 by giving   concrete examples to show the sharpness of the  results   of Theorem  1.1. Then we conclude the
 paper by some remarks and open questions.

Throughout the paper, the letter $C$ denotes  a finite constant  that is
not necessarily the same  at each occurence and that depends on $N$ and eventually other parameters.

%
%
%
%
\section{The complex manifolds    $\H_N$ and  $\M_N.$}
In this  section we fix the notations and   prove  some  preliminary  results.
For the  simplicity   of    calculations  we only  consider,  without  loss of  generality,  the case
of the domain
  $\Omega_N$  with  $N=(\underbrace{1,\ldots,1}_{l},n,m),$  where $l,m,n$ are  positive  integers  and
$n,m >1.$   In this  case  we  have
 $\vert N\vert=l+n+m$  and
 $\Omega_N$ can be  written   in the form
$$
\Omega_N:=\left\{  Z=(x,z,w)\in \C^l\times\C^n\times\C^m: 2\vert x\vert^2+\vert z\vert^2+ \vert z\bullet z\vert+
\vert w\vert^2+ \vert w\bullet w\vert  <1  \right\}.
$$
Consider  the  complex  manifold  $\H_N$ given by
$$\H_N:=\left\{  Z =(x,z,w)\in  \C^l  \times\C^{n+1}\setminus\{0\}\times\C^{m+1}             \setminus\{0\}: z\bullet z=w\bullet w=0 \right\}.
$$
Let  $\B_N$  be the  euclidean open  unit  ball in  $\C^{\vert N\vert+2}$   and  $\partial\B_N$  its  boundary. We set
$\M_N:=\H_N \cap \B_N$ and $ \partial\M_N:=\H_N \cap  \partial \B_N .$
We first point  out  that  $\H_N$
and    $\partial\B_N$   are transverse  while   the  variety
$ \{  Z=(x,z,w)\in \C^l\times\C^{n+1}\times\C^{m+1}:z\bullet z=w\bullet w=0 \}$ does not  meet
  $\partial\B_N$  transversally.
Denote  by $dV,dV_l,dV_n$ and $dV_m$    the respective  canonical  measures  on the  complex  manifolds   $\H_N,\C^l,\H_n$ and  $\H_m$.  These
measures are related by the
  following
\begin{prop}
For all compactly supported  continuous functions $f$   on  $\H_N, $ we have
   \begin{equation*}
 \int_{\H_N}f(Z) dV(Z)= C \int_{\C^l}\int_{\H_n}\int_{\H_m} f(x,z,w)dV_l(x)dV_n(z)dV_m(w).
 \end{equation*}
 \end{prop}
   \begin{proof}  Observe that
  \begin{eqnarray*}
 dV(Z) := C \left .\left (\sum_{p=1}^{l} dx_{p}\wedge
d\overline{x}_{p}+\sum_{j=1}^{n+1} dz_{j}\wedge
d\overline{z}_{j}+\sum_{k=1}^{m+1} dw_{k}\wedge
d\overline{w}_{k} \right )^{l+n+m}\right |_{\H_N}.  \\
\end{eqnarray*}
In this formula the constant  $C$   is equal to $\frac{1}{(l+n+m)!} \left ( \frac{i}{2}\right )^{l+n+m}.$  Therefore,
a direct computing  shows that
  \begin{eqnarray*}
 dV(Z) &=& C\left . \left (\sum_{p=1}^{l} dx_{p}\wedge
d\overline{x}_{p}\right )^{l}\left(\sum_{j=1}^{n+1} dz_{j}\wedge
d\overline{z}_{j}\right)^{n}\left (\sum_{k=1}^{m+1} dw_{k}\wedge
d\overline{w}_{k} \right )^{m}\right |_{\H_N}.  \\
&=&  dV_l(x)dV_n(z)dV_m(w).
\end{eqnarray*}
 This  completes  the proof.
\end{proof}
\noindent Let $    \E:=\lbrace t=(t_1,t_2,t_3) \in \rbrack 0,1\lbrack^3:\  t_1^2 + t_2^2 +t_3^2 < 1  \rbrace\ \text{and}\  \partial \E:= \lbrace  t\in \rbrack 0,1\lbrack^3  : \ t_1^2 + t_2^2 +t_3^2 = 1 \rbrace $ its boundary.      Then  the mapping
$ F:\E\times \partial\B_l\times \partial  \M_n\times \partial\M_m\longrightarrow  \M_N$
given  by $  F(t,x,z,w):=tZ=
(t_1 x,t_2 z, t_3 w),$   where
$t=(t_1,t_2,t_3)\quad \text{and}\quad  Z=(x,z,w),$
 is a diffeomorphism. Moreover, it maps  $\partial \E\times \partial\B_l\times \partial  \M_n\times \partial\M_m$  onto  $\partial\M_N.$

Let  $d\sigma_n$ be   the unique  probability measure, $SO(n+1,\R)$-invariant on $\partial\M_n.$   Similarly,
let $d\sigma_m$ be   the unique  probability measure, $SO(m+1,\R)$-invariant  on $\partial\M_m.$
Finally, let $d\sigma_l$ be the surface measure on   $\partial\B_l.$  Combining  Proposition 2.1  of  \cite{VY1}  and Lemma  2.1
of  \cite{MY}, we  obtain
\begin{cor}
 For all compactly supported  continuous functions $f$   on     $\H_n,$  we have
$$\int\limits_{\H_n}  f(z) dV_n(z)=
C\int\limits_{0}^{+\infty} t^{2n-1}\int\limits_{\partial\M_n} f(t\zeta)d\sigma_n(\zeta) dt. $$
\end{cor}
There are obviously  analogous   integral formulas in polar coordinates with $\M_m$ and $\B_l$ in place of  $\H_n.$
We now  define  a  natural  measure  $d\sigma$  on $\partial\M_N$
  by setting
$  d\sigma:= (F_{*}) \left(d\phi\wedge d\sigma_l\wedge d\sigma_n\wedge d\sigma_m\right),$
 where  $d\phi$ is  the  surface  measure of  the unit sphere $\partial \E$.
Using this,  Corollary 2.2 and  integration in polar coordinates, one can establish  the following
\begin{lem}
 For all compactly supported  continuous functions $f$   on  $\H_N,$   we  have
$$\int\limits_{\H_N}  f(Z) dV(Z)=
C(N)\int\limits_{0}^{+\infty} t^{2\vert N\vert-1}\int\limits_{\partial\M_N} f(t\Theta)d\sigma(\Theta) dt. $$
\end{lem}

In what follows we shall establish some integral formulas on $\M_N.$  To do so, we shall approximate
$\M_N$ by appropriate regular varieties which are complete intersections. Then
we  apply  to each of these varieties   the results of  Berndtsson in
\cite{BB}.

 For $ 0< r<1,$  let
 $\mathcal{D}_{r}$ be  the  domain of  $\C^{\vert N\vert+2}$  defined by
 $$  \mathcal{D}_{r}:= \left\{  Z =(x,z,w)\in  \B_N: 
  \vert z\vert
 > r, \vert w\vert >r \right\}.$$
 Note that the boundary of  $\mathcal{D}_{r}$  is  piecewisely smooth.
 We put
$\M_{r}:= \H_N \cap \mathcal{D}_{r}.$
 Let $$s:=(s_1,\ldots,s_{\vert N\vert+2}):\ \overline{ \mathcal{D}}_{r}\times
\overline{ \mathcal{D}}_{r}   \longrightarrow \C^{\vert N\vert+2}$$
be  a $\mathcal{C}^1$  function that satisfies
\begin{eqnarray}
\vert s(\Theta,Z)\vert \leq C \vert \Theta -Z\vert \quad
\text{and}\quad  \vert s(\Theta,Z)\bullet(\Theta -Z)\vert \geq C \vert \Theta -Z\vert^2
\end{eqnarray}
uniformly for  $\Theta\in  \overline{ \mathcal{D}}_{r}  $ and for  $Z$  in any  compact subset of $\mathcal{D}_{r}.$   We  shall use  the same  symbol  $s$
and  set
 $$s := \sum_{j=1}^{\vert N\vert +2} s_{j}d\Theta_{j}
 .$$
 In  the sequel, we  shall use  simultaneously  the following notations  for  $\Theta\in \C^{\vert N\vert +2}:$
 $$  \Theta \equiv \left ( \Theta_1,\ldots,\Theta_{\vert N\vert +2}  \right )
 \equiv  (\xi,\zeta,\eta)\in \C^l\times\C^{n+1}\times\C^{m+1}. $$
We next  set
$$ \Phi :=   \left( \sum_{j=1}^{n+1}
(\zeta_{j} +z_j) d\zeta_{j} \right)\wedge \left(\sum_{k=1}^{m+1}
(\eta_{k} +w_k) d\eta_{k}\right). $$
  For  every $\epsilon  >0,$
consider the differential form of  bidegree $(\vert N\vert+2,\vert N\vert+1)$
\begin{eqnarray}
K_{s}^{\epsilon} := \frac{ s\wedge (\overline{\partial}s)^{\vert N\vert-1} \wedge
( \overline{\partial} Q_{ \epsilon})^{2} }  {\lbrack s(\Theta,Z)\bullet(\Theta -Z)
\rbrack^{\vert N\vert}},
\end{eqnarray}
where   $ Q_{ \epsilon}$ is the differential form of  bidegree
 $(1,0)$   given by
    \begin{eqnarray} Q_{ \epsilon}  :=\frac{ \overline{
\zeta\bullet \zeta}  \left( \sum_{j=1}^{n+1}
(\zeta_{j} +z_j) d\zeta_{j} \right)+ \overline{ \eta\bullet \eta}  \left( \sum_{k=1}^{m+1}
(\eta_{k} +w_k) d\eta_{k} \right) }{\vert  \zeta \bullet \zeta\vert^{2} +\vert  \eta \bullet \eta\vert^{2}  + \epsilon  }  .
\end{eqnarray}
  Denote by
  $d\Theta$ the canonical  holomorphic form of  $\C^{\vert N\vert+2}$ given by
$$ d\Theta_1 \wedge \ldots \wedge d\Theta_{\vert N\vert+2    } \equiv
d\xi_1\wedge\ldots\wedge d\xi_l
\wedge d\zeta_1\wedge\ldots\wedge d\zeta_{n+1}\wedge d\eta_1\wedge\ldots\wedge d\eta_{m+1}.$$
\begin{lem}  Suppose that  $ 0<r<1.$

 1) If $u\in   \mathcal{C}^1(  \overline{\mathcal{D}}_{r}  )$
and $Z\in \overline{\M}_{r}    ,$   then   \begin{eqnarray*}
u(Z) = C(N)\lim\limits_{\epsilon\to 0} \left(
\int\limits_{\partial \mathcal{D}_{r} }  u K_{s}^{\epsilon}   -
\int\limits_{\mathcal{D}_{r} }
\overline{\partial} u\wedge  K_{s}^{\epsilon}\right ).
\end{eqnarray*}

2)
If $u \in \mathcal{C}(\mathcal{D}_{r}),$   then
\begin{eqnarray*}\lim\limits_{\epsilon \to 0}
\int\limits_{ \mathcal{D}_{r}    }
\frac{\epsilon  \vert \zeta \vert^{2}  \vert \eta \vert^{2}       u(\Theta)
 }{
\left(\vert \zeta \bullet \zeta \vert^{2}+ \vert \eta \bullet \eta\vert^{2} +     \epsilon\right)^3}
 d\overline{\Theta}\wedge d\Theta
=C(N)\int\limits_{ \M_{r}       }u(\Theta)
dV(\Theta).  \end{eqnarray*}

3)
If $u \in \mathcal{C}(\partial \B_N )$  and
 $\omega$  is the canonical volume form of $\partial\B_N,$
then
\begin{eqnarray*} \lim\limits_{\epsilon \to 0}
\int\limits_{\partial   \B_N      }
 \frac{\epsilon  \vert \zeta \vert^{2}  \vert \eta  \vert^{2}     u(\Theta) }{
 \left(\vert \zeta\bullet \zeta \vert^{2}+ \vert \eta\bullet \eta \vert^{2} +     \epsilon\right)^3
 } \omega(\Theta)
  =C(N)
\int\limits_{\partial\M_N}u(\Theta)d\sigma(\Theta). \end{eqnarray*}
\end{lem}
\begin{proof}
Part 1)   follows from formulas  (23) and (26)  in the proof of Theorem 1 in  \cite{BB}. Also,
part 2) is  an immediate consequence of  identity (25)    in  \cite{BB}.

  To prove part 3), we  may assume  without  loss of  generality  that
the support of  $u$ is contained in a sufficiently small open  neighborhood $\mathcal{U}\subset \C^{\vert N\vert +2}$
of  a point  $\Theta_0\in \partial\M_N.$   Using  local  coordinates   and Lelong theory  \cite{LG},
we  see that  there exists a  smooth $(2\vert N\vert-1)$-volume form  $d\mu$ defined  on
$\mathcal{U} \cap \partial\M_N$   such that
\begin{equation}  \lim\limits_{\epsilon \to 0}
\int\limits_{\partial   \B_N \cap  \mathcal{U}     }\frac{\epsilon  \vert \zeta \vert^{2}  \vert \eta  \vert^{2} u(\Theta)   }{
 \left(\vert \zeta\bullet \zeta \vert^{2}+ \vert \eta\bullet \eta \vert^{2} +     \epsilon\right)^3
 }\omega(\Theta) =
C
\int\limits_{\partial\M_N \cap \mathcal{U}}u(\Theta)d\mu(\Theta),
\end{equation}
for all $ u\in \mathcal{C}_{0}(\mathcal{U}).$
 Therefore, part 3)  is equivalent to the  identity    $d\mu=Cd\sigma.$

Let $\psi$ be a  function  of class  $\mathcal{C}^{\infty}_{0} (\lbrack  0,1\rbrack)$  supported in  $[\frac{1}{2},1]$ such that
$\int\limits_{0}^{1}  \rho^{2\vert N\vert -1}  \psi(\rho)d\rho  =1.$  Consider  the $\mathcal{C}^{\infty}_{0}$ extension of $u$ given by
$$  u(\rho Z):= \psi(\rho) u(Z),\quad \text{for}\ 0\leq \rho \leq 1\ \text{and} \ Z \in\partial\B_N\cap\mathcal{U}.$$
On the one  hand, using (2.4), we have that
\begin{eqnarray*}
\lim\limits_{\epsilon \to 0}
\int\limits_{\B_N
 }
\frac{\epsilon  \vert \zeta \vert^{2}  \vert \eta \vert^{2}       u(\Theta)
 }{
\left(\vert \zeta \bullet \zeta \vert^{2}+ \vert \eta \bullet \eta\vert^{2} +     \epsilon\right)^3}
 d\overline{\Theta}\wedge d\Theta=\lim\limits_{\epsilon \to 0}\int\limits_{0}^{1} \int\limits_{\partial   \B_N \cap\mathcal{U}     }
 \frac{\epsilon \rho^{2\vert N\vert +7}  \vert \zeta \vert^{2}  \vert \eta  \vert^{2}     u(\rho\Theta) }{
 \left(\rho^4\vert \zeta\bullet \zeta \vert^{2}+ \rho^4\vert \eta\bullet \eta \vert^{2} +     \epsilon\right)^3
 } \omega(\Theta)d\rho \\
 =  \lim\limits_{\epsilon \to 0}\int\limits_{0}^{1}  \rho^{2\vert N\vert -1}  \psi(\rho) d\rho \cdot \int\limits_{\partial   \B_N      \cap\mathcal{U}}
 \frac{\epsilon  \vert \zeta \vert^{2}  \vert \eta  \vert^{2}     u(\Theta) }{
 \left(\vert \zeta\bullet \zeta \vert^{2}+ \vert \eta\bullet \eta \vert^{2} +     \epsilon\right)^3
 } \omega(\Theta)= \int\limits_{\mathcal{U} \cap \partial\M_N  } ud\mu.
 \end{eqnarray*}
 On the other hand, by part 2) and Lemma 2.3,  we see that
 \begin{eqnarray*}
& & \lim\limits_{\epsilon \to 0}
\int\limits_{\B_N }
\frac{\epsilon  \vert \zeta \vert^{2}  \vert \eta \vert^{2}       u(\Theta)
 }{
\left(\vert \zeta \bullet \zeta \vert^{2}+ \vert \eta \bullet \eta\vert^{2} +     \epsilon\right)^3}
 d\overline{\Theta}\wedge d\Theta
 =C(N)\int\limits_{ \M_N  } udV\\
 &=&  C(N)\int\limits_{0}^{1}  \rho^{2\vert N\vert -1}  \psi(\rho) \int\limits_{ \mathcal{U} \cap \partial\M_N  } ud\sigma  d\rho
 = C(N) \int\limits_{\mathcal{U} \cap \partial\M_N  } ud\sigma,
 \end{eqnarray*}
  Thus  $d\mu= C(N)d\sigma$ and  thereby completes the proof.
  \end{proof}
 Next, set
\begin{eqnarray}
K_s:=\frac{ s\wedge (\overline{\partial}s)^{\vert N\vert -1} \wedge
  \Phi\wedge \overline{\partial (\zeta\bullet\zeta)}\wedge\overline{\partial (\eta\bullet\eta)}  }  {
 \vert \zeta\vert^2  \vert \eta\vert^2
  \lbrack s(\Theta,Z)\bullet(\Theta -Z)
\rbrack^{\vert N\vert}   }.
\end{eqnarray}
 In  view of  (2.2), (2.3) and the equality which  precedes
  Lemma 4 in \cite{BB},   we  see that
$K_s$ satisfies the   identity
 \begin{eqnarray}
K_{s}^{\epsilon} = \frac{\epsilon  \vert \zeta \vert^{2}  \vert \eta \vert^{2}
 }{
\left(\vert \zeta \bullet \zeta \vert^{2}+ \vert \eta \bullet \eta\vert^{2} +     \epsilon\right)^3} K_s.
\end{eqnarray}
For  every $1 \leq k \leq  \vert N\vert +2,$ denote by  $ \omega_{k}(\overline{\Theta} )     $  the
 $(0,\vert N\vert +1)$-form
 $$ (-1)^{k-1}
d\overline{\Theta}_{1}\wedge  \ldots  \wedge \widehat{d\overline{\Theta}_{k}}
\wedge \ldots \wedge d\overline{\Theta}_{\vert N\vert +2}.$$
 We can  write  $K_s$  in the form
\begin{eqnarray}
K_s=    \sum\limits_{k=1}^{\vert N\vert +2} h_k(\Theta,Z)
 \omega_{k}(\overline{\Theta})   \wedge    d\Theta, \end{eqnarray}
where
 $h_k$  are the  component functions of $K_s$   with  respect  to the forms
$\omega_1(\overline{\Theta})\wedge d\Theta,\ldots,\omega_{\vert N\vert +2}(\overline{\Theta}) \wedge d\Theta.$

Let $\overline{\M}_N$ be the closure of $\M_N$ in $\overline{\B}_N$
and  denote by  $\mathcal{C}^k(\overline{\M}_N),$ $k\in\N,$ the space of all $\mathcal{C}^k$ functions defined
in a  neighborhood of $\overline{\M}_N$  in $\overline{\B}_N.$
If $f:= \sum\limits_{j=1}^{\vert N\vert +2}  f_j d\overline{\Theta}_j$ is  a $(0,1)$-form
with coefficients in $\mathcal{C}(\overline{\M}_N),$ let $f|_{\M_N}$
denote the pull-back of $f$   under the canonical   injection of  $\M_N$  in this  neighborhood.
  Set
   \begin{eqnarray}
\Vert f\Vert_{\M_{N,\infty}}:=\underset{\Theta\in\M_N}{\sup}\sum\limits_{j=1}^{\vert N\vert +2}
\vert f_j (\Theta)\vert.
\end{eqnarray}

Let  $\overline{\partial}_{\M_N}$ be the $\overline{\partial}$-operator on
$\M_N .$ We end  this  section  by the following
\begin{prop} Consider a section   $s$  satisfying  (2.1),
  a function $u \in \mathcal{C}^{1}(\overline{\M}_N)$  and  a
$(0,1)$-form   $f:= \sum\limits_{k=1}^{\vert N\vert +2}   f_k
d\overline{\Theta}_k$   with  coefficients  in $\mathcal{C}(\overline{\M}_N)$      that  satisfy $\overline{\partial}_{\M_N}u=
f|_{\M_N}$   on  $  \M_N . $   Let $h_k$  be the functions
defined in (2.7). Then for $Z\in  \M_N ,  $
$$
u(Z) = C  \int\limits_{\partial\M_N}  u(\Theta) \left
(\sum\limits_{k=1}^{\vert N\vert +2}\Theta_k h_k(\Theta,Z)  \right ) d\sigma(\Theta) +
C\int\limits_{  \M_N   }\left
(\sum\limits_{k=1}^{\vert N\vert +2}  f_k(\Theta)h_k(\Theta,Z)  \right )dV(\Theta).
$$
 \end{prop}

\begin{proof} For every $r\in\rbrack 0,1\lbrack$ such that $Z\in \M_{r},$   consider a $\mathcal{C}^1$ extension of $u|_{\M_N}$  (which is  also  denoted by
$u$) on $\overline{\mathcal{D}}_{r}$ that  satisfies  $\overline{\partial}u=
f$  on  $ \M_{r} .$  Suppose  without loss of  generality that
$f= \overline{\partial}u$ on $\overline{\mathcal{D}}_{r}    .$
Parts 1) and 2) of   Lemma 2.4, combined with   (2.6) and  (2.7), imply  that
$  u(z)= CI^1_{r}  + CI^2_{r} ,
$
where
\begin{eqnarray*}
 I^1_{r}   &:= &
\int\limits_{\M_{r} }
\left (\sum\limits_{k=1}^{\vert N\vert +2            } f_k(\Theta)h_k(\Theta,Z)  \right )
dV(\Theta),\\
I^2_{r}    &:= &
\lim\limits_{\epsilon\to 0}
\int\limits_{\partial \mathcal{D}_{r} }  \frac{\epsilon \vert \zeta\vert^2 \vert \eta\vert^2 u(\Theta) }{
\left(\vert \zeta\bullet\zeta          \vert^{2}  + \vert \eta\bullet\eta          \vert^{2}
+\epsilon\right)^3}
\left
(\sum\limits_{k=1}^{\vert N\vert +2} h_k(\Theta,Z) \omega_{k}(\overline{\Theta})\wedge   d\Theta
  \right ). \\
\end{eqnarray*}
The proof is a consequence of   the following
two  equalities
  \begin{eqnarray}
 \lim\limits_{r\to 0       }I^1_{r}  & =& \int\limits_{\M_N }
\left (\sum\limits_{k=1}^{\vert N\vert +2} f_k(\Theta)h_k(\Theta,Z)  \right )
dV(\Theta), \\
  \lim\limits_{ r\to 0       }I^2_{r} & = & \int\limits_{\partial\M_N}  u(\Theta) \left
(\sum\limits_{k=1}^{\vert N\vert +2}\Theta_k h_k(\Theta,Z)  \right ) d\sigma(\Theta).
\end{eqnarray}

In order to prove these, fix a point $Z\in\M_N.$
By  (2.1),
(2.5) and (2.7),  there is  a constant $C$  such that
\begin{eqnarray}
h_k(\Theta,Z)  \leq  \frac{C}{ \vert \zeta\vert^2 \vert \eta\vert^2},\qquad \text{for all}\
\Theta \in\M_N\setminus \M_r,\ \text{with}\ 0<r<<1 .
\end{eqnarray}
We deduce easily from  (2.11) and the hypothesis  $\Vert f\Vert_{\M_{N,\infty}} <\infty$
that
$$ \lim\limits_{r\to 0} \int\limits_{\M_N  \setminus \M_{r} }
\left \vert\sum\limits_{k=1}^{\vert N\vert +2} f_k(\Theta)h_k(\Theta,Z)  \right \vert
dV(\Theta)  \leq  \lim\limits_{r\to 0}\int\limits_{\M_N  \setminus \M_{r} }  \frac{CdV(\Theta)}{\vert \zeta\vert^2 \vert \eta\vert^2 }=0,
$$
where   the equality follows from   Corollary 2.2. This proves (2.9).

Next, we  prove    (2.10).   Appealing to   Corollary 2.2,  Lemma 2.3,
 (2.11) and the fact that  the function $u$ is bounded,  we  see  that
\begin{eqnarray*}
 \lim\limits_{r\to 0} \int\limits_{\partial \M_N  \setminus \partial\M_{r} }
\left \vert\sum\limits_{k=1}^{\vert N\vert +2} \Theta_k h_k(\Theta,Z)  \right \vert \vert u(\Theta) \vert
d\sigma(\Theta)  \leq  \lim\limits_{r\to 0}\int\limits_{\partial \M_N  \setminus \partial\M_{r} }  \frac{Cd\sigma(\Theta)}{\vert \zeta\vert^2 \vert \eta\vert^2 }=0.
\end{eqnarray*}
This, combined with   part 3)  of  Lemma 2.4,  implies  that
\begin{equation}\begin{split}
&\lim\limits_{\epsilon\to 0}
\int\limits_{\partial \B_N \setminus  \partial \mathcal{D}_{r} }  \frac{\epsilon \vert \zeta\vert^2 \vert \eta\vert^2 u(\Theta) }{
\left(\vert \zeta\bullet\zeta          \vert^{2}  + \vert \eta\bullet\eta          \vert^{2}
+\epsilon\right)^3}
\left
(\sum\limits_{k=1}^{\vert N\vert +2} h_k(\Theta,Z) \omega_{k}(\overline{\Theta})\wedge   d\Theta
  \right ) \\
&= \int\limits_{\partial\M_N \setminus \partial\M_{r}}  u(\Theta) \left
(\sum\limits_{k=1}^{\vert N\vert +2}\Theta_k h_k(\Theta,Z)  \right ) d\sigma(\Theta)\longrightarrow 0,\ \ \text{as}\  r\to 0.
\end{split}
\end{equation}
from which it follows that (2.10) is a consequent of
\begin{equation}\lim\limits_{r\to 0}
\lim\limits_{\epsilon\to 0}
\int\limits_{\partial \mathcal{D}_r\setminus \partial \B_N }  \frac{\epsilon \vert \zeta\vert^2 \vert \eta\vert^2 u(\Theta) }{
\left(\vert \zeta\bullet\zeta          \vert^{2}  + \vert \eta\bullet\eta          \vert^{2}
+\epsilon\right)^3}
\left
(\sum\limits_{k=1}^{\vert N\vert +2} h_k(\Theta,Z) \omega_{k}(\overline{\Theta})\wedge   d\Theta
  \right ) =0.
  \end{equation}

  Next, we    prove  equality  (2.13). We first  make use of  the following remark  related to
    homogeneity properties  of  certain   differential  forms.
 Indeed, let   $\alpha,\beta  > 0$ and     write the complex manifold   $\M_N$
   as a  complete intersection  of $\B_N$    and  the  two varieties  given  by  the  equations
  $\alpha^2 \zeta\bullet\zeta=0$ and  $\beta^2 \eta\bullet \eta=0.$    Applying    Berndtsson's  formulas
  to  these  two  equations and  observing that   (2.13)   corresponds  to the  particular case  $\alpha=\beta=1,$
  then  we obtain
 \begin{equation}
 \begin{split}
&\lim\limits_{\epsilon\to 0}
\int\limits_{ \partial \mathcal{D}_r\setminus \partial \B_N              }  \frac{\epsilon\alpha^4\beta^4 \vert \zeta\vert^2 \vert \eta\vert^2 u(\Theta) }{
\left(\alpha^4\vert \zeta\bullet\zeta          \vert^{2}  +\beta^4 \vert \eta\bullet\eta          \vert^{2}
+\epsilon\right)^3}
\left
(\sum\limits_{k=1}^{\vert N\vert +2} h_k(\Theta,Z) \omega_{k}(\overline{\Theta})\wedge   d\Theta
  \right )\\
 &=\lim\limits_{\epsilon\to 0}
\int\limits_{\partial \mathcal{D}_r\setminus \partial \B_N }  \frac{\epsilon \vert \zeta\vert^2 \vert \eta\vert^2 u(\Theta) }{
\left(\vert \zeta\bullet\zeta          \vert^{2}  + \vert \eta\bullet\eta          \vert^{2}
+\epsilon\right)^3}
\left
(\sum\limits_{k=1}^{\vert N\vert +2} h_k(\Theta,Z) \omega_{k}(\overline{\Theta})\wedge   d\Theta
  \right ) ,
  \end{split}
  \end{equation}
  for all  $0<r<1.$

  We write  $\partial \mathcal{D}_r \setminus \partial \B_N        $  as a  union of the   two smooth  manifolds
\begin{eqnarray*}
 M_1^{r}         &:=& \{  Z\in \B_N:
  \vert z\vert = r,  \vert w\vert \geq r \};\\
 M_2^{r}         &:=& \{  Z\in \B_N:
  \vert z\vert \geq r,  \vert w\vert = r \}.
\end{eqnarray*}

  Let  $d\sigma_{rj} $ be the canonical  volume form on the manifold    $M_j^r,$ $j=1,2.$
  Applying equality  (3) in
 Proposition 16.4.4  of  Rudin  \cite{Ru} yields that  on  $ M_j^{r}, $
\begin{equation}
 \omega_{k}(\overline{\Theta})\wedge   d\Theta= C(N,j,k)d\sigma_{rj}.
\end{equation}
  Choosing a  function $u$  and a  section $s$ appropriately and     applying
  Lelong theory as in the proof of  (2.4),   it follows  from (2.14) and (2.15)
  that   on   $M_j^r,$  $0<r<1,\ j\in\{1,2\},$ we have
 \begin{equation}\begin{split}
 & \lim\limits_{\epsilon\to 0}\frac{\epsilon\alpha^4\beta^4 \vert \zeta\vert^2 \vert \eta\vert^2  }{
\left(\alpha^4\vert \zeta\bullet\zeta          \vert^{2}  +\beta^4 \vert \eta\bullet\eta          \vert^{2}
+\epsilon\right)^3}d\sigma_{rj}(\Theta)\\
&=\lim\limits_{\epsilon\to 0}\frac{\epsilon \vert \zeta\vert^2 \vert \eta\vert^2  }{
\left(\vert \zeta\bullet\zeta          \vert^{2}  + \vert \eta\bullet\eta          \vert^{2}
+\epsilon\right)^3}d\sigma_{rj}(\Theta)=d\mu_{rj}(\Theta),
\end{split} \end{equation}
in the distribution  sense,
  where $d\mu_{rj}$  is a $\mathcal{C}^{\infty}$ differential  form of maximal  degree  on the manifold   $M_j^r\cap \H_N.$  In view of (2.11)  and (2.16), equality (2.13) will  follow from the following  equalities
\begin{equation}
\lim\limits_{r\to 0}
\int\limits_{ M_j^{r}\cap \H_N        }  \frac{d\mu_{rj}(\Theta)}   { \vert \zeta\vert^2 \vert \eta\vert^2 } =0,\qquad j=1,2.
  \end{equation}
 We  prove   (2.17) for   $j=1$ which suffices to complete the proof. To do so, consider,
  for  every $\alpha,\beta  >0,$ the mapping  $F_{\alpha,\beta}$  given by~:
  $$  F_{\alpha,\beta}(x,z,w):=(x,\alpha z,\beta w),\qquad  \text{for}\ Z\equiv (x,z,w)\in\H_N.$$
 We  remark immediately that we have   the following  property of  homogeneity~:
  \begin{equation*}
  F^{*}_{r,s}(d\sigma_{r1} )(\Theta)= C(N)r^{2n-1} s^{2m} d\sigma_n(\zeta)\wedge dV_l(\xi)\wedge dV_m(\eta),
  \end{equation*}
  for $0<r,s\leq \frac{1}{2}$
  and   $\Theta\equiv(\xi,\zeta,\eta) \in  \C^l\times\partial\M_n\times\M_m.$
  This, combined  with  equality (2.16), implies that
  \begin{equation}
   F^{*}_{\alpha,\beta}(d\mu_{r1} )= C(N,r)\alpha^{2n-1} \beta^{2m} d\mu_{\frac{r}{\alpha},1 }\qquad
   \text{on}\ M^{\frac{r}{\alpha}     }_1.
   \end{equation}
 Take $r_0:= \frac{1}{2}.$   Since  the  differential form  $d\mu_{r_0,1}$ is in $\mathcal{C}^{\infty}(M^{r_0}_1),$
   we see that
   \begin{equation}
     \int\limits_{ \vert\zeta\vert =r_0,\frac{r_0}{2} <\vert \eta\vert <r_0}
     \frac{   d\mu_{r_0,1}(\Theta)}{  \vert \zeta\vert^2\vert \eta\vert^2}   < \infty.
     \end{equation}
Using  (2.18)  and (2.19), it is easy to  show that
  $$
  \int\limits_{ M_1^{r}\cap \H_N        }  \frac{d\mu_{r1}(\Theta)}   { \vert \zeta\vert^2 \vert \eta\vert^2 }
  \leq   C \left (\frac{r}{r_0} \right )^{2n-3}\int\limits_{ \vert\zeta\vert =r_0,\frac{r_0}{2} <\vert \eta\vert <r_0}
     \frac{   d\mu_{r_0,1}(\Theta)}{  \vert \zeta\vert^2\vert \eta\vert^2} \to 0,\ \ \text{as}\ r\to 0.
  $$
 This  implies (2.17) and thus completes the proof.
      \end{proof}
\section{Integral formulas on the manifold  $\M_N.$}

In  this section  we establish
     integral  formulas of  Martinelli-Bochner type (Theorem  3.1) and those  of Cauchy type
(Theorems  3.3 and 3.6).
 These  formulas  will allow us to construct the   $\overline\partial$-solving operators.
\begin{thm}
Suppose that  $u \in \mathcal{C}^{1}(\overline{\M}_N)$  and
  $f:= \sum\limits_{k=1}^{\vert N\vert +2}   f_k d\overline{\Theta}_k$ is a $(0,1)$-form
with coefficients in  $\mathcal{C}(\overline{\M}_N) $   such that
$\overline{\partial}_{\M_N}u=  f|_{\M_N}.$
 Then  for  every  $Z\in \M_N,$
  \begin{eqnarray*} u(Z)&= &  \left \{
\int\limits_{\partial\M_N}   \frac{A(\Theta,Z)}{ \vert  Z-\Theta \vert^{2\vert N\vert}}  u(\Theta) \frac{  d\sigma(\Theta)
}{  \vert \zeta\vert^2   \vert \eta\vert^2} \right .\\
 &+& \left .\int\limits_{\M_N}  \frac{1}{ \vert  Z-\Theta \vert^{2\vert N\vert}}  \left ( \sum\limits_{k=1}^{\vert N\vert +2} B_k(\Theta,Z)f_k(\Theta)\right )
 \frac{ dV(\Theta)}{  \vert \zeta\vert^2   \vert \eta\vert^2}  \right \} ,
 \end{eqnarray*}
where  \begin{eqnarray*}
A(\Theta,Z)&:=&C (\vert \xi \vert^2  -\overline{x}\bullet  \xi ) (\vert \zeta \vert^2  +z\bullet \overline{ \zeta} )
(\vert \eta \vert^2  +w\bullet \overline{ \eta} )\\ &  +&
C \left (
-\vert z\bullet\zeta \vert^{2}+\vert  z\bullet\overline{\zeta} \vert^{2}
  - \vert\zeta\vert^{2} (\vert\zeta\vert^{2}+
z\bullet\overline{\zeta}-\overline{z}\bullet\zeta ) \right)(\vert \eta \vert^2  +w\bullet \overline{ \eta} )\\
& +& C \left (
-\vert w\bullet\eta \vert^{2}+\vert  w\bullet\overline{\eta} \vert^{2}
  - \vert\eta\vert^{2} (\vert\eta\vert^{2}+
w\bullet\overline{\eta}-\overline{w}\bullet\eta ) \right)(\vert \zeta \vert^2  +z\bullet \overline{ \zeta} ),
\end{eqnarray*}
and $B_k$  are   polynomials  given  by   the following formulas~:
\\
(i) if $1\leq k \leq  l,$  then
$$ B_k(\Theta,Z):=  C (\overline{\xi_k} -  \overline{x_k}) (\vert \zeta \vert^2  +z\bullet \overline{ \zeta} )
(\vert \eta \vert^2  +w\bullet \overline{ \eta} );$$
\\
(ii) if $l < k\leq l+n+1$ and $j=k-l,$ then
$$ B_k(\Theta,Z):= C\left ( (\overline{z}_{j}
-\overline{\zeta}_{j})(z\bullet\overline{\zeta}+\vert\zeta\vert^{2}) -(z_{j}
+\zeta_{j})\overline{z\bullet(\zeta  -z)} \right )(\vert \eta \vert^2  +w\bullet \overline{ \eta} );$$
 \\
(iii) if $l+n+1  < k< l+n+m+2$  and $i=k-l-n-1,$   then
$$ B_k(\Theta,Z):= C\left ( (\overline{w}_{i}
-\overline{\eta}_{i})(w\bullet\overline{\eta}+\vert\eta\vert^{2}) -(w_{i}
+\eta_{i})\overline{w\bullet(\eta  -w)} \right )(\vert \zeta \vert^2  +z\bullet \overline{ \zeta} ).$$
 \end{thm}
\begin{proof}   Consider the  Martinelli-Bochner section
  $s_{b}(Z,\Theta):=
\overline{\Theta}-\overline{Z}.$
In order to prove the   theorem, we  apply Proposition
2.5   to the  section $s_b.$  Using formulas (2.5), (2.7)  and   arguing as  in the proof of
 Theorem  2.4  of \cite{VY1}, we  compute  explicitly the functions $h_k$
associated to $s_b$  and  obtain the desired formula.
\end{proof}
\begin{rem}{\rm  If $u\in \mathcal{C}^{1}(\M_N)$ is bounded, then Proposition 2.5  and Theoreme 3.1
hold for the dilated functions  $u_r(Z):=
u(rZ),\ 0 < r< 1.$   This  shows that   Theorem 3.1 remains true
if we only  assume  that  $u\in \mathcal{C}^{1}(\M_N)$   is  bounded   and
$$\lim\limits_{r\to 1^{-}}
\int\limits_{\partial\M_N} \vert u(\Theta)- u(r\Theta) \vert d\sigma(\Theta)=0.$$}
\end{rem}
Following   Charpentier \cite{Ch} let
$$  s_{0}(\Theta,Z) := \overline{\Theta}(1- \Theta\bullet \overline{Z}) -
\overline{Z}( 1-\vert \Theta \vert^{2}), \qquad \text{and}\qquad  D(\Theta,Z):=
s_{0}(\Theta,Z)\bullet(\Theta - Z).$$
In what follows,   $ \text{grad}_{Z}\ f$  denotes the gradient of  a differentiable function $f$ at a point $Z.$
\begin{thm}
There  exist  polynomials $R(\Theta,Z)$ and $P_{k}(\Theta,Z),\ Q_{k}(\Theta,Z)$ for $ 1\leq k\leq \vert N\vert  +2,$
that   satisfy  the  following properties:\\
(i) $
R(\Theta,Z)= \left( C\vert  \xi\vert^2  +  C\vert  \zeta\vert^2 +C\vert  \eta\vert^2 \right )
      (\vert \zeta \vert^2  +z\bullet \overline{ \zeta} )
(\vert \eta \vert^2  +w\bullet \overline{ \eta} ).
$ \\
(ii)   For  every $Z,\Theta\in \B_N,$  and for every $1\leq k\leq \vert N\vert +2,$
  \begin{eqnarray*}P_{k}(\Theta,Z)&=& O\left (\vert \Theta -Z\vert (\vert \zeta \vert^2  +\vert z\vert \vert \zeta\vert               )
(\vert \eta \vert^2  + \vert w\vert \vert \eta\vert                ) \right ),\\
Q_{k}(\Theta,Z)&=&  O\left (\vert \Theta -Z\vert (\vert \zeta \vert^2  + \vert z\vert \vert \zeta\vert )
(\vert \eta \vert^2  +\vert w\vert \vert \eta \vert ) \right ),\\
\left \vert  \text{grad}_{Z}\ P_{k}(\Theta,Z)\right\vert &=& O\left ( (\vert \zeta \vert^2  +\vert \zeta\vert \vert z\vert)
(\vert \eta \vert^2  +\vert\eta\vert \vert w\vert ) \right), \\
\left \vert  \text{grad}_{Z}\ Q_{k}(\Theta,Z)\right\vert &=&O\left ( (\vert \zeta \vert^2  +\vert \zeta\vert \vert z\vert)
(\vert \eta \vert^2  +\vert\eta\vert \vert w\vert ) \right) .
\end{eqnarray*} (iii)
  Given a  function   $u \in
\mathcal{C}^{1}(\overline{\M}_N)$   and  a $(0,1)$-form  $f:=
\sum\limits_{k=1}^{\vert N\vert +2}   f_k d\overline{\Theta}_k \in\mathcal{C}(\overline{\M}_N)
    $    that  satisfy $ \overline{\partial}_{\M_N}u=  f|_{\M_N},$     then  for  every $Z\in \M_N,$
\begin{eqnarray*}& &\quad  u(Z)  = \int\limits_{\partial\M_N} \frac{R(\Theta,Z)}{  \left (1-Z\bullet \overline{\Theta}\right )^{\vert N\vert}} u(\Theta)
\frac{ d\sigma(\Theta)}{  \vert \zeta\vert^2   \vert \eta\vert^2} + \\
 & &\int \limits_{\M_N} \sum_{k=1}^{\vert N\vert +2}
\frac{(1- \Theta\bullet\overline{Z})^{\vert N\vert -2} }{ D(\Theta,Z)^{\vert N\vert}}\lbrack (
1-\Theta\bullet \overline{Z}) P_{k}(\Theta,Z) +( 1-\vert \Theta
\vert^{2})Q_{k}(\Theta,Z)\rbrack
   f_k(\Theta)\frac{ dV(\Theta)}{  \vert \zeta\vert^2   \vert \eta\vert^2}.
\end{eqnarray*}
\end{thm}
\begin{proof}
From the proof of   Proposition 2.5  we  may  assume without loss of  generality that
 there  is  a  $\mathcal{C}^1$ extension of  $u|_{\M_N}$, denoted again  by
$u$,   such that  $  \overline{\partial}u=  f$ on $\B_N.$
Let    $K_0$ be the kernel  associated   to the  section  $  s_{0}$  by formula (2.5).  By  virtue of   (2.6), when  we
integrate    $uK^{\epsilon}_{0}$ over  $\partial\B_N$, all  terms    which  contain
$\overline\partial \vert \Theta \vert^{2}$   vanish. In addition    we  have
$1-\vert \Theta \vert^2=0$ and $D(\Theta,Z)=\left\vert
1-Z\bullet\overline{\Theta}\right\vert^{2}$   so that
  \begin{eqnarray*}  & & \qquad  \lim\limits_{\epsilon \to 0}
\int\limits_{\partial\B_N  } u K_{0}^{\epsilon} = \lim\limits_{\epsilon \to 0}
\int\limits_{ \partial\B_N   }u(\Theta)
\frac{\epsilon\vert \zeta\vert^2 \vert \eta\vert^2
 }{
\left(\vert \zeta \bullet \zeta \vert^{2}+ \vert \eta \bullet \eta\vert^{2} +     \epsilon\right)^3} \\
 &\cdot &
\frac{1}{\vert \zeta\vert^2 \vert \eta\vert^2(1-Z\bullet\overline{\Theta})^{\vert N\vert} }\left\lbrace
\left ( \sum\limits_{k=1}^{\vert N\vert +2} \overline{\Theta}_k d\Theta_k \right )
\left ( \sum_{k=1}^{\vert N\vert +2}     d\overline{\Theta}_{k}
\wedge d\Theta_{k}  \right )^{\vert N\vert -1} \wedge \Psi \wedge\overline{ \partial (\zeta\bullet\zeta)}\wedge
 \overline{ \partial (\eta\bullet\eta)}\right\rbrace.
\end{eqnarray*}
Rewriting  the differential form in braces in the form
$
\sum\limits_{k=1}^{\vert N\vert +2} h_k(\Theta,Z) \omega_k( \overline{\Theta})\wedge d\Theta$
and  applying part 3) of Lemma 2.4, we  obtain
$$ \lim\limits_{\epsilon \to 0}
\int\limits_{\partial\B_N  } u K_{0}^{\epsilon} =
\int\limits_{\partial\M_N} \frac{ \Theta_k h_k(\Theta,Z)}
{\vert \zeta\vert^2 \vert \eta\vert^2(1-Z\bullet\overline{\Theta})^{\vert N\vert} }                       d\sigma(\Theta).$$
A straightforward  calculation of the   functions $h_k(\Theta,Z)$  shows that
\begin{equation}
R(\Theta,Z):= \sum\limits_{k=1}^{\vert N\vert +2} \Theta_k h_k(\Theta,Z)   \end{equation}  satisfies  assertion  (i)   of the theorem.

Write the kernel  $K_0$  in the  form (2.7) as
$
K_0 =\sum\limits_{k=1}^{\vert N\vert +2} h_k(\Theta,Z) \omega_k(\overline{\Theta}) \wedge d\Theta.
$
Then  we have
\begin{eqnarray}
  I:= \overline{\partial} u\wedge K_0 =\sum\limits_{k=1}^{\vert N\vert +2}f_k(\Theta) h_k(\Theta,Z)
  d\Theta   \wedge d\overline{\Theta}.
\end{eqnarray}
To finish  the proof of the theorem, it suffices  to  prove  the following lemma~:
\begin{lem}
The functions $h_k$  in the formula (3.2) can  be  rewritten  in the  form
\begin{equation}
 h_k(\Theta,Z)=
\frac{(1- \Theta\bullet\overline{Z})^{\vert N\vert -2} }{ \vert \zeta
\vert^{2}  \vert \eta
\vert^{2}     D(\Theta,Z)^{\vert N\vert}}\lbrack (
1-\Theta\bullet \overline{Z}) P_{k}(\Theta,Z) +( 1-\vert \Theta
\vert^{2})Q_{k}(\Theta,Z)\rbrack ,
\end{equation}
where  $P_k$ and $Q_k$ are  some  polynomials  that satisfy  assertion (ii)   of the theorem.
\end{lem}

\noindent{\it End  of the proof  of Theorem  3.3.}
Suppose that  the  lemma  above  is  proved.  Applying  Proposition 2.5 and using  (3.1)--(3.3),
the theorem follows.
\end{proof}
\noindent{\it  Proof of Lemma 3.4.}
By virtue of (2.5) and (3.2), we  can  write
$
I=  I_1 + I_2,$
where
\begin{eqnarray*}
I_1 &:= &   \left ( \sum\limits_{k=1}^{\vert N\vert +2}f_k
   d\overline{\Theta}_k \right )
\wedge \left \lbrace \frac{(1- \Theta\bullet\overline{Z})^{\vert N\vert -1} }{\vert \zeta
\vert^{2}  \vert \eta
\vert^{2}
D(\Theta,Z)^{\vert N\vert}} \sum_{j=1}^{\vert N\vert +2} \left \lbrack\overline{\Theta}_{j} (1-
\Theta\bullet \overline{Z}) - \overline{Z}_{j}( 1-\vert \Theta \vert^{2})
\right\rbrack \right .\\ & & \left .d\Theta_{j} \wedge    \Psi  \wedge \overline{\partial(\zeta\bullet\zeta)}
  \wedge \overline{\partial(\eta\bullet\eta)}\wedge
\left \lbrack   \sum_{q=1}^{\vert N\vert +2}     d\overline{\Theta}_{q} \wedge
d\Theta_{q}  \right\rbrack^{\vert N\vert-1} \right\rbrace,
\end{eqnarray*}
and
\begin{eqnarray*}
I_2   & := &\left ( \sum\limits_{k=1}^{\vert N\vert +2}f_k
   d\overline{\Theta}_k \right )\wedge
\left\lbrace
 \frac{(1- \Theta\bullet\overline{Z})^{\vert N\vert -2} }{\vert \zeta
\vert^{2}  \vert \eta
\vert^{2}
D(\Theta,Z)^{\vert N\vert}} \sum_{j=1}^{\vert N\vert +2} \left \lbrack\overline{\Theta}_{j} (1-
\Theta\bullet \overline{Z}) - \overline{Z}_{j}( 1-\vert \Theta \vert^{2})
\right\rbrack  \right .\\&  & \left. d\Theta_{j}
\wedge  \overline{\partial}  \vert \Theta\vert^2   \wedge
\left \lbrack \sum_{k=1}^{n+1} \overline{Z}_{k} d\Theta_{k} \right\rbrack\wedge \Psi \wedge \overline{\partial(\zeta\bullet\zeta)}
  \wedge \overline{\partial(\eta\bullet\eta)}\wedge
\left \lbrack    \sum_{q=1}^{\vert N\vert +2}     d\overline{\Theta}_{q} \wedge
d\Theta_{q}  \right\rbrack^{\vert N\vert-2} \right\rbrace.
 \end{eqnarray*}
A  straightforward  computation   shows that
\begin{equation}\begin{split}
I_1 &:=
\frac{(1- \Theta\bullet\overline{Z})^{\vert N\vert -1} }{ \vert \zeta
\vert^{2}  \vert \eta
\vert^{2} D(\Theta,Z)^{\vert N\vert}} \left \lbrace C
\sum_{k=1}^{l}  f_k
\left \lbrack  -\overline{\xi}_{k} (1- \Theta\bullet
\overline{Z}) + \overline{x}_{k}( 1-\vert \Theta \vert^{2}) \right\rbrack
\right. \\ &  (z\bullet\overline{\zeta}+\vert \zeta\vert^{2}) (w\bullet\overline{\eta}+\vert \eta\vert^{2})
+ C
\sum_{k=1}^{n+1}  f_{k+l} \left\lbrace
\left \lbrack  -\overline{\zeta}_{k} (1- \Theta\bullet
\overline{Z}) + \overline{z}_{k}( 1-\vert \Theta \vert^{2}) \right\rbrack \right.
  \\ &   \left.
(z\bullet\overline{\zeta}+\vert \zeta\vert^{2})
-(1- \vert
\zeta\vert^{2})(z_{k}+\zeta_{k}) \overline{z\bullet\zeta}\right\rbrace
 (w\bullet\overline{\eta}+\vert \eta\vert^{2})\\
  & +\left.  C
\sum_{k=1}^{m+1}  f_{k+l+n+1} \left\lbrace
\left \lbrack  -\overline{\eta}_{k} (1- \Theta\bullet
\overline{Z}) + \overline{w}_{k}( 1-\vert \Theta \vert^{2}) \right\rbrack
(w\bullet\overline{\eta}+\vert \eta\vert^{2}) \right
.   \right
.  \\
&-\left.  \left.(1- \vert
\eta\vert^{2})(w_{k}+\eta_{k}) \overline{w\bullet\eta}\right\rbrace
 (z\bullet\overline{\zeta}+\vert \zeta\vert^{2})\right\rbrace d\Theta\wedge d\overline\Theta
.\end{split}
\end{equation}
 Hence the functions $h_k$ associated with $I_1$  are  of the  form (3.3).

To  simplify  notations we set
\begin{eqnarray*}
\omega_{\xi} :=  \sum\limits_{k=1}^{l} d\xi_k  \wedge d\overline{\xi_k },\quad
\omega_{\zeta} :=  \sum\limits_{k=1}^{n+1} d\zeta_k  \wedge  d\overline{\zeta_k },\quad
\omega_{\eta} :=  \sum\limits_{k=1}^{m+1} d\eta_k  \wedge d\overline{\eta_k },
 \end{eqnarray*}
 and  we  set  for  every   form $\omega$ and  every positive  integer $k,$
 $$\omega^k:=\underbrace{\omega\wedge \ldots\wedge \omega}_{k}.$$
 Then  a  simple  calculation  gives that
 \begin{equation}
 \begin{split}
\Psi &\wedge \overline{\partial(\zeta\bullet\zeta)}
  \wedge \overline{\partial(\eta\bullet\eta)}\wedge
\left \lbrack    \sum\limits_{k=1}^{\vert N\vert +2}     d\overline{\Theta}_{k} \wedge
d\Theta_{k}  \right\rbrack^{\vert N\vert-2}= \Psi \wedge \overline{\partial(\zeta\bullet\zeta)}
  \wedge \overline{\partial(\eta\bullet\eta)}\\& \wedge
\left (C \omega_{\xi}^{l-2}\wedge \omega_{\zeta}^{n} \wedge\omega_{\eta}^{m}
+C \omega_{\xi}^{l}\wedge \omega_{\zeta}^{n-2}\wedge \omega_{\eta}^{m} +C \omega_{\xi}^{l}\wedge \omega_{\zeta}^{n} \wedge \omega_{\eta}^{m-2} \right. \\
&+ \left. C \omega_{\xi}^{l-1}\wedge \omega_{\zeta}^{n-1}\wedge \omega_{\eta}^{m}
+C \omega_{\xi}^{l}\wedge \omega_{\zeta}^{n-1}\wedge \omega_{\eta}^{m-1}
+C \omega_{\xi}^{l-1}\wedge \omega_{\zeta}^{n}\wedge \omega_{\eta}^{m-1} \right ) \\
&\equiv  \Psi \wedge \overline{\partial(\zeta\bullet\zeta)}
  \wedge \overline{\partial(\eta\bullet\eta)}\wedge  \left (\sum\limits_{k=1}^{6} J_{k}\right ).
\end{split}\end{equation}
To  conclude the proof of Lemma 3.4, it  suffices  to prove  the following lemma~:
\begin{lem}  For  every $1\leq  k\leq 6,$   the differential form
 \begin{eqnarray*}
 I_{2k}&:=&\left \lbrack\sum\limits_{j=1}^{ \vert N\vert +2 }f_j  d\overline{\Theta}_{j} \right\rbrack
\wedge  \sum_{j=1}^{\vert N\vert +2} \left \lbrack\overline{\Theta}_{j} (1-
\Theta\bullet \overline{Z}) - \overline{Z}_{j}( 1-\vert \Theta \vert^{2})
\right\rbrack d\Theta_{j}  \wedge  \overline{\partial}  \vert \Theta\vert^2  \\ &\quad& \wedge
\left \lbrack \sum_{j=1}^{\vert N\vert +2} \overline{Z}_{j} d\Theta_{j} \right\rbrack\wedge  J_k\wedge  \Psi  \wedge \overline{\partial(\zeta\bullet\zeta)}
 \wedge \overline{\partial(\eta\bullet\eta)}
\end{eqnarray*}
can be  expressed as the  product of the canonical volume form  $ d\Theta\wedge d\overline{\Theta }$  and  a function of the form
 \begin{eqnarray*}
 \sum\limits_{j=1}^{ \vert N\vert +2 }f_j\left (  (1-
\Theta\bullet \overline{Z}) P_j(\Theta,Z)+( 1-\vert \Theta \vert^{2}) Q_j(\Theta,Z)  \right ),
 \end{eqnarray*}
 where $P_j,Q_j$  are some   polynomials satisfying
  assertion (ii)  of  Theorem 3.3.
\end{lem}
\noindent{\it End  of the proof of Lemma 3.4.}
Suppose that Lemma  3.5 is  proved.  We deduce  from   the definition of $I_2, I_{2k}$ and   (3.5)  that
 $$ I_2=
\frac{(1- \Theta\bullet\overline{Z})^{\vert N\vert -2} }{     \vert \zeta\vert^2 \vert \eta\vert^2    D(\Theta,Z)^{\vert N\vert}}\cdot
\left ( \sum\limits_{k=1}^{6} I_{2k} \right ).$$
Therefore Lemma 3.4  follows from  Lemma 3.5.  \hfill
$\square$\\

\smallskip

\noindent{\it  Proof of Lemma 3.5.}
We break   the proof  into 6 cases according to the integer  $k,$  $1\leq k\leq 6.$\\
\noindent{\bf Case 1:} $J_1=\omega_{\xi}^{l-2}\wedge \omega_{\zeta}^{n}\wedge \omega_{\eta}^{m} .$  In this  case  a  direct
computation  shows that
 \begin{eqnarray*}  I_{21} & = &( z\bullet\overline{\zeta}+\vert \zeta\vert^{2})
 ( w\bullet\overline{\eta}+\vert \eta\vert^{2})
  \left \lbrack\sum\limits_{j=1}^{ l }f_j  d\overline{\xi_{j}} \right\rbrack  \\
& \quad  & \wedge \left \lbrack\sum_{j=1}^{l} \left \lbrack\overline{\xi}_{j} (1-
\Theta\bullet \overline{Z}) - \overline{x}_{j}( 1-\vert \Theta \vert^{2})
\right\rbrack d\xi_{j}  \right\rbrack
        \wedge  \overline{\partial}  \vert \xi\vert^2   \wedge
\left \lbrack \sum_{j=1}^{l} \overline{x}_{j} d\xi_{j} \right\rbrack \wedge\omega_{\xi}^{l-2}\\
& =&  (1- \Theta\bullet\overline{Z})( z\bullet\overline{\zeta}+\vert \zeta\vert^{2})
 ( w\bullet\overline{\eta}+\vert \eta\vert^{2})
  \left \lbrace \left \lbrack\sum\limits_{j=1}^{ l }f_j d\overline{\xi_{j}} \right\rbrack \right.  \\
&\quad  &  \wedge\left. \left \lbrack\sum_{j=1}^{l} \overline{\xi}_{j} d\xi_{j}  \right\rbrack
        \wedge
\left \lbrack \sum_{j=1}^{l} \overline{x}_{j} d\xi_{j} \right\rbrack   \wedge \overline{\partial}  \vert \xi\vert^2   \wedge  \omega_{\xi}^{l-2}\right\rbrace.
\end{eqnarray*}
Since
\begin{eqnarray*}
 \left \lbrack\sum_{j=1}^{l} \overline{\xi}_{j} d\xi_{j}  \right\rbrack
        \wedge
\left \lbrack \sum_{j=1}^{l} \overline{x}_{j} d\xi_{j} \right\rbrack  =\sum_{j,k=1, j<k}^{l} \left(\overline{\xi}_{j}
\overline{x}_{k}-\overline{\xi}_{k}
\overline{x}_{j}\right )   d\xi_{j}\wedge d\xi_{k} ,
\end{eqnarray*}
 we  see easily that $I_{21}$ satisfies the conclusion of the  lemma.
 \\
 \noindent{\bf Case 2:} $J_2=\omega_{\xi}^{l}\wedge \omega_{\zeta}^{n-2}\wedge \omega_{\eta}^{m} .$
 In this case  we can  rewrite $I_{22}$  in the form
 \begin{eqnarray*}
 ( w\bullet\overline{\eta}+\vert \eta\vert^{2})\left \lbrace \left \lbrack\sum\limits_{j=1}^{ n +1 }f_{j+l} d\overline{\zeta_{j}} \right\rbrack
\wedge  \sum_{j=1}^{n +1} \left \lbrack\overline{\zeta}_{j} (1-
\Theta\bullet \overline{Z}) - \overline{z}_{j}( 1-\vert \Theta \vert^{2})
\right\rbrack d\zeta_{j}  \wedge  \overline{\partial}  \vert \zeta\vert^2 \right . \\  \left .\wedge
\left \lbrack \sum_{j=1}^{n+1} \overline{z}_{j} d\zeta_{j} \right\rbrack\wedge  \left \lbrack\sum_{k,p=1}^{n+1}(z_k+\zeta_k) \overline{\zeta_p}   d\zeta_k\wedge d\overline{\zeta_p}   \right \rbrack\wedge \omega_{\zeta}^{n-2} \right\rbrace
 \wedge \omega_{\xi}^{l} \wedge \omega_{\eta}^{m+1}.
\end{eqnarray*}
 In view of the proof of Lemma 2.7 in \cite{VY1}, the differential form  in braces can be  expressed as  the product
 of $d\zeta \wedge d \overline{\zeta}$ and a function of the form
 \begin{eqnarray*}
 \sum\limits_{j=1}^{ n +1 }f_{j+l}\left (  (1-
\Theta\bullet \overline{Z}) S_j(\Theta,Z)+( 1-\vert \Theta \vert^{2}) T_j(\Theta,Z)  \right ),
 \end{eqnarray*}
 where  $S_j,T_j$  are some  polynomials such that
 \begin{eqnarray*}
 S_j(\Theta,Z)&=& O\left ( \vert z -\zeta\vert  \vert \zeta \vert^2
 \right ),\ \
 T_j(\Theta,Z)= O\left ( \vert z -\zeta\vert    \vert \zeta \vert^2
 \right ),\\
 \text{grad}_{Z}  S_j(\Theta,Z)&=&O\left ( \vert \zeta \vert^2 +\vert \zeta\vert \vert z\vert
 \right)
 ,\ \ \text{grad}_{Z}  T_j(\Theta,Z)=O\left ( \vert \zeta \vert^2 +\vert \zeta\vert \vert z\vert \right)
 .
\end{eqnarray*}
Combining what we have proved  so far, we obtain  that $I_{22}$ satisfies the conclusion of the lemma.
\\
 \noindent{\bf Case 3:} $J_3=\omega_{\xi}^{l}\wedge \omega_{\zeta}^{n}\wedge \omega_{\eta}^{m-2} .$   This case can be
 treated in the same way  as  the previous case.\\
 \noindent{\bf Case 4:} $J_4=\omega_{\xi}^{l-1}\wedge \omega_{\zeta}^{n-1}\wedge \omega_{\eta}^{m} .$
 Then we have
   \begin{eqnarray*}
 I_{24}= \left \lbrack\sum\limits_{j=1}^{ l+n +1 }f_j d\overline{\Theta}_{j} \right\rbrack
\wedge \left \lbrace \sum_{t=1}^{l+n +1} \left \lbrack\overline{\Theta}_{t} (1-
\Theta\bullet \overline{Z}) - \overline{Z}_{t}( 1-\vert \Theta \vert^{2})
\right\rbrack d\Theta_{t}  \wedge  \overline{\partial}(  \vert \xi\vert^2 + \vert \zeta\vert^2 )  \right . \\  \left .\wedge
\left \lbrack \sum_{s=1}^{l+n+1} \overline{Z}_{s} d\Theta_{s} \right\rbrack\wedge
 \left \lbrack\sum_{k,p=1}^{n+1}(z_k+\zeta_k) \overline{\zeta_p}   d\zeta_k\wedge d\overline{\zeta_p}
   \right \rbrack \wedge \omega_{\zeta}^{n-1} \wedge \omega_{\xi}^{l-1}\right\rbrace
         \wedge ( w\bullet\overline{\eta}+\vert \eta\vert^{2})\omega_{\eta}^{m+1}.
\end{eqnarray*}
By splitting  $ \sum\limits_{j=1}^{ l+n +1 }f_j d\overline{\Theta}_{j}$ into a sum of the  two  parts
$\sum\limits_{j=1}^{ l }f_j d\overline{\xi_{j}}$  and $\sum\limits_{k=1}^{ n +1 }f_{k+l} d\overline{\zeta_{k}},$ we  also
split   $I_{24}$ into  two corresponding  parts as
  $I_{24}=I_{241}+ I_{242}.$
 A little  calculation  gives that
 \begin{eqnarray*}
  && I_{241}=\sum\limits_{j=1}^{l} \sum\limits_{k,p=1}^{n+1} f_j(z_k+\zeta_k)\overline{\zeta_{p}}
     d\overline{\xi_j}\wedge d\zeta_k\wedge d\overline{\zeta_p}  \\
   & \wedge&\left\lbrace \left\lbrack \sum_{t=1}^{l+n+1} \left \lbrack\overline{\Theta}_{t} (1-
\Theta\bullet \overline{Z}) - \overline{Z}_{t}( 1-\vert \Theta \vert^{2})
\right\rbrack d\Theta_{t}    \right\rbrack \wedge  \overline{\partial} \vert \zeta\vert^2  \wedge
\left \lbrack \sum_{s=1}^{l+n+1} \overline{Z}_{s} d\Theta_{s} \right\rbrack \right\rbrace\\
&\wedge & \omega_{\xi}^{l-1} \wedge \omega_{\zeta}^{n-1}
         \wedge ( w\bullet\overline{\eta}+\vert \eta\vert^{2}) \omega_{\eta}^{m+1}\\
 &=&C(1-
\Theta\bullet \overline{Z})  ( w\bullet\overline{\eta}+\vert \eta\vert^{2})  \left\lbrace
\sum\limits_{j=1}^{l}f_j \left (\sum\limits_{k,p=1}^{n+1}
(z_k+\zeta_k)\overline{\zeta_p}\zeta_k \left (\overline{\xi_j}\overline{z_p}-
\overline{x_j}\overline{\zeta_p}  \right ) \right ) \right\rbrace \\
&\wedge & \omega_{\xi}^{l} \wedge \omega_{\zeta}^{n+1}
         \wedge  \omega_{\eta}^{m+1}.
 \end{eqnarray*}
Similarly,   since $ \overline{\partial} \vert \xi\vert^2= \sum\limits_{s=1}^{l} \xi_s d\overline{\xi_s},$
 we  obtain
 \begin{eqnarray*}
  && I_{242}=\sum\limits_{s=1}^{l} \sum\limits_{\overset{k,p=1}{k\not =p}}^{n+1} f_{k+l}(z_k+\zeta_k)\overline{\zeta_p}    \xi_s
      d\zeta_k\wedge  d\overline{\zeta_k}\wedge d\overline{\zeta_p}\wedge d\overline{\xi_s} \\
   & \wedge&\left\lbrace \left\lbrack \sum_{t=1}^{l+n+1} \left \lbrack\overline{\Theta}_{t} (1-
\Theta\bullet \overline{Z}) - \overline{Z}_{t}( 1-\vert \Theta \vert^{2})
\right\rbrack d\Theta_{t}    \right\rbrack   \wedge
\left \lbrack \sum_{t=1}^{l+n+1} \overline{Z}_{t} d\Theta_{t} \right\rbrack \right\rbrace\\
&\wedge & \omega_{\xi}^{l-1} \wedge \omega_{\zeta}^{n-1}
         \wedge ( w\bullet\overline{\eta}+\vert \eta\vert^{2}) \omega_{\eta}^{m+1}\\
 &=&C(1-
\Theta\bullet \overline{Z}) ( w\bullet\overline{\eta}+\vert \eta\vert^{2})\sum\limits_{k=1}^{n+1}  f_{k+l} \left\lbrace
\sum\limits_{s=1}^{l} \sum\limits_{\overset{p=1}{p\not = k}}^{n+1}
\left \lbrack (z_k+\zeta_k)\overline{\zeta_p}\xi_s \left (\overline{x_s}\overline{\zeta_p}-
\overline{\xi_s}\overline{z_p}  \right ) \right . \right . \\
& &\left. \left.  -(z_p+\zeta_p)\overline{\zeta_p}\xi_s \left (\overline{x_s}\overline{\zeta_k}-
\overline{\xi_s}\overline{z_k}  \right )  \right \rbrack \right \rbrace
\wedge  \omega_{\xi}^{l} \wedge \omega_{\zeta}^{n+1}
         \wedge  \omega_{\eta}^{m+1}.
 \end{eqnarray*}
It can be checked that  $I_{241},I_{242}$ and   $I_{24}$ satisfy the conclusion of the lemma.
 \\
 \noindent{\bf Case 5:} $J_5=\omega_{\xi}^{l}\wedge \omega_{\zeta}^{n-1}\wedge \omega_{\eta}^{m-1} .$ Observe that
 \begin{eqnarray*}
 I_{25}=\left \lbrack\sum\limits_{j=l+1}^{ \vert N\vert +2}f_j d\overline{\Theta}_{j} \right\rbrack
\wedge \left \lbrace  \sum_{t=l+1}^{\vert N\vert +2} \left \lbrack\overline{\Theta}_{t} (1-
\Theta\bullet \overline{Z}) - \overline{Z}_{t}( 1-\vert \Theta \vert^{2})
\right\rbrack d\Theta_{t}    \right . \\ \wedge  \overline{\partial}(  \vert \zeta\vert^2 + \vert \eta\vert^2 ) \wedge
\left \lbrack \sum_{s=l+1}^{\vert N\vert +2} \overline{Z}_{s} d\Theta_{s} \right\rbrack\wedge  \left \lbrack\sum_{k,p=1}^{n+1}(z_k+\zeta_k) \overline{\zeta_p}   d\zeta_k\wedge d\overline{\zeta_p}   \right \rbrack
 \\  \left . \left\lbrack\sum_{r,s=1}^{m+1}(w_r+\eta_r) \overline{\eta_s}   d\eta_r\wedge d\overline{\eta_s}   \right \rbrack
\wedge \omega_{\zeta}^{n-1} \wedge \omega_{\eta}^{m-1}     \right\rbrace
         \wedge \omega_{\xi}^{l}.
\end{eqnarray*}
Rewriting  $\sum\limits_{j=l+1}^{ \vert N\vert +2}f_j d\overline{\Theta}_{j}$ as the sum  of two  differential forms
$\sum\limits_{k=1}^{n+1}f_{k+l} d\overline{\zeta_{k}}$  and $\sum\limits_{j=1}^{m+1}f_{j+l+n+1} d\overline{\eta_{j}},$
we thus  divide $I_{25}$ into   two  corresponding  terms: $I_{25}=I_{251}+ I_{252}.$
 A straightforward  computation shows that
 \begin{equation}\begin{split}
    & I_{251}=\sum\limits_{k,p,q=1}^{n+1} \sum\limits_{r,s=1}^{m+1} f_k(z_q+\zeta_q)\overline{\zeta_{p}} (\eta_r+ w_r)\overline{\eta_s}
     d\overline{\zeta_k}\wedge d\zeta_q\wedge d\overline{\zeta_p} \wedge d\eta_r\wedge d\overline{\eta_s} \\
   & \wedge\left\lbrace \left\lbrack \sum_{t=l+1}^{\vert N\vert +2} \left \lbrack\overline{\Theta}_{t} (1-
\Theta\bullet \overline{Z}) - \overline{Z}_{t}( 1-\vert \Theta \vert^{2})
\right\rbrack d\Theta_{t}    \right\rbrack \wedge  \overline{\partial} \vert \eta\vert^2  \wedge
\left \lbrack \sum_{s=l+1}^{\vert N\vert +2} \overline{Z}_{s} d\Theta_{s} \right\rbrack \right\rbrace\\
&\wedge \omega_{\zeta}^{n-1} \wedge \omega_{\eta}^{m-1}
         \wedge \omega_{\xi}^{l}\\
&=
 (1-
\Theta\bullet \overline{Z})\vert  \eta\bullet(w-\eta)\vert^2\left\lbrace \sum\limits_{k=1}^{n+1} f_k\left\lbrack C \sum\limits_{p\not = k}(z_p+\zeta_p)
\overline{\zeta_{p}}\overline{\zeta_{k}}+ C
(z_k+\zeta_k) \vert \zeta_k\vert^2  \right\rbrack
 \right \rbrace\\
&\wedge  \omega_{\zeta}^{n+1} \wedge \omega_{\eta}^{m+1}
         \wedge \omega_{\xi}^{l}.
 \end{split}\end{equation}
 We obtain  in exactly  the  same  way  an explicit expression for $I_{252}.$  Finally, we deduce from these
 expressions that
  $I_{251},I_{252}$ and $I_{25}$ satisfy the conclusion of the lemma.\\
  \noindent{\bf Case 6:} $J_6=\omega_{\xi}^{l-1}\wedge \omega_{\zeta}^{n}\wedge \omega_{\eta}^{m-1} .$   This last  case can be  treated in the same  way as  Case 4. The proof of  Lemma 3.5 is  therefore  complete.
\hfill
$\square$\\

\smallskip

 We end  this  section  with  the study of the particular case $N=(2,2).$  In this case  we   write
 for $Z,\Theta\in\B_N:$
 \begin{eqnarray*}
 Z\equiv (z,w)\equiv(z_1,z_2,z_3,w_1,w_2,w_3),\qquad \text{and}\qquad \Theta\equiv (\zeta,\eta)\equiv(\zeta_1,\zeta_2,\zeta_3,\eta_1,\eta_2,\eta_3).
 \end{eqnarray*}

To  establish  optimal Lipschitz estimates for the domain
  $\Omega_{ (2,2)},$  we need a  more precise  formulation of the Cauchy type formula
   given in   Theorem 3.3.
\begin{thm}  Let  $N:=(2,2).$
There are polynomials $R(\Theta,Z)$ and  $P_{jk}(\Theta,Z),$ \\ $Q_{jk}(\Theta,Z),\ 1\leq j\leq 2,\ 1\leq k\leq 4,$
that  satisfy  the following properties:\\
(i)  $
R(\Theta,Z)= \left(   C\vert  \zeta\vert^2 +C\vert  \eta\vert^2 \right )
      (\vert \zeta \vert^2  +z\bullet \overline{ \zeta} )
(\vert \eta \vert^2  +w\bullet \overline{ \eta} ).
$\\
(ii)   For  every $Z,\Theta\in \B_N,$  and  for  every $1\leq j\leq 2$ and $1\leq k\leq 4,$
\begin{eqnarray*}P_{1k}(\Theta,Z)&=& O\left ( \vert \Theta -Z\vert(\vert \zeta_3\vert +\vert \eta_3\vert +\vert \Theta -Z\vert )
(\vert \zeta \vert^2  +\vert \zeta\vert \vert z\vert)
(\vert \eta \vert^2  +\vert\eta\vert \vert w\vert ) \right) ,\\
Q_{1k}(\Theta,Z)&=&  O\left (\vert \Theta -Z\vert
 (\vert \zeta_3\vert +\vert \eta_3\vert +\vert \Theta -Z\vert )
(\vert \zeta \vert^2  +\vert \zeta\vert \vert z\vert)
(\vert \eta \vert^2  +\vert\eta\vert \vert w\vert ) \right),
\end{eqnarray*}
\begin{eqnarray*}
\left \vert  \text{grad}_{Z}\ P_{jk}(\Theta,Z)\right\vert &=& O\left ( (\vert \zeta \vert^2  +\vert \zeta\vert \vert z\vert)
(\vert \eta \vert^2  +\vert\eta\vert \vert w\vert ) \right), \\
\left \vert  \text{grad}_{Z}\ Q_{jk}(\Theta,Z)\right\vert &=&O\left ( (\vert \zeta \vert^2  +\vert \zeta\vert \vert z\vert)
(\vert \eta \vert^2  +\vert\eta\vert \vert w\vert ) \right).
 \end{eqnarray*}
(iii)
  Let   $u \in
\mathcal{C}^{1}(\overline{\M}_N)$   and
$f:=f_1 d\overline{\zeta}_1+f_2 d\overline{\zeta}_2+f_3 d\overline{\eta}_1+f_4 d\overline{\eta}_2
 $ is  a $(0,1)$-form with coefficients in $\mathcal{C}(\overline{\M}_N)$ that  satisfy  $ \overline{\partial}_{\M_N}u=  f|_{\M_N},$     then  for  every $Z\in \M_N,$
\begin{eqnarray*}& &\quad  u(Z)  = \int\limits_{\partial\M_N} \frac{R(\Theta,Z)}{  \left (1-Z\bullet \overline{\Theta}\right )^{\vert N\vert}} u(\Theta)
\frac{ d\sigma(\Theta)}{  \vert \zeta\vert^2   \vert \eta\vert^2} + \\
 & &\int \limits_{\M_N}\sum\limits_{j=1}^{2} \sum\limits_{k=1}^{4}
\frac{(1- \Theta\bullet\overline{Z})^{1+j} }{ D(\Theta,Z)^{4}}\lbrack (
1-\Theta\bullet \overline{Z}) P_{jk}(\Theta,Z) +( 1-\vert \Theta
\vert^{2})Q_{jk}(\Theta,Z)\rbrack
   f_k(\Theta)\frac{ dV(\Theta)}{  \vert \zeta\vert^2   \vert \eta\vert^2}.
\end{eqnarray*}
\end{thm}
\begin{proof}
We return to  the arguments used in the proof of Theorem  3.3.  By the  hypothesis on $f$  and   (3.2), we have that
\begin{eqnarray}
  I:= \overline{\partial} u\wedge K_0 =(f_1 H_1+ f_2 H_2+f_3 H_3 +f_4 H_4)
  d\Theta   \wedge d\overline{\Theta},
\end{eqnarray}
with  $ H_1:=h_1,H_2:=h_2,H_3:=h_4$ and $H_4:=h_5.$
 To  complete the proof, it suffices to  prove  the following
\begin{lem}
The functions $H_k$  in formula (3.7) can be expressed in the form
\begin{equation}
 H_k(\Theta,Z)=\sum\limits_{j=1}^{2}
\frac{(1- \Theta\bullet\overline{Z})^{1+j} }{ \vert \zeta\vert^2 \vert \eta\vert^2      D(\Theta,Z)^{4}}\lbrack (
1-\Theta\bullet \overline{Z}) P_{jk}(\Theta,Z) +( 1-\vert \Theta
\vert^{2})Q_{jk}(\Theta,Z)\rbrack ,
\end{equation}
where  $P_{jk}$ and $Q_{jk}$ are some   polynomials satisfying  assertion  (ii)  of
the theorem.
\end{lem}
\noindent{\it End of the proof  of Theorem  3.6.}
Suppose that the lemma is proved.   Using the  arguments  that    precede  Lemma 3.4  in the proof of  Theorem
3.3 and  applying  Proposition 2.5,
the theorem follows.
\end{proof}
\noindent{\it  Proof of Lemma 3.7.}
Following    the  proof of   Lemma 3.4,  we  write  $I= I_1+I_2.$   By virtue of   (3.4), the functions $H_k$  associated  to
$I_1$   (similarly to those    associated to   $I$ in (3.7)) are in the form
(3.8) with $j=2.$

Since $l=0$  and $\omega_{\xi}=0,$   formula (3.5) becomes
\begin{equation}
 \begin{split}
\Psi &\wedge \overline{\partial(\zeta\bullet\zeta)}
  \wedge \overline{\partial(\eta\bullet\eta)}\wedge
\left \lbrack    \sum\limits_{k=1}^{6}     d\overline{\Theta}_{k} \wedge
d\Theta_{k}  \right\rbrack^{2}= \Psi \wedge \overline{\partial(\zeta\bullet\zeta)}
  \wedge \overline{\partial(\eta\bullet\eta)}\\& \wedge
\left (  2 \omega_{\zeta} \wedge\omega_{\eta}
+  \omega_{\eta}^{2} + \omega_{\zeta}^{2}  \ \right ) \\
&\equiv  \Psi \wedge \overline{\partial(\zeta\bullet\zeta)}
  \wedge \overline{\partial(\eta\bullet\eta)}\wedge  \left (  J_1 + J_2+  J_3 \right ).
\end{split}\end{equation}
Therefore,   to conclude the proof of  Lemma 3.7, it suffices to prove the following
\begin{lem}  For every $1\leq  k\leq 3,$   the differential form
 \begin{eqnarray*}
 I_{2k}&:=&\left \lbrack  f_1 d\overline{\zeta}_1 +  f_2 d\overline{\zeta}_2+
  f_3 d\overline{\eta}_1+ f_4 d\overline{\eta}_2         \right\rbrack
\wedge  \sum_{r=1}^{6} \left \lbrack\overline{\Theta}_{r} (1-
\Theta\bullet \overline{Z}) - \overline{Z}_{r}( 1-\vert \Theta \vert^{2})
\right\rbrack d\Theta_{r}   \\
&\wedge  & \overline{\partial}  \vert \Theta\vert^2   \wedge
\left \lbrack \sum_{s=1}^{6} \overline{Z}_{s} d\Theta_{s} \right\rbrack\wedge  J_k\wedge  \Psi  \wedge \overline{\partial(\zeta\bullet\zeta)}
 \wedge \overline{\partial(\eta\bullet\eta)}
\end{eqnarray*}
can  be  expressed as the   product of the  canonical  volume  form   $ d\Theta\wedge d\overline{\Theta }$  and a function of
the form
 \begin{eqnarray*}
 \sum\limits_{t=1}^{ 4 }f_t\left (  (1-
\Theta\bullet \overline{Z}) P_{1t}(\Theta,Z)+( 1-\vert \Theta \vert^{2}) Q_{1t}(\Theta,Z)  \right ),
 \end{eqnarray*}
 where $P_{1t},Q_{1t}$  are some polynomials satisfying
  assertion (ii)  of the
theorem.
\end{lem}
\noindent{\it End of the proof of  Lemma  3.7.}
Suppose  that   Lemma 3.8 is proved.  In view of  (3.9) and     the expression of    $I_2$ given  at the beginning
of the proof of Lemma 3.4,  we see that $$ I_2=
\frac{(1- \Theta\bullet\overline{Z})^{2} }{ \vert \zeta\vert^2 \vert \eta\vert^2 D(\Theta,Z)^{ 4}}\cdot
\left ( \sum\limits_{k=1}^{3} I_{2k} \right ).$$ Therefore,
 Lemma 3.7 follows from Lemma 3.8. \hfill
$\square$\\

\noindent{\it  Proof of  Lemma 3.8.}
We  first  remark that  the case  $k=1$  corresponds to the  case 5  in the proof of   Lemma 3.5.
Hence, by  virtue  of  identity   (3.6),
 $I_{21}$  satisfies the conclusion of  the lemma.
Consider the case  $k=2$  which corresponds to  case 2 in the proof of   Lemma  3.5.  Then
we have
 \begin{eqnarray*}
 I_{22}= \left \lbrace \left \lbrack  f_1  d\overline{\zeta}_1 +       f_2 d\overline{\zeta}_2 \right\rbrack
\wedge  \sum_{r=1}^{3} \left \lbrack\overline{\zeta}_{r} (1-
\Theta\bullet \overline{Z}) - \overline{z}_{r}( 1-\vert \Theta \vert^{2})
\right\rbrack d\zeta_{r}  \wedge  \overline{\partial}  \vert \zeta\vert^2 \right . \\  \left .\wedge
\left \lbrack \sum_{s=1}^{3} \overline{z}_{s} d\zeta_{s} \right\rbrack\wedge  \left \lbrack\sum_{t,p=1}^{3}(z_t+\zeta_t) \overline{\zeta}_p d\zeta_t\wedge d\overline{\zeta}_p \right \rbrack  \right\rbrace
  \wedge ( w\bullet\overline{\eta}+\vert \eta\vert^{2})\omega_{\eta}^{3}.
\end{eqnarray*}
A simple  calculation   gives that
\begin{eqnarray*}
I_{22}&=& \left  \lbrace f_1 (1-\Theta\bullet \overline{Z})(C\overline{\zeta}_3\zeta_2 +C\zeta_3\overline{\zeta}_2)
\left \lbrack \sum (-1)^{\epsilon(r,s,t)} (\overline{\zeta}_r \overline{z}_s - \overline{\zeta}_s \overline{z}_r)
(z_t+\zeta_t) \right\rbrack \right . \\
&+&  \left.f_2 (1-\Theta\bullet \overline{Z})(C\overline{\zeta}_3\zeta_1 +C\zeta_3\overline{\zeta}_1)
\left \lbrack \sum (-1)^{\epsilon(r,s,t)} (\overline{\zeta}_r \overline{z}_s - \overline{\zeta}_s \overline{z}_r)
(z_t+\zeta_t) \right\rbrack \right\rbrace \\
&\cdot  &   ( w\bullet\overline{\eta}+\vert \eta\vert^{2})  \omega_{\zeta}^{3} \wedge \omega_{\eta}^{3},
 \end{eqnarray*}
where  the sum is taken   over  all   permutations $( r,s,t)$  of $\{1,2,3\}$   such that
$r<s$ and  where  $\epsilon(r,s,t)$ is the sign  of  such  permutations.
It  follows    from this    that   $I_{22}$ satisfies  the  conclusion
of the lemma. Similarly, we have the  same  conclusion  for  $I_{23},$  which completes the proof. \hfill
$\square$\\

\section{Local coordinate systems on  the complex manifolds $\H_n$ and $\H_m$}

In the  next theorem,  we  construct an open neighborhood  $\mathcal{U}_n$ of $\H_n$ in  $\C^{n +1},$ and  for  every $z\in  \mathcal{U}_n, $        a coordinate  chart $\Phi^z$ defined  on
a coordinate patch $\mathcal{U}(z)$ of $\H_n$  that
possess  some  interesting properties of  homogeneity.
The  same construction will be applied  to the complex  manifold   $\H_m.$   These local coordinate systems  will allow us
in the next section  to
 reduce  certain types of integral  estimates  over  $\M_N$ to   simpler integral estimates over $\C^{\vert N\vert}.$
 \begin{thm}
There are an open  neighborhood  $\mathcal{U}_n$ of $\H_n$ in $\C^{n +1}$ and constants  $ C_1, C_2,C_3 >1$
 that  satisfy  the  following properties~:
 \\
 1)  If $z\in \C^{n+1}\setminus\mathcal{U}_n$ then $\text{dist}(z,\H_n) > \frac{\vert z\vert }{C_1}$  with  the understanding
 that
  $\text{dist}(.,.)$ is   the euclidean distance.
\\
 2) If $z\in\mathcal{U}_n$ and if  the open set
   $\mathcal{U}(z):=\left\lbrace    \zeta\in\H_n: \vert \zeta -z\vert <\frac{\vert z\vert}{C_1}\right\rbrace$ is  non-empty,
   then there exists a diffeomorphism
   $\Phi^{z}$  mapping  $\mathcal{U}(z)$ into  the open neighborhood
    $\widetilde{\mathcal{U}}(\widetilde{z}):=\left\lbrace
      \widetilde{\zeta}\in\Cn: \vert \widetilde{\zeta} -\widetilde{z}\vert
       <\frac{\vert \widetilde{z}\vert}{C_2}\right\rbrace$  of a point
       $\widetilde{z}\in\C^n$ which is exactly $
  \Phi^{z}(z) $ in case  $z\in\H_n$   such that

(i)  $\zeta\bullet\overline{z}=
\Phi^z(\zeta)\bullet\overline{\Phi^z(z)},$      for  all $ \zeta\in\mathcal{U}(z) .$

(ii) $\vert \Phi^{z}(z) \vert=\vert z\vert $ and
$\frac{\vert \zeta\vert}{2} \leq \vert \Phi^{z}(\zeta) \vert
 \leq  \vert \zeta\vert,$         for  all $ \zeta\in\mathcal{U}(z).$

(iii) For  all  $ \zeta\in\mathcal{U}(z),$   we have
$\frac{\vert \zeta\vert}{\vert \zeta_{n+1}\vert}\leq C_3 \frac{\vert \Phi^{z}(\zeta)\vert}
{\vert \Phi^{z}_{n}(\zeta)\vert},$    where $\Phi^{z}:=\left(
   \Phi_{1}^{z},\ldots,\Phi_{n}^{z}\right).$

(iv) For  all compactly supported  functions $f\in\mathcal{C}_{0} (\mathcal{U}(z))$  such that $f\geq 0,$      we  have
$$\int\limits_{\mathcal{U}(z) }  fdV_n  \leq  C_3 \int\limits_{\widetilde{\mathcal{U}}(\widetilde{z}) }
(\Phi^{z}_{*} f)(\widetilde{\zeta})dV_n(\widetilde{\zeta}),$$
where  $dV_n(\widetilde{\zeta})$ denotes the Lebesgue  measure on $\C^n$  and   $\Phi^{z}_{*} f$ is the  pushforward
of  $f$  under the  diffeomorphism $\Phi^z.$
\end{thm}
\begin{rem} {\rm
We construct in the  same  way an open  neighborhood $\mathcal{U}_m$ of $\H_m$ in $\C^{m +1}, $
and   for  every $w\in \mathcal{U}_m ,$ a  coordinate chart  $\Phi^w$ defined on a coordinate patch   $\mathcal{U}(w)$ of $\H_m$   that  possess the same   properties  as
 $\mathcal{U}_n, \Phi^z$ and $\mathcal{U}(z).$ }
\end{rem}

 To prove  Theorem 4.1, we need the  following
 \begin{lem}
There exists a constant $C_0 >0$ such that
  $ \underset{1\leq j< k\leq n+1}{\max} \left\vert\text{Im} (z_j\overline{z_k})\right\vert > C_0,$
  for  all   $z:=(z_1,\ldots,z_{n+1})\in\partial \M_n.$  Here $\text{Im} \lambda$ denotes  the  imaginary part of  $\lambda\in\C.$
  \end{lem}
 \begin{proof}
 Since the   function  $z \mapsto \underset{1\leq j< k\leq n+1}{\max} \left\vert\text{Im} (z_j\overline{z_k})\right\vert $ is  continuous on the   compact set $\partial \M_n,$  it attains  its minimum at  a  point $z.$
 Therefore it  suffices to prove that there exist $1\leq j< k\leq n+1$  such that $\text{Im} ( z_j\overline{z_k} )\not = 0.$
 Suppose the contrary. Since $\vert z\vert =1,$  there is an $k$ such that $z_k\not =0.$
 Hence  for every $1\leq j\leq n+1,$ we have $ z_j=\lambda_j z_k$ with $ \lambda_j\in\R,$  from which it follows that
  $0= \left (\sum\limits_{j=1}^{n+1} \lambda_j^{2} \right ) z_k^{2}.$  Thus $z_{k}=0$ and we  obtain
 a contradiction. This completes the proof of the    lemma.
 \end{proof}
 We now  turn to the proof of   Theorem 4.1.

 \smallskip

 \noindent{\it Proof of  Theorem 4.1.}
The  construction of the  open neighborhood $\mathcal{U}_n,$ the coordinate patches  $\mathcal{U}(z)$ and the coordinate charts   $\Phi^{z}:\mathcal{U}(z)\longrightarrow \C^n$ for  every $ z\in \mathcal{U}_n,     $   will be done within two steps. First,
by   Lemma 4.3, we divide $\partial \M_n$ into $\frac{n(n+1)}{2}$ compact sets
 $E_{jk},\ j<k,$ where
 $
  E_{jk} := \left\lbrace  z\in\partial \M_n     : \left\vert\text{Im}( z_j\overline{z_k})\right\vert \geq C_0   \right\rbrace.
  $

   Fix a sufficiently small number   $\delta >0.$    The exact value  of $\delta$  will be  clear     in the course of the proof. Let  $z$ be a point of $\C^{n+1}.$
   \\ {\bf Step 1: } $\text{dist}(z,\partial\M_n)  <\delta.$

   According to the discussion above,  suppose without loss of  generality that
   there exist  $j<k$ and $\hat{z}\in E_{jk} $ such that  $\vert z -\hat{z}  \vert  <\delta.$
    Define the  diffeomorphism $\Phi^{z}$ as  follows~:
$\Phi^{z}:=\left(
   \Phi_{1}^{z},\ldots,\Phi_{n}^{z}\right),$  where
 \begin{eqnarray*}\left\lbrace
 \begin{array}{l}
\Phi_{k-1}^{z}(\zeta):=\frac{\zeta_j\bullet\overline{z}_j  +
\zeta_k\bullet\overline{z}_k
}{\sqrt{\vert z_j\vert^2+\vert z_k\vert^2}  };\\
\Phi_{l}^{z}(\zeta):=\zeta_l,\quad \text{if}\ l <j;\\
\Phi_{l}^{z}(\zeta):=\zeta_{l+1}, \quad \text{if}\  j\leq l< k-1\ \text{or }\ k \leq l\leq n.
\end{array} \right.
\end{eqnarray*}
 We can choose the   functions $\zeta_l,\ l\not= j,$ as
  the  $n$-local  coordinate  functions of  $\H_n$ at the point  $\hat{z}.$   Substituting $ \zeta_j$ by
  $i\sqrt{\sum\limits_{l\not= j}  \zeta_l^{2} }$ in the expression of $\Phi^z,$
   straightforward  computations show that  the real Jacobian
   of $\Phi^{z}$ at the point $\zeta$   corresponding to this local   coordinate   system is equal to
   $\frac{\left\vert \zeta_j\overline{z_k}-  \zeta_k\overline{z_j} \right\vert^{2}}
   { \vert \zeta_j\vert^2 (\vert z_j\vert^2+\vert z_k\vert^2)       }.$
    This  quantity  is    uniformly bounded from above  and from below by  some positive constants as $\zeta\in\H_n$ and $z$ are very near to
    $\hat{z}\in E_{jk}.$
   Therefore, when  $C_2$ is  sufficiently large, there exists a sufficiently  small $\delta$
   so that  for every $\hat{z}\in E_{jk}$ and every
    $z$  such that $\vert z -\hat{z}\vert <\delta,$    $\Phi^{z}$ is a diffeomorphism from
     $\left\lbrace  \zeta\in\H_n:  \vert \zeta - z\vert < 2\delta \right\rbrace$  to
     $\left\lbrace    \widetilde{\zeta}\in\Cn: \vert \widetilde{\zeta} - \Phi^{z}(z)
         \vert <\frac{1}{2C_2}\right\rbrace.$

    Taking $C_1 > \frac{1}{2\delta}$ and observing that $\vert \Phi^{z}(z)    \vert=\vert z\vert \approx 1,$
    it follows  from  the previous  discussion that
      $\Phi^{z}$ is a  diffeomorphism from $\mathcal{U}(z)$   onto an open neighborhood      $\widetilde{\mathcal{U}}(\widetilde{z})$ of the point $\widetilde{z}:=
  \Phi^{z}(z) \in \C^{n}.$

  To finish part 2) of the  theorem, it  remains to prove assertions (i)-(iv).

  Assertions (i) and (ii) can be  checked direcly. In particular,
  the estimate
  $\vert \Phi^{z}(\zeta) \vert
 \leq  \vert \zeta\vert$  follows      from the Cauchy-Schwarz inequality.

We prove now  assertion (iii). Consider two  cases  according to  $k~:$
\\
{\bf  Case }  $k< n+1.$ In this case, in view of the definition of $\Phi^z,$ we have
$ \zeta_{n+1}  = \Phi_{n}^{z}(\zeta) .$ This, combined with (ii), implies assertion (iii).
\\
{\bf  Case } $k=n+1.$  If $\zeta\in \mathcal{U}(z),$  then when    $C_1$  is sufficiently large,
  we  have   $ 1>\vert \zeta_{n+1}\vert \approx \vert \hat{z}_{n+1}\vert \geq C_0.$  Hence  assertion (iii) is
 almost obvious.

  It now remains to prove assertion  (iv).  By Proposition 2.1 in \cite{VY1}, for  $\zeta\in \mathcal{U}(z) $   we have the following identity:
\begin{eqnarray*}
 dV_n(\zeta)= C\frac{\vert \zeta\vert^{2}}{\vert
\zeta_{j}\vert^{2}}
d\zeta_{1}\wedge
d\overline{\zeta}_{1}\wedge  \ldots \wedge  \widehat{d\zeta_{j}}\wedge
\widehat{d\overline{\zeta}_{j}} \wedge \ldots \wedge d\zeta_{n+1}\wedge
d\overline{\zeta}_{n+1}.
\end{eqnarray*}
Since $ 2 > \vert \zeta \vert >\vert \zeta_{j}\vert \approx \vert \hat{z}_{j}\vert \geq C_0,$ it follows that
$ dV_n(\zeta) \approx (\Phi^{z})^{\ast} \left (dV_n(\widetilde{\zeta}) \right )$ for
$\widetilde{\zeta}= \Phi^{z}(\zeta).$  This implies  assertion (iv).
\\
{\bf  Step 2:  General case.}

Set $\mathcal{U}_n:=\left\lbrace
rz:  r>0\ \text{and}\ \text{dist}(z,\partial\M_n)  <\delta\right\rbrace.$

If $z\in \mathcal{U}_n,$ then according to the    definition above, there exist $\hat{z}\in\partial\M_n$ and $r >0$ such that $\left\vert rz - \hat{z} \right \vert  <\delta.$  Therefore, the construction given in  Step 1
can then be applied to   the point $r z .$ Hence,  we can define
\begin{eqnarray*}
\mathcal{U}(z)&:=&\frac{1}{r}\cdot\mathcal{U}(rz);\\
\Phi^{z}(\zeta)&:= &\frac{1}{r} \cdot  \Phi^{rz} ( r\zeta),
\qquad \forall \zeta\in\mathcal{U}(z).
\end{eqnarray*}
 Using the  homogeneous invariance  of the  complex  manifold $\H_n$ with respect to the dilations,  we conclude that
 for  every   $z\in \mathcal{U}_n,$ the function
$\Phi^{z}$ just defined  satisfies part 2) of the theorem.
To finish the proof of the theorem, it only remains to  check part 1). Let $z\not\in  \mathcal{U}_n.$ Then there exists
 a point $\hat{z}\in \H_n $  such that $\vert z -\hat{z}\vert =\text{dist}(z,\H_n).$
Since  $z\not\in  \mathcal{U}_n,$  we deduce that $\vert z -\hat{z}\vert \geq\delta\vert \hat{z}\vert.$  Hence
$$ \left ( \frac{1}{\delta}+1\right) \text{dist}(z,\H_n) > \vert z -\hat{z}\vert+\vert \hat{z}\vert\geq \vert z\vert.$$
Thus, if we choose $C_1 > \frac{1}{\delta}+1,$  then  part 1) is satisfied.      This completes  the proof of the   theorem.
\hfill
$\square$\\

\section{Reduction of estimates from $\M_N$ to $\B_{\vert N\vert }$}
This  section proves  the Theorem  of reduction of estimates.  We use
the   notations and the  constants introduced  in the previous section.
In order to state this theorem,  we  need some more notations and definitions.

We denote by $\B_{\vert N\vert}$ the euclidean unit ball of   $\C^{\vert N\vert}.$  We often use the following notations
for $\widetilde{\Theta},\widetilde{Z}\in \C^{\vert N\vert}:     $
\begin{eqnarray*}\widetilde{\Theta} \equiv   (\tilde{\xi},\tilde{\zeta},\tilde{\eta}) \in \C^l\times \C^n\times\C^m\  \ \text{and}\ \
\widetilde{Z} \equiv   (\tilde{x},\tilde{z},\tilde{w}) \in \C^l\times \C^n\times\C^m.
\end{eqnarray*}
Let  $dV(\widetilde{\Theta})$ be the
 Lebesgue measure on   $\C^{\vert N\vert}.$
 For every  $i\in\{1,2\},$ note $\B_{i,\vert N\vert}$ the euclidean ball of  $\C^{\vert N\vert}$ centered at the origin
 with radius $i.$ Thus $\B_{\vert N\vert}=\B_{1,\vert N\vert}.$

We shall   define  various    notions of comparability.
\begin{defi}  Consider  two points $Z\equiv   (x,z,w)\in \B_N$    and  $\widetilde{Z} \equiv   (\tilde{x},\tilde{z},\tilde{w})\in \C^{\vert N\vert}.$   $Z$ is  said  to be  {\it comparable} with $\widetilde{Z}$ if the  following  conditions are
true~:

 (1) $x=\widetilde{x}. $

 (2) If  $z\in\mathcal{U}_n$ and $\mathcal{U}(z)\not =\varnothing,$  then      $\widetilde{z}=
 \Phi^{z}(z),$ if not  $\vert \widetilde{z}\vert=\vert z\vert.$

 (3)  If $w\in\mathcal{U}_m$ and $\mathcal{U}(w)\not =\varnothing,$ then      $\widetilde{w}=
 \Phi^{w}(w),$ if not  $\vert \widetilde{w}\vert=\vert w\vert.$
\end{defi}
\begin{rem} {\rm
It should  be noted  that  by this  definition and    Theorem 4.1 (ii), we have
$\vert x\vert =\vert  \widetilde{x}      \vert, \vert z\vert=\vert  \widetilde{z}      \vert,
 \vert w\vert=\vert  \widetilde{w}      \vert.$  Hence $\widetilde{Z}\in \B_{\vert N\vert}.$ }
\end{rem}
\begin{defi} Let $i\in\{1,2\}$  and fix two  comparable  points   $Z\equiv   (x,z,w)\in \B_N$    and $\widetilde{Z} \equiv   (\tilde{x},\tilde{z},\tilde{w})\in \C^{\vert N\vert}.$

We say that  $\xi\in \C^l$  is {\it i-comparable} with $\widetilde{\xi}\in \C^l$  if  $\xi=\widetilde{\xi}. $

We say  that  $\zeta\in \H_n$  is {\it i-comparable} with $\widetilde{\zeta}\in \C^{n}$  if  the following  conditions  are
true~:

(1) If $\vert \zeta\vert >   \sqrt{2}\vert z\vert,$ then   $\vert \widetilde{\zeta}\vert=\vert \zeta\vert.$

 (2) If   $\vert \zeta -z \vert  < \frac{\vert z\vert}{C_1},$  then
 $\widetilde{\zeta}=\Phi^{z}(\zeta).$

 (3) If  $\vert \zeta\vert\leq   \sqrt{2}\vert z\vert $ and $\vert \zeta -z \vert  \geq \frac{\vert z\vert}{C_1},$ then
 $  \vert \widetilde{\zeta}\vert\leq   \sqrt{2}\vert \widetilde{ z}\vert $ and $\vert \widetilde{\zeta} -\widetilde{z} \vert  > \frac{\vert \widetilde{z}\vert}{C_2};$ if moreover $i=1,$  then we have   $\vert \widetilde{\zeta}\vert \leq \vert \zeta\vert.$

We can define  in the  same way  the  notion of i-comparability between $\eta\in \H_m$ and
  $\widetilde{\eta}\in \C^{m}$  upon   substituting $n$ by $m$ and $\Phi^{z}$ by $\Phi^{w}.$

  Finally,
two  points $\Theta\equiv (\xi,\zeta,\eta)\in\overline{\M}_N$ and
$\widetilde{\Theta}\equiv \left(\widetilde{\xi},\widetilde{\zeta},\widetilde{\eta}\right)\in \C^{\vert N\vert}$ are  said   to be  {\it i-comparable }  if $\xi$ (resp. $\zeta$ and $\eta$) is i-comparable with
$\widetilde{\xi}$ (resp.  $\widetilde{\zeta}$ and $\widetilde{\eta}$).
\end{defi}
\begin{rem} {\rm We deduce easily from  this definition and    Theorem 4.1 (ii) that
if $\Theta\in\overline{\M}_N$ is i-comparable  with $\widetilde{\Theta}\in\C^{\vert N\vert},      $
then  $\widetilde{\Theta}\in\overline{\B}_{i,\vert N\vert}.$ }
\end{rem}
\begin{defi} Let  $i\in\{1,2\}$ and fix  two  comparable points   $Z\in \B_N$    and $\widetilde{Z} \in \B_{\vert N\vert}.$
 Consider    two non-negative   measurable   functions  $ K,\widetilde{K}$      defined respectively
 on   $\overline{\M}_N$ and $\overline{\B}_{i,\vert N\vert}.$
 \begin{itemize}
 \item[$\bullet$]
  We  write  $K\lesssim  C\widetilde{K}$ (respectively,  $\widetilde{K}\lesssim  CK$ ) at $(Z,\widetilde{Z})$ for a positive  constant $C$
   if   for  all
 points   $\Theta\in\overline{\M}_N$  i-comparable with
$\widetilde{\Theta}\in \overline{\B}_{i,\vert N\vert},$
\begin{eqnarray*} K(\Theta)\leq C\widetilde{K}( \widetilde{\Theta}  )\qquad \text{(respectively,}\
\widetilde{K}( \widetilde{\Theta}  )
 \leq C  K(\Theta)) .
    \end{eqnarray*}
    \item[$\bullet$]   We  write  $K\approx\widetilde{K}$   at $(Z,\widetilde{Z})$
   if there exists  $C >0 $   such that $K\lesssim C \widetilde{K}\lesssim C^{2} K.$
   \end{itemize}
   \end{defi}

Now we  are in a position   to   state  the  main  theorem  of this  section.
\begin{thm} Let $i\in\{1,2\}$ and fix   two  comparable    points   $Z\in \B_N$    and $\widetilde{Z} \in \B_{\vert N\vert}.$
 Let  $C$ be  a positive constant. Consider non-negative   measurable    functions  $ K,L$   defined  on $\overline{\M}_N$  and $\widetilde{K},\widetilde{L}$  defined on  $\overline{\B}_{i,\vert N\vert}$  such  that
 $$K\lesssim  C\widetilde{K} \quad\text{and}\quad L\lesssim  C\widetilde{L}\lesssim  C^2 L
 \quad\text{at}\  (Z,\widetilde{Z}).$$
 For  every  $ \alpha:= \left ( \alpha_{1},\alpha_{2},\alpha_{3},\alpha_{4}  \right)$    such that
 $0\leq \alpha_1 < 2n,$ $0\leq\alpha_2 < 2m$ and  $0\leq \alpha_3,\alpha_4 < 2,$   we  set
 \begin{eqnarray*} K_{ \alpha    }(\Theta)&:=&K(\Theta)
  \left (  1 +  \frac{ \vert z \vert}{ \vert\zeta \vert}
 \right )^{\alpha_1} \left (  1 +  \frac{ \vert w \vert}{ \vert\eta \vert}
\right )^{\alpha_2}
      \left \vert\frac{ \zeta}
{ \zeta_{n+1}  }\right \vert^{\alpha_3}
 \left \vert\frac{ \eta}
{ \eta_{m+1}  }\right\vert^{\alpha_4}, \\
\widetilde{K}_{1, \alpha    }(\widetilde{\Theta})&:=&  \widetilde{K}( \widetilde{\Theta})
 \left (  1 +  \frac{ \vert\widetilde{  z} \vert}{ \vert\widetilde{\zeta} \vert}
 \right )^{\alpha_1} \left (  1 +  \frac{ \vert\widetilde{w} \vert}{ \vert\widetilde{\eta} \vert}
 \right )^{\alpha_2}
\left \vert\frac{ \widetilde{\zeta}}
{ \widetilde{\zeta}_{n}  }\right \vert^{\alpha_3}
\left \vert\frac{ \widetilde{\eta}}
{ \widetilde{\eta}_{m}  }\right\vert^{\alpha_4},\\
\widetilde{K}_{2, \alpha    }(\widetilde{\Theta})&:=&  \widetilde{K}( \widetilde{\Theta})\left \vert\frac{ \widetilde{\zeta}}
{ \widetilde{\zeta}_{n}  }\right \vert^{\alpha_3}
\left \vert\frac{ \widetilde{\eta}}
{ \widetilde{\eta}_{m}  }\right\vert^{\alpha_4}.
\end{eqnarray*}
   Then  there   exists a constant $C_4$  that   depends only on  $N,\alpha$ and $C,C_1,C_2,C_3 , $ (in particular  this constant is independant of $Z$ and $\widetilde{Z}$), such that \\
1)
\begin{eqnarray*}
\int\limits_{\M_N} K_{ \alpha    }(\Theta) dV(\Theta) \leq C_4
\int\limits_{\B_{i,\vert N\vert}}  \widetilde{K}_{ i,\alpha    }(\widetilde{\Theta})      dV(\widetilde{\Theta});
\end{eqnarray*}
2) for $\delta >  0,$
\begin{eqnarray*}
\int\limits_{\Theta\in \M_N,\ L( \Theta)  \leq \delta} K_{ \alpha    }(\Theta) dV(\Theta) \leq C_4
\int\limits_{ \widetilde{\Theta}\in   \B_{i,\vert N\vert},\  \widetilde{L}(\widetilde{\Theta})
\leq C_4\delta
}  \widetilde{K}_{ i,\alpha    }(\widetilde{\Theta})      dV(\widetilde{\Theta});
\end{eqnarray*}
3) for $0 < \delta_1 < \delta_2,$
\begin{eqnarray*}
\int\limits_{\Theta\in \M_N,\ \delta_1 \leq  L( \Theta)  \leq \delta_2
} K_{ \alpha    }(\Theta) dV(\Theta) \leq C_4
\int\limits_{ \widetilde{\Theta}\in   \B_{i,\vert N\vert},\ \frac{\delta_1}{C_4}  \leq  \widetilde{L}(\widetilde{\Theta})
\leq C_4\delta_2
}  \widetilde{K}_{ i,\alpha    }(\widetilde{\Theta})      dV(\widetilde{\Theta}).
\end{eqnarray*}
\end{thm}
%
%
%
%
\begin{proof}
We  shall  only prove  part 3).  The two  other  assertions can be  shown in exactly the same  way.
Firstly, we  extend the domain of definition of the  functions $K,L,\widetilde{K},\widetilde{L}$ by setting
\begin{eqnarray*}
K(\Theta)&=&L(\Theta):= 0,\qquad \text{if}\ \Theta\in \H_N\setminus\overline{\M}_N;\\
\widetilde{K}(\widetilde{\Theta})&=&\widetilde{L}(\widetilde{\Theta}):= 0,\qquad \text{if}\ \widetilde{\Theta}\in \C^{\vert N\vert} \setminus\overline{\B}_{i,\vert N\vert}.
\end{eqnarray*}
By the hypothesis on   $L$ and $\widetilde{L},$  for  every $\Theta\in\overline{\M}_N$  such that
 $\delta_1 \leq  L( \Theta)  \leq \delta_2$ and for  every
$\widetilde{\Theta}\in \overline{\B}_{i,\vert N\vert}$ i-comparable with  $\Theta,$ we  have
\begin{eqnarray}
\frac{\delta_1}{C}  \leq  \widetilde{L}(\widetilde{\Theta})
\leq C\delta_2.
\end{eqnarray}
For  every $\xi,\widetilde{\xi}\in\B_l$ and $\eta,\widetilde{\eta}\in \H_m,$
  consider  the following  integrals
\begin{eqnarray*}
R(\xi,\eta)  &:=&\int\limits_{\zeta\in \H_n,\ \delta_1 \leq  L(\xi,\zeta,\eta      )  \leq \delta_2
} K(\xi,\zeta,\eta)  \left (  1 +  \frac{ \vert z \vert}{ \vert\zeta \vert}
 \right )^{\alpha_1}
      \left \vert\frac{ \zeta}
{ \zeta_{n+1}  }\right \vert^{\alpha_3}
  dV_n(\zeta),  \\
 \widetilde{R}_1( \widetilde{\xi},\widetilde{\eta})  &:=&
\int\limits_{ \widetilde{\zeta}\in   \C^n,\ \frac{\delta_1}{C}  \leq  \widetilde{L}\left(\widetilde{\xi},\widetilde{\zeta},\widetilde{\eta}        \right  )
\leq C\delta_2
}  \widetilde{K}\left(\widetilde{\xi},\widetilde{\zeta},\widetilde{\eta}        \right  ) \left (  1 +  \frac{ \vert\widetilde{ z} \vert}{ \vert\widetilde{\zeta} \vert}
 \right )^{\alpha_1}
\left \vert\frac{ \widetilde{\zeta}}
{ \widetilde{\zeta}_{n}  }\right \vert^{\alpha_3}      dV_n(\widetilde{\zeta}),\\
\widetilde{R}_2( \widetilde{\xi},\widetilde{\eta})  &:=&
\int\limits_{ \widetilde{\zeta}\in   \C^n,\ \frac{\delta_1}{C}  \leq  \widetilde{L}  \left(\widetilde{\xi},\widetilde{\zeta},\widetilde{\eta}        \right  )
\leq C\delta_2
}  \widetilde{K}\left(\widetilde{\xi},\widetilde{\zeta},\widetilde{\eta}        \right  )
\left \vert\frac{ \widetilde{\zeta}}
{ \widetilde{\zeta}_{n}  }\right \vert^{\alpha_3}      dV_n(\widetilde{\zeta}).
\end{eqnarray*}
where $dV_n(\widetilde{\zeta})$ denotes  the   Lebesgue measure on   $\C^n.$

Next, consider  the  following  integrals
 \begin{eqnarray*}
S(\xi)  &:=&\int\limits_{ \H_m
} R(\xi,\eta)  \left (  1 +  \frac{  \vert w \vert}{ \vert\eta \vert}
 \right )^{\alpha_2}
      \left \vert\frac{ \eta}
{ \eta_{m+1}  }\right \vert^{\alpha_4}
  dV_m(\eta),  \\
 \widetilde{S}_1( \widetilde{\xi})  &:=&
\int\limits_{  \C^m
}  \widetilde{R}_1(\widetilde{\xi}, \widetilde{\zeta}           ) \left (  1 +  \frac{ \vert\widetilde{ w} \vert}{ \vert\widetilde{\eta} \vert}
 \right )^{\alpha_2}
\left \vert\frac{ \widetilde{\eta}}
{ \widetilde{\eta}_{m}  }\right \vert^{\alpha_4}      dV_m(\widetilde{\eta}),\\
\widetilde{S}_2( \widetilde{\xi})  &:=&
\int\limits_{ \C^m
}  \widetilde{R}_2(\widetilde{\xi},  \widetilde{\zeta}           )
\left \vert\frac{ \widetilde{\eta}}
{ \widetilde{\eta}_{m}  }\right \vert^{\alpha_4}      dV_m(\widetilde{\eta}),
\end{eqnarray*}
where $dV_m(\widetilde{\eta})$ is the     Lebesgue measure on  $\C^m.$

We outline  the  main  ideas  of  the  proof. Suppose that   $\xi$ (resp. $\eta$)
is i-comparable with  $\widetilde{\xi}$ (resp. $\widetilde{\eta}$). Using  the  hypothesis that
  $K\lesssim  C\widetilde{K},$  we  shall  prove that
\begin{eqnarray}
R(\xi,\eta) \leq  C_4 \widetilde{R}_i( \widetilde{\xi},\widetilde{\eta}),\quad i=1,2.
\end{eqnarray}
Next,   we  shall  establish  in  the  same  way  as in the proof of   (5.2) the  following  estimate~: (note  that $\xi=\widetilde{\xi}$)
 \begin{eqnarray}  S(\xi)\leq  C_4  \widetilde{S}_i(\widetilde{\xi}),\quad i=1,2.
\end{eqnarray}
Finally, an application of   Fubini's theorem  gives that
\begin{eqnarray*}
\int\limits_{\Theta\in \M_N,\ \delta_1 \leq  L( \Theta)  \leq \delta_2
} K_{ \alpha    }(\Theta) dV(\Theta)
 =\int\limits_{\B_l} S(\xi)dV_l(\xi),
\end{eqnarray*}
and
\begin{eqnarray*}
\int\limits_{\widetilde{\Theta}\in   \B_{i,\vert N\vert},\frac{\delta_1}{C} \leq  \widetilde{L}(\widetilde{\Theta})
\leq C\delta_2
}  \widetilde{K}_{ i,\alpha    }(\widetilde{\Theta})      dV(\widetilde{\Theta})= \int\limits_{\B_l} \widetilde{S}_i(\xi)dV_l(\xi).                 \end{eqnarray*}
 Part 3)  now  follows  by combining  (5.3)  with  the  latter  two  estimates. It now  remains  to prove inequality   (5.2).

To do so, divide the  domain of integration  $\{ \zeta\in   \H_n:\ \frac{\delta_1}{C}  \leq  L( \xi,\zeta,\eta   )
\leq C\delta_2 \}$ of $ R(\xi,\eta)$ into the three  subsets~:
\begin{eqnarray*}
E_1 &:= & \left \{  \vert \zeta- z\vert <\frac{\vert z\vert}{C_1}             \right\};\qquad
E_2 :=  \left \{ \vert \zeta\vert >  \sqrt{2}\vert  z\vert         \right\};\\
E_3 &:= & \left \{ \vert \zeta\vert  \leq \sqrt{2}\vert  z\vert  \ \text{and}\             \vert \zeta- z\vert \geq \frac{\vert z\vert}{C_1}             \right\}
.
\end{eqnarray*}
Also, divide the domain of integration    $\left\{\widetilde{\zeta}\in\Cn:\  \frac{\delta_1}{C}  \leq  \widetilde{L} \left(\widetilde{\xi},\widetilde{\zeta},\widetilde{\eta}\right)
\leq C\delta_2
 \right\}$ of  $ \widetilde{R}_i( \widetilde{\xi},\widetilde{\eta})$   into three  corresponding subsets~:
\begin{eqnarray*}
\widetilde{E}_1 &:= & \left \{  \vert \widetilde{\zeta}-   \widetilde{z}\vert <\frac{\vert \widetilde{z}\vert}{C_2}   \right\};\qquad
 \widetilde{E}_2 :=  \left \{ \vert \widetilde{\zeta}\vert >  \sqrt{2}\vert \widetilde{ z}\vert     \right\};\\
\widetilde{E}_3 &:= & \left \{ \vert \widetilde{\zeta}\vert  \leq \sqrt{2}\vert \widetilde{ z}\vert  \ \text{and}\             \vert \widetilde{\zeta}- \widetilde{z}\vert \geq \frac{\vert \widetilde{z}\vert}{C_2}             \right\}.
\end{eqnarray*}
 Estimate (5.2) will  follow by combining  three integral estimates of  the  form
$\int\limits_{E_j}\quad \leq C_4 \int\limits_{\widetilde{E}_j}\quad$   with  some appropriate  integrands  and $j=1,2,3.$
 Therefore, we  may  assume  without loss of
 generality  that   $E_j\not =\varnothing,$ $j=1,2,3.$

 Combining  Theorem 4.1,  definition 5.3 and estimate (5.1), we  see that
  $\left(\xi,\zeta,\eta        \right  )$ is i-comparable with $\left(\widetilde{\xi},\Phi^z(\zeta),\widetilde{\eta}        \right  )\in \widetilde{E}_1,$
 for  every  $\zeta\in E_1.$ Hence, the hypothesis  $K\lesssim  C\widetilde{K}$ implies that
 $ K\left(\xi,\zeta,\eta        \right  )\leq C \widetilde{K} \left(\widetilde{\xi},\Phi^z(\zeta),\widetilde{\eta}        \right  ).$ Moreover, the fact that  $\zeta\in E_1$ gives that   $\vert \zeta\vert >\left (1 -\frac{1}{C_1}\right )\vert z\vert.$
 Therefore, applying  Theorem 4.1 (iii)-(iv)   gives that
\begin{equation}
\int\limits_{E_1
} K \left(\xi,\zeta,\eta        \right  )           \left (  1 + \left \vert \frac{z}{\zeta}
\right \vert \right )^{\alpha_1}
      \left \vert\frac{ \zeta}
{ \zeta_{n+1}  }\right \vert^{\alpha_3} dV_n(\zeta) \leq C_4
\int\limits_{ \widetilde{E}_1
}  \widetilde{K} \left(\widetilde{\xi},\Phi^z(\zeta),\widetilde{\eta}        \right  )                      \left\vert \frac{ \widetilde{\zeta}}
{ \widetilde{\zeta}_{n}  }\right \vert^{\alpha_3}     dV_n(\widetilde{\zeta}).
\end{equation}
 Next, we  prove  the estimate of  the  form  $\int\limits_{E_2}\ \leq C_4 \int\limits_{\widetilde{E}_2}\ .$
 Set $I:=\left\{ \vert \zeta\vert :\ \zeta\in E_2\right\}.$
 We  remark  that $\frac{\vert z\vert}{\vert \zeta\vert} <\frac{1}{ \sqrt{2}},$ for  every $\zeta\in E_2.$
 Therefore, by integration in  polar coordinates (Corollary 2.2), we  obtain
\begin{eqnarray*}
& &\int\limits_{E_2
} K   \left(\xi,\zeta,\eta        \right  )            \left (  1 +  \frac{\vert z\vert}{\vert \zeta\vert}
\right )^{\alpha_1}
      \left \vert\frac{ \zeta}
{ \zeta_{n+1}  }\right \vert^{\alpha_3}
  dV_n(\zeta) \\ & \lesssim &
\int\limits_{I}
\sup_{\zeta\in E_2, \vert \zeta\vert =r} K \left(\xi,\zeta,\eta        \right  ) r^{2n-1}dr\cdot  \int\limits_{\partial\M}
      \left \vert\frac{ \zeta}
{ \zeta_{n+1}  }\right \vert^{\alpha_3}
  d\sigma_n(\zeta) \\
 & \lesssim&  \int\limits_{I}
\sup_{\zeta\in E_2, \vert \zeta\vert =r} K \left(\xi,\zeta,\eta        \right  ) r^{2n-1}dr,
\end{eqnarray*}
where the  latter   inequality  holds by an application of   Lemma 4.1 in  \cite{V} with
$ \alpha_3 <2.$

On account of  definition 5.3 and  estimate (5.1),
  $\left(\xi,\zeta,\eta        \right  )$ is i-comparable with $\left(\widetilde{\xi},\widetilde{\zeta},\widetilde{\eta}        \right  )\in \widetilde{E}_2,$
 for every $\zeta\in E_2 $ and  $\widetilde{\zeta}\in\widetilde{E}_2$ such that
 $\vert \widetilde{\zeta} \vert= \vert\zeta\vert.$ This, combined with the hypothesis
   $K\lesssim  C\widetilde{K},$ implies that
   \begin{eqnarray*}
\int\limits_{I}
\sup_{\zeta\in E_2, \vert \zeta\vert =r} K \left(\xi,\zeta,\eta        \right  ) r^{2n-1}dr
   \lesssim \int\limits_{I}
\inf_{\widetilde{\zeta}\in \widetilde{E}_2, \vert \widetilde{\zeta}\vert =r} \widetilde{K}
\left(\widetilde{\xi},\widetilde{\zeta},\widetilde{\eta}        \right  )
 r^{2n-1}dr.
\end{eqnarray*}
The right side  of  the  latter  estimate is   majorized by   $ C_4 \int\limits_{\widetilde{E}_2} \widetilde{K}
 \left(\widetilde{\xi},\widetilde{\zeta},\widetilde{\eta}        \right  )dV_n(\widetilde{\zeta}).$  In summary,
 we have that
\begin{equation}
\int\limits_{E_2
} K   \left(\xi,\zeta,\eta        \right  )            \left (  1 +  \frac{\vert z\vert}{\vert \zeta\vert}
\right )^{\alpha_1}
      \left \vert\frac{ \zeta}
{ \zeta_{n+1}  }\right \vert^{\alpha_3}
  dV_n(\zeta)
\leq  C_4  \int\limits_{\widetilde{E}_2} \widetilde{K} \left(\widetilde{\xi},\widetilde{\zeta},\widetilde{\eta}        \right  )               dV_n(\widetilde{\zeta}).
\end{equation}
It now  remains to prove  the  estimate of the  form  $\int\limits_{E_3}\ \leq C_4 \int\limits_{\widetilde{E}_3}\ .$
 Consider  two  cases according to the value of $i:$\\
{\bf Case} $i=1.$  We set $R:=\sup_{\zeta\in E_3} \vert \zeta\vert.$
In view of  definition 5.3, Remark 5.2 and  estimate   (5.1), we  see that $\left(\xi,\zeta,\eta        \right  )$ is 1-comparable with $\left(\widetilde{\xi},\widetilde{\zeta},\widetilde{\eta}        \right  )\in \widetilde{E}_3$  for  every $\zeta\in E_
3$ and $\widetilde{\zeta}\in \Cn$ such that $\vert \widetilde{\zeta} - \widetilde{z} \vert  \geq  \frac{\widetilde{z}    }{C_2} $  and $\vert\widetilde{\zeta}\vert
 \leq  \vert\zeta\vert.$
   Therefore,  using  integration in polar  coordinates, we  obtain
\begin{equation}
\begin{split}
& \int\limits_{E_3
} K  \left(\xi,\zeta,\eta        \right  )
\left (  1 +  \frac{\vert z\vert}{\vert\zeta\vert}
 \right )^{\alpha_1}
      \left \vert\frac{ \zeta}
{ \zeta_{n+1}  }\right \vert^{\alpha_3}
  dV_n(\zeta) \\ & \lesssim
\int\limits_{0}^{R}
\sup_{\zeta\in E_3, \vert \zeta\vert =r} K \left(\xi,\zeta,\eta        \right  )              r^{2n-1} \left (  1 +  \frac{\vert z\vert}{r}
 \right )^{\alpha_1}  \cdot  \int\limits_{\partial\M}
                        \left \vert\frac{ \zeta}
{ \zeta_{n+1}  }\right \vert^{\alpha_3}
  d\sigma_n(\zeta)dr \\
 & \lesssim
\int\limits_{0}^{R}
\inf_{ \widetilde{  \zeta}\in \widetilde{E}_3, \vert \widetilde{\zeta}\vert =r} \widetilde{K}
\left(\widetilde{\xi},\widetilde{\zeta},\widetilde{\eta}        \right  )r^{2n-1}
\left (  1 + \frac{\vert \widetilde{z}\vert}{r}
 \right )^{\alpha_1}
\cdot  \int\limits_{\partial\B_n}
                           \left \vert\frac{ \widetilde{\zeta}}
{ \widetilde{\zeta}_{n}  }\right \vert^{\alpha_3}
  d\sigma_n(\widetilde{\zeta})dr    \\
& \leq  C_4
\int\limits_{ \widetilde{E}_3  }
 \widetilde{K}  \left(\widetilde{\xi},\widetilde{\zeta},\widetilde{\eta}        \right  )
  \left (  1 +  \frac{\vert \widetilde{z}\vert}{\vert \widetilde{\zeta} \vert}
 \right )^{\alpha_1}                         \left \vert\frac{ \widetilde{\zeta}}
{ \widetilde{\zeta}_{n}  }\right \vert^{\alpha_3}
  dV_n(\widetilde{\zeta}),
\end{split}
\end {equation}
where on the third line, $d\sigma_n(\widetilde{\zeta})$ is the surface measure of  the euclidean unit  sphere
 $\partial\B_n$  of $\Cn.$
\\
{\bf Case} $i=2.$  We  see easily  that
$$ \widetilde{E}_3=
\left\{\widetilde{\zeta}  \in\Cn:\ \vert \widetilde{\zeta}\vert  \leq \sqrt{2}\vert \widetilde{ z}\vert  \ \text{and}\             \vert \widetilde{\zeta}- \widetilde{z}\vert \geq \frac{\vert \widetilde{z}\vert}{C_2}             \right\} .$$
Moreover,  $\left(\xi,\zeta,\eta        \right  )$ is 2-comparable with $\left(\widetilde{\xi},\widetilde{\zeta},\widetilde{\eta}        \right  )$
for  every   $\zeta\in  E_3$  and     $\widetilde{\zeta}\in  \widetilde{E}_3.
$    On the other hand, by  Remark  5.2, we have  $\vert z\vert=\vert\widetilde{z}\vert.$  Thus,
\begin{equation}
\begin{split}
& \int\limits_{E_3
} K \left(\xi,\zeta,\eta        \right  )           \left (  1 +  \frac{\vert z\vert}{\vert\zeta\vert}
 \right )^{\alpha_1}
      \left \vert\frac{ \zeta}
{ \zeta_{n+1}  }\right \vert^{\alpha_3}
  dV_n(\zeta) \\ & \lesssim \sup_{\zeta\in E_3} K\left(\xi,\zeta,\eta        \right  )
\int\limits_{  \vert \zeta\vert \leq \sqrt{2}\vert z\vert}
  \left (  1 + \frac{\vert z\vert}{\vert \zeta\vert}
\right )^{\alpha_1}                         \left \vert\frac{ \zeta}
{ \zeta_{n+1}  }\right \vert^{\alpha_3}
  dV_n(\zeta)\\  &\lesssim
\vert z\vert^{2n}
\sup_{  \zeta\in E_3} K\left(\xi,\zeta,\eta        \right  )
 \lesssim
\vert \widetilde{z}\vert^{2n}
\inf_{ \widetilde{  \zeta}\in \widetilde{E}_3}  \widetilde{K}\left(\widetilde{\xi},\widetilde{\zeta},\widetilde{\eta}        \right  )
\leq  C_4  \int\limits_{\widetilde{E}_3}  \widetilde{K}\left(\widetilde{\xi},\widetilde{\zeta},\widetilde{\eta}        \right  )
dV_n(\widetilde{\zeta}).
\end{split}
\end {equation}

 Now estimate (5.2)   follows from    (5.4)-(5.7).
 This  completes the proof of part 3).
\end{proof}
To conclude  this section, we  give without proof  some  properties  of  the   relations  $"\lesssim"$
and $"\approx".$
\begin{prop}
Let  $Z,\widetilde{Z}$ and $K,L, \widetilde{K}, \widetilde{L}$  be     as in the
statement of      Theorem 5.6.
Suppose that
 $   K\lesssim  \widetilde{K} \quad\text{and}\quad L\lesssim \widetilde{L} \quad\text{at}\ (Z, \widetilde{Z}).$
Then  $K+L \lesssim \widetilde{K}+ \widetilde{L}$  and
          $K^{\alpha} L^{\beta}\lesssim
\widetilde{K}^{\alpha} \widetilde{L}^{\beta},$ for  every $\alpha,\beta \geq 0.$

If in addition    $\widetilde{K} \approx K$  and  $\widetilde{L} \approx L$ then $K+L\approx \widetilde{K}+
\widetilde{L}$ and
$K^{\alpha} L^{\beta}  \approx
\widetilde{K}^{\alpha} \widetilde{L}^{\beta},$ for  every $\alpha,\beta\in\R.$
 \end{prop}

\section{Integral  kernels}
  The  pairs of  integral  kernels $K,\widetilde{K}$ satisfying the condition $  K\approx \widetilde{K} $ that we  shall use
 are  studied here.
Recall the  function $D$  introduced  by Charpentier \cite{Ch}~:
\begin{eqnarray*}
D(\Theta,Z):=  \vert 1- \Theta\bullet \overline{Z}\vert^2 - (1-\vert \Theta\vert^2)(1-\vert Z\vert^2), \qquad\text{for all}\
\Theta,Z\in \C^{k}\ \text{and}\ k\in \N.
\end{eqnarray*}
\begin{thm}  Let  $i\in\{1,2\}$ and fix two    comparable  points  $Z\in \B_N$    and $\widetilde{Z} \in \B_{\vert N\vert}.$
  Consider two  functions     $ K,\widetilde{K}$   defined respectively on    $\overline{\M}_N$ and  $\overline{\B}_{i,\vert N\vert}$ that  correspond to  one  of the    following three cases~:

 (1) $i=2$ and $ K(\Theta):= \vert \Theta - Z\vert,\quad \widetilde{K}(\widetilde{\Theta}):= \vert \widetilde{\Theta} - \widetilde{Z}\vert.$

 (2) $i=1$ and $ K(\Theta):= \vert  1- \Theta\bullet\overline{Z}\vert,\quad\widetilde{K}(\widetilde{\Theta}):= \vert 1-\widetilde{\Theta} \bullet \overline{\widetilde{Z}} \vert.$

 (3) $i=1$ and  $ K(\Theta):=D(\Theta,Z),\quad\widetilde{K}(\widetilde{\Theta}):=D(\widetilde{\Theta},\widetilde{Z}).$

 Then $K \approx \widetilde{K} $ at $(Z,\widetilde{Z}).$
\end{thm}
\begin{proof}
Using the  definitions  5.1, 5.3  and  5.5, it can be  easily   checked that
  $$\vert z -\zeta\vert \approx \vert \widetilde{z} -\widetilde{\zeta}\vert \quad \text{and}\quad
    \vert w -\eta\vert \approx \vert \widetilde{w} -\widetilde{\eta}\vert \quad \text{at}\ (Z,\widetilde{Z}) .$$ Applying
     Proposition 5.7 to the latter  two  relations,  assertion (1) follows.

    To prove   assertions (2) and (3), we  need the following  estimates   of
    Bonami-Charpentier \cite[p. 67]{BCh}~:
   \begin{equation}
  \vert  1- \Theta\bullet\overline{Z}\vert  \approx  (1-\vert Z\vert^2) + ( 1- \vert \Theta\vert^2   )  +\left\vert  \text{Im}\ \Theta\bullet \overline{Z} \right\vert +\vert  \Theta -Z\vert^{2},
     \end{equation}
 and
   \begin{equation}
    D(\Theta,Z)\approx  (1-\vert Z\vert^{2}) \vert  \Theta -Z\vert^{2}+\left (  \vert \Theta\vert^2 -\vert Z\vert^2 \right )^{2}  +\left\vert  \text{Im}\ \Theta\bullet \overline{Z} \right\vert^{2} +\vert   \Theta -Z\vert^{4} ,
     \end{equation}
for  every  $\Theta, Z\in\B_k,$ where  $\B_k$ is as usual  the  euclidean  ball of $\C^k.$

 Write $Z\equiv   (x,z,w)\in \B_N$    and $\widetilde{Z} \equiv   (\tilde{x},\tilde{z},\tilde{w})\in \B_{\vert N\vert}.$
 Let $\Theta\equiv (\xi,\zeta,\eta)\in\overline{\M}_N$ be  1-comparable with
$\widetilde{\Theta}\equiv \left(\widetilde{\xi},\widetilde{\zeta},\widetilde{\eta}\right)\in \B_{\vert N\vert}.$   We break the  proof into    four   cases.
\\
{\bf Case 1:}
 $\widetilde{z}=\Phi^{z}(z),  \widetilde{\zeta}=\Phi^{z}(\zeta)$ and          $\widetilde{w}=\Phi^{w}(w),$  $\widetilde{\eta}=\Phi^{w}(\eta).$

 In this case by Theorem 4.1 (i)-(ii), we  have that
 \begin{eqnarray}
 \zeta\bullet\overline{z}  =\widetilde{\zeta}\bullet\overline{\widetilde{z}},\
\eta\bullet\overline{w}  =\widetilde{\eta}\bullet\overline{\widetilde{w}},\quad\text{and}\
\vert z\vert =\vert\overline{z}\vert,\ \vert w\vert =\vert\overline{w}\vert.
\end{eqnarray}
We deduce  easily  from the  first two  equalities   of (6.3) that $\vert  1- \Theta\bullet\overline{Z}\vert= \vert 1-\widetilde{\Theta} \bullet \overline{\widetilde{Z}} \vert,$    which proves assertion  (2).

On the  other hand, we have  the following identity~:
 \begin{eqnarray*}
D(\Theta,Z)=  (1-\vert Z\vert^{2}) \vert  \Theta -Z\vert^{2}   +\left\vert \overline{Z}  \bullet (\Theta-Z)  \right\vert^{2}.  \end{eqnarray*}
 By  assertion (1), we have  $\vert  \Theta -Z\vert^{2}\approx  \vert \widetilde{ \Theta} -\widetilde{Z}\vert^{2}.$
This, combined with (6.3),    implies that      $ D(\Theta,Z)\approx  D\left(\widetilde{\Theta},\widetilde{Z}\right),$ which
  completes the  proof of assertion (3).
\\
{\bf Case 2:}
 $\widetilde{z}=\Phi^{z}(z),  \widetilde{\zeta}=\Phi^{z}(\zeta)$ and         $\vert w-\eta\vert  \geq \frac{ \vert w\vert}{C_1},  \vert  \widetilde{w}-\widetilde{\eta}\vert\geq \frac{ \vert \widetilde{w}\vert}{C_2} .$

 In this case  it is easy  to  check that
\begin{eqnarray}
\max\left\{\vert  w\vert,\vert  \eta\vert\right\}\lesssim  \vert w-\eta\vert,\quad\text{and}\
\left\vert  \text{Im}\ \eta\bullet \overline{w} \right\vert  \lesssim  \vert w-\eta\vert^{2}.
\end{eqnarray}
Now  we set
$$ Z^{'}:=(x,z),\  \widetilde{Z}^{'}:=(\widetilde{x},\widetilde{z}),\ \Theta^{'}:=(\xi,\zeta), \widetilde{\Theta}^{'}:=(\widetilde{\xi},\widetilde{\zeta}).$$
Combining (6.1),(6.3) and (6.4), we  obtain
 \begin{equation}\begin{split}
  \vert  1- \Theta\bullet\overline{Z}\vert & \approx  \vert \eta -w\vert^{2}+ (1-\vert Z^{'}\vert^2) + ( 1- \vert \Theta^{'}\vert^2   )  +\left\vert  \text{Im}\ \Theta^{'}\bullet \overline{Z^{'}} \right\vert +\vert  \Theta^{'} -Z^{'}\vert^{2}\\
  &\approx \vert \eta -w\vert^{2}+ \vert  1- \Theta^{'}\bullet\overline{Z^{'}}\vert.
   \end{split}  \end{equation}
     On the  one hand, we  have $\vert \eta-w\vert^{2}\approx   \vert \widetilde{\eta}-\widetilde{w}\vert^{2}.$ On the other
     hand,
    proceeding as in the first case, we get  $\vert  1- \Theta^{'}\bullet\overline{Z^{'}}\vert= \vert 1-\widetilde{\Theta^{'}} \bullet \overline{\widetilde{Z^{'}}} \vert.$
   This, combined with  (6.5),  shows that   $\vert  1- \Theta\bullet\overline{Z}\vert \approx \vert 1-\widetilde{\Theta} \bullet \overline{\widetilde{Z}} \vert,$  which completes the  proof of  assertion (2).

We  now  come to  the proof of assertion (3).   Applying (6.2),(6.3) and (6.4), we  see easily that
\begin{equation*}
    D(\Theta,Z)\approx  (1-\vert Z\vert^{2}) \vert  \Theta -Z\vert^{2}+\left (  \vert \Theta^{'}\vert^2 -\vert Z^{'}\vert^2 \right )^{2}  +\left\vert  \text{Im}\ \Theta^{'}\bullet \overline{Z^{'}} \right\vert^{2} +\vert   \Theta^{'} -Z^{'}\vert^{4}
    +\vert \eta -w\vert^{4}.
     \end{equation*}
     Since $ (1-\vert Z\vert^{2}) \vert  \Theta -Z\vert^{2}= \left (1-\vert Z^{'}\vert^{2}-\vert w\vert^{2}\right) \left(\vert
      \Theta^{'} -Z^{'}\vert^{2}+\vert
      \eta -w\vert^{2}\right ),$  we obtain
\begin{equation}\begin{split}
    D(\Theta,Z)&\approx  \left\lbrack (1-\vert Z^{'}\vert^{2}) \vert  \Theta^{'} -Z^{'}\vert^{2}+\left (  \vert \Theta^{'}\vert^2 -\vert Z^{'}\vert^2 \right )^{2}  +\left\vert  \text{Im}\ \Theta^{'}\bullet \overline{Z^{'}} \right\vert^{2} +\vert   \Theta^
{'} -Z^{'}\vert^{4}\right\rbrack\\
    &+\vert \eta -w\vert^{2}\left ( 1-    \vert Z^{'}\vert^2  \right ) \\
    &\approx  D(\Theta^{'},Z^{'})+\vert \eta -w\vert^{2}\left ( 1-  \vert Z^{'}\vert^2         \right ).
    \end{split} \end{equation}
   On the one hand, we  have $\vert \eta-w\vert^{2}\approx   \vert \widetilde{\eta}-\widetilde{w}\vert^{2}$ and
   $1-\vert Z^{'}\vert^2=1-\vert \widetilde{Z}^{'}\vert^2.$   On the other  hand,
    proceeding as in  the  first case, we  can show  that  $ D(\Theta^{'},Z^{'})\approx D\left(\widetilde{\Theta}^{'}, \widetilde{Z}^{'}\right ) .$
   This,  combined with (6.6), shows that $ D(\Theta,Z)\approx D\left(\widetilde{\Theta}, \widetilde{Z}\right ) $ and
   the  proof of assertion (3) is thereby completed.
    \\
    {\bf  Case 3:}
 $\widetilde{w}=\Phi^{w}(w),  \widetilde{\eta}=\Phi^{w}(\eta)$ and
         $\vert z-\zeta\vert  \geq \frac{ \vert z\vert}{C_1},
         \vert  \widetilde{z}-\widetilde{\zeta}\vert\geq \frac{ \vert \widetilde{z}\vert}{C_2} .$

 This  case  can be  treated analogously  as  the previous case.
 \\
 {\bf  Case  4:}
   $\vert z-\zeta\vert  \geq \frac{ \vert z\vert}{C_1},  \vert  \widetilde{z}-\widetilde{\zeta}\vert\geq
    \frac{ \vert \widetilde{z}\vert}{C_2} $  and   $\vert w-\eta\vert  \geq \frac{ \vert w\vert}{C_1},
      \vert  \widetilde{w}-\widetilde{\eta}\vert\geq \frac{ \vert \widetilde{w}\vert}{C_2} .$

We repeat the arguments used in  the proof of  the second  case. More precisely, proceeding as in the proof of  (6.5) and (6.6), one can show that
\begin{equation}\begin{split} \vert  1- \Theta^{'}\bullet\overline{Z^{'}}\vert
  &\approx \vert \zeta -z\vert^{2}+ \vert  1- \xi\bullet\overline{x}\vert\\
   D(\Theta^{'},Z^{'})  &\approx  D(\xi,x)+\vert \zeta -z\vert^{2}\left ( 1-\vert x \vert^{2}  \right ).\end{split}
    \end{equation}
    On the other hand,  it is clear that
    $\vert 1-  \xi\bullet\overline{ x     } \vert =\vert 1-  \widetilde{\xi}\bullet\overline{\widetilde{ x}     } \vert     $
    and  $D(\xi,x)= D(\widetilde{\xi},\widetilde{x}).$

    These  equalities, combined with  (6.5), (6.6)  and (6.7), imply that
        $\vert  1- \Theta\bullet\overline{Z}\vert \approx \vert 1-\widetilde{\Theta} \bullet \overline{\widetilde{Z}} \vert$  and $ D(\Theta,Z)\approx D\left(\widetilde{\Theta}, \widetilde{Z}\right ) .$   The proof of the theorem is   complete in this  la
st  case.
\end{proof}
\section{Integral  estimates}
 In this section, we  prove, with the help of  Theorem 5.6, two important  integral  estimates  that will   be  used  repeatedly throughout  the paper.

For  every  $\lambda >0$ and $ Z\in \B_{N},$ consider the  function
\begin{eqnarray*}
K_{\lambda, Z}(\Theta):= \frac{1}{ \left\vert 1 - \Theta \bullet \overline{Z}  \right \vert^{  \vert N\vert +1+\lambda}},\qquad
\text{for all}\  \Theta\in\B_{N}.
\end{eqnarray*}
The  first   result of  this section  is the following
\begin{thm}  For  every  $ \alpha:= \left ( \alpha_{1},\alpha_{2},\alpha_{3},\alpha_{4}  \right) \in \R^4$ such that
$0\leq \alpha_3, \alpha_4 < 2$ and
 $0\leq \alpha_1 +\alpha_3 < 2n,$  $0\leq \alpha_2 +\alpha_4 < 2m,$              there exists  a constant $C$ independant of $Z\in\B_N$  such that
 \begin{eqnarray*} \int\limits_{\M_N}  K_{\lambda, Z}(\Theta)
  \left (  1 +  \frac{\vert z\vert }{\vert \zeta\vert }
 \right )^{\alpha_1} \left (  1 +  \frac{\vert w\vert }{\vert \eta\vert }
 \right )^{\alpha_2}
      \left \vert\frac{ \zeta}
{ \zeta_{n+1}  }\right \vert^{\alpha_3}
 \left \vert\frac{ \eta}
{ \eta_{m+1}  }\right\vert^{\alpha_4}dV(\Theta)\leq  C\left (1-\vert Z\vert^2\right )^{-\left (
\lambda +\frac{\alpha_3}{2}+ \frac{\alpha_4}{2}\right )}.
\end{eqnarray*}
\end{thm}
 In order  to state  the second result  of this  section, we introduce some more notations.
Let $\Theta :=(\xi,\zeta,\eta),\ Z:=(x,z,w)$ and $Z^{'}:=(x^{'},z^{'},w^{'}) $  be  three points of $\M_N.$ Define
$$
\Delta(\Theta, Z,Z^{'}):=\frac{1}{\vert\zeta\vert^2 \vert\eta\vert^2}
\sum\limits_{j=1}^{\vert N\vert +2}
\left \vert \frac{B_j(\Theta,Z)}{ \vert \Theta - Z\vert^{2\vert N\vert}}
- \frac{B_j(\Theta,Z^{'})}{ \vert \Theta - Z^{'}\vert^{2\vert N\vert}}
\right \vert,$$
where  $B_j(\Theta,Z)$ are  the    polynomials given  in the statement of   Theorem 3.1.
\begin{thm}
Let  $q$ be  a  real number such that  $ 1\leq q <\frac{ 2\vert N\vert+4}{ 2\vert N\vert+3}.$ Then we have the estimate
\begin{eqnarray*}
\int\limits_{\M_N} \Delta(\Theta, Z,Z^{'})^q \left (\frac{\vert \zeta\vert \vert \eta\vert}
{\vert \zeta_{n+1}\vert \vert \eta_{m+1}\vert }\right )^{ 2(q-1)} dV(\Theta)
\leq   \left\lbrace \begin{array}{l}
             \vert Z -Z^{'} \vert^{ 2\vert N\vert +4 - ( 2\vert N\vert +3)q},\ \text{if}\ q >1;\\
             \vert Z -Z^{'} \vert \vert \log {\vert Z -Z^{'} \vert }\vert,\ \text{if}\ q=1.
             \end{array} \right.
\end{eqnarray*}
\end{thm}
To prove  these  theorems, we  need   some  preparatory  lemmas.
\begin{lem}
Given $0< R_1\leq  R_2,$   $\alpha <1$ and  $0\leq\beta,\gamma,$  then
$$
\int\limits_{0}^{R_1} \frac{dx}{x^{\alpha} (x+R_1)^{\beta}(x+R_2)^{\gamma}}
\leq C(\alpha,\beta,\gamma) \int\limits_{0}^{R_1} \frac{ dx}{ (x+R_1)^{\alpha +\beta}(x+R_2)^{\gamma}}.
$$
\end{lem}
\begin{lem}  Consider $\alpha:=\left ( \alpha_1,\alpha_2,\alpha_3,\alpha_4,\alpha_5 \right ) \in\R^{5}$   such that
$0\leq \alpha_4,\alpha_5 <2,$  and $\alpha_2 < 2n,$   $\alpha_3  < 2m.$
For $a:=(a_1,a_2)\in\C^{2},$
we  set
\begin{eqnarray*}
I_{\alpha,a}( \widetilde{ \Theta}):= \frac{1}{ \vert \widetilde{ \Theta}\vert^{\alpha_1}   \vert \widetilde{ \zeta}\vert^{\alpha_2} \vert\widetilde{  \eta}\vert^{\alpha_3} }
\left \vert \frac{\widetilde{ \zeta}}
{ \widetilde{ \zeta}_{n}-a_1 }\right \vert^{ \alpha_4}\left \vert \frac{\widetilde{ \eta}}
{ \widetilde{ \eta}_{m}-a_2 }\right \vert^{ \alpha_5},\quad \widetilde{ \Theta}\equiv (
 \widetilde{\xi},\widetilde{\zeta},\widetilde{\eta})     \in \C^{\vert N\vert}.\end{eqnarray*}
Then
\begin{eqnarray*}
\int\limits_{ \vert \widetilde{ \Theta}\vert < \delta}  I_{\alpha,a}(\widetilde{ \Theta}) dV(\widetilde{ \Theta}) \leq  C(N,\alpha)\delta^{2\vert N\vert -(\alpha_1+\alpha_2 +\alpha_3)},\  \text{if}\ \alpha_1+\alpha_2 +\alpha_3 <
2\vert N\vert;\\
\int\limits_{\delta< \vert \widetilde{ \Theta}\vert }  I_{\alpha,a}(\widetilde{ \Theta}) dV(\widetilde{ \Theta}) \leq  C(N,\alpha)\delta^{2\vert N\vert -(\alpha_1+\alpha_2 +\alpha_3)},\  \text{if}\ \alpha_1+\alpha_2 +\alpha_3 >
2\vert N\vert.
\end{eqnarray*}
\end{lem}
\begin{proof}
 Consider the  case $a=(0,0).$   We use  integration in polar  coordinates and  then    apply Lemma 7.3  four times to  obtain
\begin{eqnarray*}
\int\limits_{ \vert \widetilde{ \Theta}\vert < \delta}  I_{\alpha,0}(\widetilde{ \Theta}) dV(\widetilde{ \Theta}) \lesssim
\int\limits_{0< \vert \widetilde{ \Theta}\vert < \delta}\frac{1}{ \vert \widetilde{ \Theta}\vert^{\alpha_1+
 \alpha_2 +\alpha_3 } } dV(\widetilde{ \Theta})\leq  C(N,\alpha)\delta^{2\vert N\vert -(\alpha_1+\alpha_2 +\alpha_3)},
\end{eqnarray*}
which  is the first estimate in the lemma.
The  second estimate  can be proved  in the same  way.
Now  we consider the  general case $a\in \C^2.$  We  write  $\widetilde{ \Theta}\equiv \left(\widetilde{ \Theta}^{'},
\widetilde{ \eta}_{m}\right )$  and observe  that if   $ \vert \widetilde{ \Theta}\vert < \delta$ and $
 \vert \widetilde{ \eta}_{m}-a_2 \vert <\frac{\vert a_2\vert}{2},$      then we have
$\left\vert (\widetilde{ \Theta}^{'},
 t )\right  \vert  <  3\delta,$  for  every $t\in\C$  such that $ \vert t- a_2\vert <\frac{\vert a_2\vert}{2}.$ In addition,
 it is clear that
 $$
 \int\limits_{\vert t- a_2\vert <\frac{\vert a_2\vert}{2}} \frac{dt\wedge d\overline{t}}{ \vert t- a_2\vert^{\alpha_5} }
 \leq  C(\alpha_5) \int\limits_{\vert t- a_2\vert <\frac{\vert a_2\vert}{2}} \frac{dt\wedge d\overline{t}}
 {\vert t\vert^{\alpha_5}}.$$
 On the other hand,  if $  \vert \widetilde{ \eta}_{m}-a_2 \vert \geq \frac{\vert a_2\vert}{2},$ then
 $\vert \widetilde{ \eta}_{m} \vert  \leq 3  \vert \widetilde{ \eta}_{m}-a_2 \vert.$
 It follows from these considerations   that
\begin{eqnarray*}
\int\limits_{ \vert \widetilde{ \Theta}\vert < \delta}  I_{\alpha,a}(\widetilde{ \Theta}) dV(\widetilde{ \Theta})
\leq C(\alpha_5) \int\limits_{ \vert \widetilde{ \Theta}\vert < 3\delta}
\frac{1}{ \vert \widetilde{ \Theta}\vert^{\alpha_1}   \vert \widetilde{ \zeta}\vert^{\alpha_2} \vert\widetilde{  \eta}\vert^{\alpha_3} }
\left \vert \frac{\widetilde{ \zeta}}
{ \widetilde{ \zeta}_{n}-a_1 }\right \vert^{ \alpha_4}\left \vert \frac{\widetilde{ \eta}}
{ \widetilde{ \eta}_{m} }\right \vert^{ \alpha_5}
.\end{eqnarray*}
The same reasoning, applied to    the  variables $\widetilde{ \zeta}$ and $\widetilde{ \zeta}_n,$ shows that
\begin{eqnarray*}
\int\limits_{ \vert \widetilde{ \Theta}\vert < \delta}  I_{\alpha,a}(\widetilde{ \Theta}) dV(\widetilde{ \Theta})
\leq   C(\alpha_4,\alpha_5) \int\limits_{ \vert \widetilde{ \Theta}\vert < 9\delta} I_{\alpha,0}(\widetilde{ \Theta}) dV(\widetilde{ \Theta}).
\end{eqnarray*}
The proof of the first estimate is now complete. The second estimate can be established in a similar way.
\end{proof}

 \noindent{\bf  Proof of   Theorem 7.1.}
Let $\widetilde{Z}\equiv  \left ( \widetilde{x},\widetilde{z},\widetilde{w}\right )        $ be a point of $\B_{\vert N\vert}$ which is  comparable with  $Z.$   Consider the function
\begin{eqnarray*}
\widetilde{K}_{\lambda, \widetilde{Z}}(\widetilde{\Theta}):= \frac{1}{ \left\vert 1 - \widetilde{\Theta} \bullet \overline{\widetilde{Z}}  \right \vert^{  \vert N\vert +1+\lambda}},\qquad \text{for  all}\ \widetilde{\Theta}\equiv
\left ( \widetilde{\xi}, \widetilde{\zeta}, \widetilde{\eta}\right )            \in\B_{\vert N\vert}.
\end{eqnarray*}
By  Theorem  6.1(2) we have  $ K_{\lambda, Z}\approx \widetilde{K}_{\lambda, \widetilde{Z}}.$
By Remark 5.2 we obtain  $\vert Z\vert= \vert \widetilde{Z}\vert.$   Therefore, in view of    Theorem  5.6 we see that   Theorem  7.1  will follow from the  estimate
\begin{equation}
  \int\limits_{\B_{\vert N\vert}}  \widetilde{K}_{\lambda, \widetilde{Z}}(\widetilde{\Theta})
  \left (  1 +  \frac{\vert \widetilde{z}\vert }{\vert \widetilde{\zeta}\vert }
 \right )^{\alpha_1} \left (  1 +  \frac{\vert \widetilde{w}\vert }{\vert \widetilde{\eta}\vert }
 \right )^{\alpha_2}
      \left \vert\frac{ \widetilde{\zeta}}
{\widetilde{ \zeta}_{n}  }\right \vert^{\alpha_3}
 \left \vert\frac{ \widetilde{\eta}}
{ \widetilde{\eta}_{m}  }\right\vert^{\alpha_4}dV(\widetilde{\Theta})  \leq  C\left (1-\vert \widetilde{Z}\vert^2\right )^{-\left (
\lambda + \frac{\alpha_3}{2}+ \frac{\alpha_4}{2}\right )}.
\end{equation}
 Now we go back to  the proof of   Lemma I.5 in the work of Bonami-Charpentier \cite[p. 68-69]{BCh}. We may  assume  without loss of  generality that
 $\vert \widetilde{w}_{1}\vert \geq  \frac{1}{2\sqrt{ \vert N\vert } }$   and   set
 $$  \widetilde{w}^{'}:= \left (  \widetilde{w}_{2},\ldots, \widetilde{w}_{m}\right )\quad\text{and}\quad
 A:= 1 -\vert \widetilde{Z}\vert^{2} .$$
 As in  \cite[p. 68]{BCh},  observe  that
 $$  u:= 1-\vert\widetilde{ \Theta}\vert^{2},\ v:= \text{Im}( \widetilde{\Theta}\bullet\overline{\widetilde{Z}}),\ \widetilde{ \xi}, \widetilde{\zeta}\quad \text{and}
 \quad
  \widetilde{\eta}^{'}:= \left (  \widetilde{\eta}_{2},\ldots, \widetilde{\eta}_{m}\right )$$
 form a system of  coordinates whose  jacobian is bounded from  above and from below by positive  constants uniformly in  $\widetilde{Z}\in \B_{\vert N\vert}$ and
 that satisfies  $ \vert\widetilde{Z}-\widetilde{\Theta}\vert  \leq  \frac{1}{4\sqrt{ \vert N\vert } }.  $
Using the following  estimate    (see  \cite[p. 68]{BCh})~:
$$  \vert 1-\widetilde{\Theta}\bullet\overline{\widetilde{Z}} \vert \approx  A+\vert u\vert +\vert v\vert +
\vert\widetilde{x}-\widetilde{\xi}\vert^{2}+\vert\widetilde{z}-\widetilde{\zeta}\vert^{2}+
 \vert\widetilde{w}^{'}-\widetilde{\eta}^{'}\vert^{2},$$
we see that  in order to prove  (7.1), it suffices to establish
 \begin{eqnarray*}
& & A^{\lambda+ \frac{ \alpha_{3}}{2} +\frac{ \alpha_{4}}{2}  }\int\limits_{\C^{\vert N\vert }}
 \frac{
 \left (  1 +  \frac{\vert \widetilde{z}\vert }{\vert \widetilde{\zeta}\vert }
 \right )^{\alpha_1}
      \frac{ 1}
{\vert \widetilde{ \zeta}_{n} \vert^{\alpha_3} }
 \frac{ 1}
{ \vert \widetilde{ \eta}_{m}\vert^{\alpha_4}}    \cdot dudv dV(\widetilde{\xi})dV(\widetilde{\zeta})
dV(\widetilde{\eta}^{'})
}
{ \left (A+\vert u\vert +\vert v\vert +  \vert\widetilde{x}-\widetilde{\xi}\vert^{2}+\vert\widetilde{z}-\widetilde{\zeta}\vert^{2}+
 \vert\widetilde{w}^{'}-\widetilde{\eta}^{'}\vert^{2}
\right )^{
\vert N\vert +1+\lambda}}\\
&\leq&  C(N,\alpha,\lambda)    \int\limits_{0}^{\infty}  \frac{du}{ ( 1+u)^{1+\frac{\lambda}{5}}}
 \int\limits_{-\infty}^{\infty}  \frac{dv}{ ( 1+\vert v\vert)^{1+\frac{\lambda}{5}}}
  \int\limits_{\C^{l}}  \frac{ dV(\widetilde{\xi})
  }
  {  ( 1+ \vert\widetilde{x}-\widetilde{\xi}\vert^{2})^{
l+\frac{\lambda}{5} }} \\
&      &
  \int\limits_{\C^{n}}  \frac{
\left (  1 +  \frac{\vert \widetilde{z}\vert }{\vert \widetilde{\zeta}\vert }
 \right )^{\alpha_1}
      \frac{ 1}
{\vert \widetilde{ \zeta}_{n} \vert^{\alpha_3}    }
  }
  {  ( 1+ \vert\widetilde{z}-\widetilde{\zeta}\vert^{2} )^{
n+\frac{\lambda}{5} }}
dV(\widetilde{\zeta})
\int\limits_{\C^{m -1}}  \frac{
 \frac{ 1}
{ \vert \widetilde{ \eta}_{m}\vert^{\alpha_4}}
  }
  { ( 1+ \vert\widetilde{w}^{'}-\widetilde{\eta}^{'}\vert^{2} )^{
m -1+\frac{\lambda}{5} }}
dV(\widetilde{\eta}^{'}).
\end{eqnarray*}
 To finish the proof we only need show that
 the last  two   integrals in the last line  are finite.
 But this  will  follow   from the following
 \begin{lem}
 For  every  real numbers  $\alpha,\beta,\gamma$  such that $0<\gamma ,$ $0\leq \beta <2$ and
 $0\leq \alpha +\beta < 2n,$  there  is a  constant $C:= C(n,\alpha,\beta,\gamma)$  such that
 \begin{eqnarray*}
 I:=\int\limits_{\C^{n}}  \frac{
\left (  1 +  \frac{\vert z\vert }{\vert \zeta\vert }
 \right )^{\alpha}
      \frac{ 1}
{\vert \zeta_{n} \vert^{\beta}  }
  }
  {  ( 1+ \vert z-\zeta\vert^{2} )^{
n+\gamma }}
dV(\zeta) < C,\qquad \text{for all}  \ z\in\C^{n},
\end{eqnarray*}
where $dV$ is the   Lebesgue measure of  $\Cn.$
 \end{lem}
{\it Proof.}   Dividing the domain of  integration $\C^n$  of  $I$  into the
 three subsets  $\left\{ \vert \zeta \vert <  \frac{\vert z\vert }{2} \right\},$
   $\left\{ \vert \zeta \vert >  2\vert z\vert  \right\}$  and  $\left\{  \frac{\vert z\vert }{2}\leq \vert \zeta \vert \leq
   2\vert z\vert \right\},$ we  thus  divide  $I$ into  three corresponding  terms    $I_1,I_2$ and  $I_3.$
We now estimate   each of  these  terms.  On the one hand, we have
 \begin{eqnarray*}
 I_1 \lesssim  \int\limits_{ \left\{ \vert \zeta \vert <  \frac{\vert z\vert }{2} \right\}          }  \frac{
  \frac{\vert z\vert^{\alpha} }{\vert \zeta\vert^{\alpha} }
  \frac{ 1}
{\vert \zeta_{n} \vert^{\beta}  }
  }
  {  ( 1+ \vert z\vert^{2} )^{
n+\gamma }}
dV(\zeta) \lesssim   \frac{
\vert z \vert^{\alpha}
  }
  {  ( 1+ \vert z\vert^{2} )^{
n+\gamma } } \cdot \vert z \vert^{2n-\alpha-\beta}  < C,
\end{eqnarray*}
 where  the  second inequality  holds by  applying  a variant of   Lemma 7.4. On the other hand,
 \begin{eqnarray*}
 I_2 \lesssim  \int\limits_{ \left\{ \vert \zeta \vert >  2\vert z\vert  \right\}                       }  \frac{
  \frac{ 1}
{\vert \zeta_{n} \vert^{\beta}  }
  }
  {  ( 1+ \vert \zeta\vert^{2} )^{
n+\gamma }}
dV(\zeta) \lesssim
\int\limits_{\vert \zeta\vert < 1} \frac{
dV(\zeta)
  }
  {  \vert \zeta_n\vert^{\beta}} +\int\limits_{\vert \zeta\vert \geq 1} \frac{
dV(\zeta)
  }
  { \vert \zeta\vert^{2n+2\gamma}  \vert \zeta_n\vert^{\beta}}
      < C,
\end{eqnarray*}
 where the  last  inequality  follows from  applying  twice  a variant of   Lemma 7.4. Finally,
\begin{eqnarray*}
 I_3\lesssim  \int\limits_{ \left\{  \frac{\vert z\vert }{2}\leq \vert \zeta \vert \leq
   2\vert z\vert \right\}                      }  \frac{
  \frac{ 1}
{\vert \zeta_{n} \vert^{\beta}  }
  }
  {  ( 1+ \vert z- \zeta\vert^{2} )^{
n+\gamma }}
dV(\zeta) \lesssim
\int\limits_{\vert z-\zeta\vert < 1} \frac{
dV(\zeta)
  }
  {  \vert \zeta_n\vert^{\beta}} +\int\limits_{\vert z- \zeta\vert \geq 1} \frac{
dV(\zeta)
  }
  { \vert z -\zeta\vert^{2n+2\gamma}  \vert \zeta_n\vert^{\beta}}
      < C,
\end{eqnarray*}
 where  the  last    inequality  holds by applying twice  a variant of Lemma 7.4 and  an obvious change of variable.
The   lemma is  now proved.
 \hfill $\square$\\

 In order to  prove  Theorem 7.2, we  need  the following
\begin{lem}  There is a  constant $C=C(n)>1$ such that for all points $z,z^{'}  \in\H_n,$ there is
   a  smooth curve  $\gamma=\gamma_{ z,
z^{'}}:
\lbrack 0,1\rbrack\to \H_n $ satisfying
$$\gamma(0)=z,\quad \gamma(1)=z^{'},\quad  \vert\gamma(t)\vert \leq \max{ \{\vert z\vert, \vert z^{'}\vert\} },\quad
\vert\gamma^{'}(t)\vert \leq C\vert z- z^{'}\vert.$$
 \end{lem}
 \begin{proof}   Suppose without loss of  generality that $\vert z\vert \geq  \vert z^{'} \vert.$  We  set
 $\hat{z}:=\vert  z^{'}\vert  \frac{ z}{\vert z\vert }.$  Then a little  geometric argument  shows that
 $\vert z^{'}- \hat{z} \vert  \leq \vert z -z^{'} \vert $ and
  $\vert z- \hat{z} \vert  \leq \vert z -z^{'} \vert. $
  Since  the
 group  $SO(n+1,\R)$ acts transitively  on   $\partial\M_n$, there  exists  a curve  $\gamma_1: \lbrack 0,1\rbrack \to \H_n$
 satisfying  $$\gamma_1(0)=\hat{z},\quad \gamma_1(1)=z^{'},\quad \vert \gamma_1(t)\vert =\vert z^{'}\vert,\quad
 \vert \gamma_1^{'}(t)\vert \leq  C(n)\vert z^{'}-\hat{z} \vert.$$
 Define
 $$
 \gamma_2(t):=  \left \lbrace
 \begin{array}{l}
 (1-2t)z+ 2t\hat{z},\quad\text{if}\  0\leq t\leq  \frac{1}{2};\\
 \gamma_1(2t-1),\quad\text{if}\   \frac{1}{2} \leq t\leq 1.\\
 \end{array}
 \right. $$
 It is  easy  to  see that  for every    $t\not =\frac{1}{2},$ the curve $\gamma_2(t)$   satisfies  all  the  properties stated  in  the lemma.
 To  conclude  the proof, it suffices to  approximate  in  $\H_n$  the curve $\gamma_2$ by a smooth
 curve $\gamma.$
  \end{proof}

\noindent{\bf  Proof of   Theorem 7.2.}
 We only  give  the proof for   the case $q>1.$
   For  every points $Z,Z^{'}  \in\H_N,$      consider   the  smooth  curve
   $$\gamma=\gamma_{Z,Z^{'}}:=\left (\gamma_{ x,x^{'}}, \gamma_{ z,z^{'}},\gamma_{ w,w^{'}} \right ):
\lbrack 0,1\rbrack\to \H_N=\B_l \times \H_n\times \H_m,$$
where $\gamma_{ z,z^{'}},\gamma_{ w,w^{'}}$ are given by   Lemma 7.6 and $\gamma_{ x,x^{'}}(t):=
(1-t)x+tx^{'}.$ Then it follows from
  Lemma 7.6 that  there is  a constant $C:=C(N)$  such that
  $$\gamma(0)=Z,\ \gamma(1)=Z^{'},\
\vert\gamma^{'}(t)\vert \leq C\vert Z- Z^{'}\vert,\   \vert\zeta\vert \leq \max{ \{\vert z\vert, \vert z^{'}\vert\} },\
\vert\eta\vert \leq \max{ \{\vert w\vert, \vert w^{'}\vert\} },     $$
where $ \Theta\equiv (\xi,\zeta,\eta)=\gamma(t) .$   Set
$$E:= \left \{ \Theta \in \H_N: \vert \Theta -Z \vert \geq 2C\vert Z^{'} -Z \vert\right \}.$$
On the one hand, for $\Theta \not \in E,$ using the definition of  $\Delta$
and Theorem 3.1, we  check easily  that
\begin{eqnarray*}
I_1   &:= &\int\limits_{\M_N \setminus E}
\Delta(\Theta,Z,Z^{'})^{q}   \left (\frac{\vert \zeta\vert \vert \eta\vert}
{\vert \zeta_{n+1}\vert \vert \eta_{m+1}\vert }\right )^{ 2(q-1)} dV(\Theta)\\
&\lesssim &
\int\limits_{\vert \Theta- Z\vert \leq 2C\vert Z^{'}  -Z\vert }
    \frac{ \left (  1 +  \frac{\vert z\vert }{\vert \zeta\vert}
 \right )^{2q} \left (  1 +  \frac{\vert w\vert}{\vert \eta\vert}
 \right )^{2q} }
{\vert  \Theta -Z \vert^{(2\vert N\vert -1)q}}\left (\frac{\vert \zeta\vert \vert \eta\vert}
{\vert \zeta_{n+1}\vert \vert \eta_{m+1}\vert }\right )^{ 2(q-1)} dV(\Theta) \\
&+ &
\int\limits_{\vert \Theta- Z^{'}\vert \leq 3C\vert Z^{'}  -Z\vert }
    \frac{ \left (  1 +  \frac{\vert z^{'}\vert}{\vert \zeta\vert}
 \right )^{2q} \left (  1 +  \frac{\vert w^{'}\vert }{\vert \eta\vert}
 \right )^{2q}}
{\vert  \Theta -Z^{'} \vert^{(2\vert N\vert -1)q}} \left (\frac{\vert \zeta\vert \vert \eta\vert}
{\vert \zeta_{n+1}\vert \vert \eta_{m+1}\vert }\right )^{ 2(q-1)} dV(\Theta) \\
&=& I_{11}+I_{12}.
\end {eqnarray*}
To estimate  $I_{11}$ and $I_{12},$ it suffices to  apply  part 2) of   Theorem 5.6   with  $i=2$ and
 Theorem  6.1 (1).   This can be  reduced to    majorizing  $I_{11}$ and $I_{12}$   by
$\int\limits_{ \vert \widetilde{ \Theta}\vert < C\vert Z  -Z^{'}  \vert }  I_{\alpha,a}(\widetilde{ \Theta}) dV(\widetilde{ \Theta}), $
where $$\alpha:= \left ( (2\vert N\vert -1)q,2(q-1),2(q-1),2(q-1),2(q-1) \right ).$$
 An application of  Lemma 7.4 shows that  the latter integral  is bounded from above  by
  $ C \vert Z  -Z^{'}  \vert^{ 2\vert N\vert +4 - ( 2\vert N\vert +3)q}.$  Hence
\begin{eqnarray}
I_1  \lesssim  \vert Z  -Z^{'}  \vert^{ 2\vert N\vert +4 - ( 2\vert N\vert +3)q}.
\end{eqnarray}
 On the  other hand, if $\Theta\in E,$    then  for  every   $ 0 \leq t \leq 1 $  and $\gamma:=\gamma_{ Z,
Z^{'}},$ we have that  $\vert \gamma(t)
-\Theta\vert \approx \vert \Theta-Z\vert.$      Therefore,  using the explicit  formula  of  $B_j(\Theta,Z)$
and  taking into  account the  properties   of the curve $\gamma$ stated  at the beginning of the   proof,
the Mean Value Theorem, applied  to  the functions  of  variable $Z:$  $\frac{B_j(\Theta,Z)}{    \vert \Theta - Z\vert^{2\vert N\vert}}, $     shows that

\begin{eqnarray}
  \Delta\left(\Theta, Z ,Z^{'}\right ) \lesssim
 \vert Z-Z^{'}\vert \frac{
 \left (  1 +  \frac{\vert z\vert}{\vert \zeta\vert}+
 \frac{\vert w\vert}{\vert \eta\vert}
 \right )^3
}{\vert \Theta  -Z\vert^{2\vert N\vert }}
+\vert Z-Z^{'}\vert \frac{
 \left (  1 +  \frac{\vert z^{'}\vert}{\vert \zeta\vert}+
 \frac{\vert w^{'}\vert}{\vert \eta\vert}
 \right )^3
}{\vert \Theta  -Z^{'}\vert^{2\vert N\vert }}
.
\end {eqnarray}
Proceeding exactly  as in estimating    $I_{11}$ and $I_{12},$ we  get
\begin{equation}\begin{split}
I_{21}&=   \left\vert Z -Z^{'}  \right\vert ^{q}
\int\limits_{E }
    \frac{ \left (  1 +  \frac{\vert z\vert}{\vert \zeta\vert}+
 \frac{\vert w\vert}{\vert \eta\vert}
 \right )^{3q} }
{\vert  \Theta -Z \vert^{2\vert N\vert q}}\left (\frac{\vert \zeta\vert \vert \eta\vert}
{\vert \zeta_{n+1}\vert \vert \eta_{m+1}\vert }\right )^{ 2(q-1)} dV(\Theta) \\
&\lesssim   \vert Z -Z^{'} \vert^{ 2\vert N\vert +4 - ( 2\vert N\vert +3)q}.
\end{split}
\end{equation}
Also,
\begin{equation}\begin{split}
 I_{22}&=   \left\vert Z -Z^{'}  \right\vert ^{q}
\int\limits_{E }
    \frac{ \left (  1 +  \frac{\vert z^{'}\vert}{\vert \zeta\vert}+
 \frac{\vert w^{'}\vert}{\vert \eta\vert}
 \right )^{3q} }
{\vert  \Theta -Z^{'} \vert^{2\vert N\vert q}}\left (\frac{\vert \zeta\vert \vert \eta\vert}
{\vert \zeta_{n+1}\vert \vert \eta_{m+1}\vert }\right )^{ 2(q-1)} dV(\Theta) \\
&\lesssim   \vert Z -Z^{'} \vert^{ 2\vert N\vert +4 - ( 2\vert N\vert +3)q}.
\end{split}
\end{equation}
Therefore,  it follows   from (7.3), (7.4) and (7.5) that
\begin{eqnarray*}
\int\limits_{ E}
\Delta(\Theta,Z,Z^{'})^{q}   \left (\frac{\vert \zeta\vert \vert \eta\vert}
{\vert \zeta_{n+1}\vert \vert \eta_{m+1}\vert }\right )^{ 2(q-1)} dV(\Theta)\lesssim I_{21}+I_{22}
\lesssim   \vert Z -Z^{'} \vert^{ 2\vert N\vert +4 - ( 2\vert N\vert +3)q}.
\end{eqnarray*}
This, combined with  estimate (7.2), completes the proof of the   theorem.\hfill
$\square$
\section{Lipschitz  estimates on the  complex manifold $\M_N$}
Let    $u$ be a function in $ \mathcal{C}^{1}(\M_N).$  For  every $Z\in\M_N,$   define
\begin{eqnarray*}
\left (\text{grad}_{\M_N} u\right )(Z):=\sup{ \vert (f\circ \gamma)^{'}(0) \vert},
\end{eqnarray*}
where the supremum being taken  over  all  smooth curves   $\gamma:\lbrack 0,1\rbrack\longrightarrow \M_N$
such that $\gamma(0)=Z$ and  $\vert \gamma^{'}(t) \vert \leq  1.$

We begin this  section with  the following   Hardy-Littlewood type lemma.
\begin{lem}
  For  every $ 0 < \alpha \leq 1,$   there  exists a constant $C= C(N,\alpha)$   with the following property:
    Suppose $u\in\mathcal{C}^{1}(\M_N)$  and  $K$ is  some  finite constant    such that
\begin{eqnarray*}
\left (\text{grad}_{\M_N} u\right )(Z)                      \leq            K\left (1 -\vert Z \vert \right )^{\alpha -1}\qquad
\text{for all}\ Z\in \M_N.
\end {eqnarray*}
Then $\vert u(Z) -u(Z^{'})\vert \leq CK\vert Z -Z^{'}\vert^{\alpha}$  for all $Z,Z^{'}\in \M_N.$
 \end{lem}
\begin{proof}
First  we make  the following   remark~:

{\it Write } $Z=(x,z,w),$ $Z^{'}=(x^{'},z^{'},w^{'}),$ $X:=(x,z),$ $Y:=w$ {\it and}  $X^{'}:=(x^{'},z^{'}),$ $ Y^{'}:= w^{'}.$
{\it  Suppose without loss of  generality that} $\vert Z^{'}\vert \leq  \vert Z\vert.$ {\it Then  } $\vert X^{'} \vert^2+ \vert Y^{'} \vert^2\leq \vert X \vert^2 +\vert Y \vert^2 <  1  .$
\\
 $\bullet$ {\it If} $\vert X^{'} \vert\leq\vert X \vert,$ {\it  by noticing that}
 $\vert (X^{'},Y) \vert \leq \vert Z\vert,$ {\it  then we write}
 $$ \vert u(Z) -u(Z^{'})\vert  =\vert u(X,Y) -u(X^{'},Y^{'})\vert \leq\vert u(X,Y) -u(X^{'},Y)\vert+
 \vert u(X^{'},Y) -u(X^{'},Y^{'})\vert.$$
 $\bullet$  {\it If } $\vert Y^{'} \vert\leq\vert Y \vert,$ {\it   by noticing that}
 $\vert (X,Y^{'}) \vert \leq \vert Z\vert,$ {\it  then we write}
 $$ \vert u(Z) -u(Z^{'})\vert  =\vert u(X,Y) -u(X^{'},Y^{'})\vert \leq\vert u(X,Y) -u(X,Y^{'})\vert+
 \vert u(X,Y^{'}) -u(X^{'},Y^{'})\vert.$$

 Let  $Z,Z^{'}$  be two  points of $ \M_N$   such that $0 <\vert Z^{'}\vert \leq \vert Z\vert <1$ and   set $\delta:= \vert Z -Z^{'} \vert.$

 First  assume that $\delta < 1 -\vert Z\vert.$
 Applying  the  previous remark  three times, we  only need   prove the lemma in one of  the
 following three cases:\\
  1)  $x=x^{'},z=z^{'};$ \qquad\qquad 2) $x=x^{'},w=w^{'};$ \qquad \qquad  3) $z=z^{'},w=w^{'}.$

 Suppose  for  example  we  are in the first case  $x=x^{'},z=z^{'}.$   In this case,
    take the curve $\gamma=\gamma_{Z,Z^{'}}.$
  According to the  hypothesis of the lemma and the properties of the curve   $\gamma$ given in the proof of
  Theorem 7.2, we  have
    \begin{eqnarray*}
\left (\text{grad}_{\M_N} u\right )(\Theta)                      \leq  K\delta^{\alpha-1},\qquad
\text{for all}\ \Theta\in \gamma(\lbrack 0,1\rbrack ).
\end {eqnarray*}
Therefore,
$$\vert u(Z) -u(Z^{'})\vert \leq  CK \delta^{\alpha-1} \vert Z -Z^{'}\vert        =CK \vert Z -Z^{'}\vert^{\alpha}.$$
  The remaining cases   $1-\vert Z \vert <\delta \leq 1-\vert Z^{'}\vert $ and
$1-\vert Z^{'}\vert < \delta $ can be checked  using the same argument as in Lemma
6.4.8 of \cite{Ru}.
\end{proof}

 In order to  state the main result of this section, we consider, for  $1\leq p <\infty,$ the space
\begin{equation*}
     L^p(\M_N):=\left\{ f:  \left (\int\limits_{\M_N}
     \vert f(\Theta)\vert^{p}   \frac{\vert \zeta_ {n+1}\vert^2  \vert \eta_ {m+1}\vert^2 }
     {\vert \zeta\vert^2  \vert \eta\vert^2 }   dV(\Theta)\right )^{\frac{1}{p} } \equiv  \Vert f
\Vert_{\M_{N,p}      }
 <\infty        \right\}.
    \end{equation*}
    If $f:= \sum\limits_{j=1}^{\vert N\vert +2}  f_j d\overline{\Theta}_j$ is a $(0,1)$-form defined in a neighborhood
     of $\overline{\M}_N$  in  $\overline{\B}_N,$  we  set
    \begin{equation}
\Vert f\Vert_{\M_{N,p}}:=\sum\limits_{j=1}^{\vert N\vert +2}
\Vert f_j \Vert_{\M_{N,p}      }      .
\end{equation}
Recall that the  norm  $\Vert \ \Vert_{\M_{N,\infty}      } $  was defined by  formula (2.8).

Next, for  every  $0 <\alpha < \beta \leq 1$  and  for $X= \M_N$ or $\partial\M_N,$  we define
\begin{eqnarray*}
     \Gamma_{\alpha,\beta }(X)  :=\left\{ f\,\, : \Vert f
\Vert_{\Lambda_{\alpha  } (X)} +
\underset{\underset{\gamma \subset X}{\gamma \in \mathcal{C}^{2}_{1}(\B_N)}  }{\sup} \Vert f\circ \gamma
\Vert_{\Lambda_{\beta }(\lbrack 0,1\rbrack)}   \equiv  \Vert f
\Vert_{ \Gamma_{\alpha,\beta}(X)     }
 <\infty        \right\}.
  \end{eqnarray*}
 We can say  informally that   $ \Gamma_{\alpha,\beta }(X) $ is the trace  of  the  non-isotropic Lipchitz
 space  $ \Gamma_{\alpha,\beta }(\B_N) $ (see Definition 1.1 in Krantz \cite{Kr2})  on the  manifold $X.$
 \begin{thm}
Suppose that $u \in \mathcal{C}^{1}(\overline{\M}_N)$  and consider   a
$(0,1)$-form
$$ f:= \left\lbrace
\begin{array}{l}
      \sum\limits_{k=1}^{\vert N\vert +2}   f_k
d\overline{\Theta}_k          ,\qquad\qquad\qquad\ \ \  \  \ \ \  \ \text{if}\  N \not =  (0,2,2);\\
  f_1 d\overline{\zeta}_1+f_2 d\overline{\zeta}_2+f_3 d\overline{\eta}_1+f_4 d\overline{\eta}_2,  \  \ \text{if}\  N =  (0,2,2);
  \end{array}  \right.$$
 with  coefficients in  $\mathcal{C}(\overline{\M}_N)$  such that  $\overline{\partial}_{\M_N}u=
f|_{\M_N}$   on  $  \M_N . $
  Define $T_1f$ on
$\partial\M_N$  as follows:
\\ $\bullet$ for $N\not= (0,2,2),$
 \begin{eqnarray*} (T_1f)(Z) :=  \int\limits_{\M_N}
\sum\limits_{k=1}^{\vert N\vert +2} \left\lbrack \frac {(1- \Theta\bullet\overline{Z})
P_{k}(\Theta,Z) +( 1-\vert \Theta\vert^{2})Q_{k}(\Theta,Z)}
{(1- \overline{\Theta}\bullet Z)^{\vert N\vert} (1-
\Theta\bullet\overline{Z})^{2}}\right\rbrack  f_k(\Theta)
\frac{dV(\Theta)}{  \vert \zeta \vert^2 \vert \eta \vert^2} ;
\end{eqnarray*}
 $\bullet$  for $N= (0,2,2),$ 
 \begin{eqnarray*}
 (T_1f)(Z) :=\int \limits_{\M_N}\sum\limits_{j=1}^{2} \sum\limits_{k=1}^{4}
\frac{(1- \Theta\bullet\overline{Z})^{1+j} }{  \left\vert 1-\Theta\bullet\overline{Z} \right\vert^{8}                      }\lbrack (
1-\Theta\bullet \overline{Z}) P_{jk}(\Theta,Z) +( 1-\vert \Theta
\vert^{2})Q_{jk}(\Theta,Z)\rbrack
   f_k(\Theta)\frac{ dV(\Theta)}{  \vert \zeta\vert^2   \vert \eta\vert^2},
\end{eqnarray*}
where    $P_k,Q_k$ and $P_{jk},Q_{jk}$ are the  polynomials given by Theorems 3.3 and 3.6.

Then the   definition of $T_1f$  can be extended to    $\M_N$  by setting
\begin{eqnarray}
(T_1f)(Z):=J_1(Z)+J_2(Z),
\end{eqnarray}
where
\begin{eqnarray*}
J_1(z)&:=&  \int\limits_{\partial\M_N}   \frac{A(\Theta,Z)}{ \vert  Z-\Theta \vert^{2\vert N\vert}}  (T_1f)(\Theta) \frac{  d\sigma(\Theta)
}{  \vert \zeta\vert^2   \vert \eta\vert^2},  \\
J_2(z)&:=& \int\limits_{\M_N}  \frac{1}{ \vert  Z-\Theta \vert^{2\vert N\vert}}  \left ( \sum\limits_{k=1}^{\vert N\vert +2} B_k(\Theta,Z)f_k(\Theta)\right )
 \frac{ dV(\Theta)}{  \vert \zeta\vert^2   \vert \eta\vert^2}
 \end{eqnarray*}
 and the operator $T_1f$ satisfies\\
 (i)   $\overline{\partial}_{\M_N}(T_1f)=f|_{\M_N}.$

 Moreover,  for  every    $p > 0,$  we  set  (as in the statement of  Theorem 1.1)~:
 $$\alpha=\alpha (N,p):= \left\lbrace
\begin{array}{l}
 \frac{1}{2}-\frac{\vert N\vert + 3}{p},\  \qquad \text{if}\   N \not =  (0,2,2)\  \text{and} \ p>2(\vert N\vert +3);\\
 \frac{1}{2}-\frac{6}{p},\qquad \ \  \   \  \ \text{if}\  N =  (0,2,2)\  \text{and} \ p> 12.
  \end{array}  \right.$$
Then there exists a constant C:= C(N,p)  such that
 \\
(ii)$$\left \lbrace
\begin{array}{l}
T_1f|_{\partial\M_N} \in \Gamma_{\alpha,2\alpha}(\partial\M_N)\ \text{and} \ \Vert
T_1f|_{\partial\M_N}        \Vert_{ \Gamma_{\alpha,2\alpha}(\partial\M_N)               }\leq  C\Vert f\Vert_{\M_{N,p}},\quad
\text{if}\ p<\infty; \\
T_1f|_{\partial\M_N} \in \Gamma_{\frac{1}{2},\tilde{1}}(\partial\M_N)\ \text{and} \ \Vert
T_1f|_{\partial\M_N}        \Vert_{ \Gamma_{\frac{1}{2},\tilde{1}         }(\partial\M_N)               }\leq  C\Vert f\Vert_{\M_{N,\infty}},\quad
\text{if}\ p=\infty;
\end{array}
\right.
$$
(iii)  $T_1f \in \Lambda_{\alpha}(\M_N)$ and $\Vert
T_1f\Vert_{ \Lambda_{\alpha}(\M_N)               }\leq  C\Vert f\Vert_{\M_{N,p}}.$
\end{thm}
\begin{proof}
We only  give  the proof  in the case     $N\not = (0,2,2)$  and $p <\infty.$  The first
 remaining case $N=(0,2,2)$ can be proved
in exactly the same way by  applying
   Theorem 3.6  instead  of  Theorem 3.3. The second remaining case
   $p=\infty$ follows essentially along the same lines as in our previous work
   \cite{VY1} basing on the work of Greiner-Stein \cite{GS}.

 We first introduce    two  new integral  operators  $T_2$ and  $T_3:$
 \begin{eqnarray*}
 (T_2f)(Z)  &:=&  \int \limits_{\M_N} \sum_{k=1}^{\vert N\vert +2}
\frac{(1- \Theta\bullet\overline{Z})^{\vert N\vert -2} }{ D(\Theta,Z)^{\vert N\vert}}\lbrack (
1-\Theta\bullet \overline{Z}) P_{k}(\Theta,Z) +   \\
&  &\qquad ( 1-\vert \Theta
\vert^{2})Q_{k}(\Theta,Z)\rbrack
   f_k(\Theta)\frac{ dV(\Theta)}{  \vert \zeta\vert^2   \vert \eta\vert^2}, \\
 (T_3f)(Z)  &:=& \int\limits_{\M_N}
\sum\limits_{k=1}^{\vert N\vert +2} \left\lbrack \frac {(1- \Theta\bullet\overline{Z})
P_{k}(\Theta,Z) +( 1-\vert \Theta\vert^{2})Q_{k}(\Theta,Z)}
{(1- \overline{\Theta}\bullet Z)^{\vert N\vert} (1-
\Theta\bullet\overline{Z})^{2}}\right\rbrack  f_k(\Theta)
\frac{dV(\Theta)}{  \vert \zeta \vert^2 \vert \eta \vert^2},\end{eqnarray*}
for all $Z\in\B_N.$

   Applying    Theorem 3.3
to the function  $u$ gives that
\begin{equation}
 (T_2f)(Z)=u(Z)-\int\limits_{\partial\M_N} \frac{R(\Theta,Z)}{  \left (1-Z\bullet \overline{\Theta}\right )^{\vert N\vert}}
u(\Theta)\frac{ d\sigma(\Theta)}{  \vert \zeta\vert^2   \vert \eta\vert^2}  , \qquad \text{for all}\ Z\in\M_N.
\end{equation}
Moreover, we note that
\begin{equation}
 (T_1f)(Z)=(T_2f)(Z)=(T_3f)(Z),
\qquad\text{for all}\  Z\in\partial\M_N.
\end{equation}
Arguing as in the proof of    Lemma 3.5 in \cite{VY1} and using Theorems  5.6, 6.1   and 7.1, one  can show that
    \begin{equation*}\lim\limits_{r\to 1^{-}}  \int\limits_{\partial\M_N}  \vert
(T_2f)(\Theta)-(T_2f)(r\Theta)\vert d\sigma(\Theta) =0.\end{equation*}
  Therefore, in view of   Remark 3.2, we can  apply  Theorem 3.1 to the  function $T_2f.$
Next  observe that   (8.2) is just   the Martinelli-Bochner formula.  Hence  by virtue of  (8.3), the hypothesis and
 the fact that   $R(\Theta,Z)$ is holomorphic  in the variable $Z,$ we obtain
\begin{eqnarray} T_1f=T_2f|_{\M_N}\quad  \text{and}\quad
\overline{\partial}_{\M_N}(T_1f)= \overline{\partial}_{\M_N} u     =f|_{\M_N}.
\end{eqnarray}
This completes the  proof of assertion (i).
In view of  (8.4), assertion (ii) will follow from the  following lemma.
\begin{lem}
\begin{equation*}
\left \lbrace
\begin{array}{l}
 \Vert
T_3f|_{\M_N}        \Vert_{ \Gamma_{\alpha,2\alpha}(\M_N)               }\leq  C\Vert f\Vert_{\M_{N,p}}\quad
\text{if}\ p<\infty; \\
 \Vert
T_3f|_{\M_N}        \Vert_{ \Gamma_{\frac{1}{2},\tilde{1}         }(\M_N)               }\leq  C\Vert f\Vert_{\M_{N,\infty}}\quad
\text{if}\ p=\infty;
\end{array}
\right.
\end{equation*}
\end{lem}
\begin{proof}
 Using the properties of the   polynomials $P_k(\Theta,Z)$  and  $Q_k(\Theta,Z)$ in  Theorem
3.3(ii), we  see that
\begin{eqnarray*}
\left \vert (\text{grad}\ T_3 f)(Z)  \right \vert  & \lesssim & \sum\limits_{k=1}^{\vert N\vert +2}   \int\limits_{\M_N}
 \frac{\vert \Theta - Z\vert }{ \left\vert 1 - \Theta \bullet \overline{Z}  \right \vert^{  \vert N\vert +2}}
   \left (  1 +  \frac{\vert z\vert }{\vert \zeta\vert }
 \right ) \left (  1 +  \frac{\vert w\vert }{\vert \eta\vert }
 \right ) f_k(\Theta)
       dV(\Theta) \\
       &+& \sum\limits_{k=1}^{\vert N\vert +2}   \int\limits_{\M_N}
 \frac{1 }{ \left\vert 1 - \Theta \bullet \overline{Z}  \right \vert^{  \vert N\vert +1}}
   \left (  1 +  \frac{\vert z\vert }{\vert \zeta\vert }
 \right ) \left (  1 +  \frac{\vert w\vert }{\vert \eta\vert }
 \right ) f_k(\Theta)
       dV(\Theta).
\end{eqnarray*}
Since $\vert \Theta  -Z\vert \leq 2\sqrt{ \left\vert 1 - \Theta \bullet \overline{Z}  \right \vert            },$ this implies by
  H\"{o}lder's inequality that   $ \left \vert (\text{grad}\ T_3 f)(Z)  \right \vert$
is bounded from above by
$$
C
\Vert f
\Vert_{\M_{N,p}      }
\left (\int\limits_{\M_N}
 \frac{1 }{ \left\vert 1 - \Theta \bullet \overline{Z}  \right \vert^{ ( \vert N\vert +\frac{3}{2})q}}
   \left (  1 +  \frac{\vert z\vert }{\vert \zeta\vert }
 \right )^{q} \left (  1 +  \frac{\vert w\vert }{\vert \eta\vert }
 \right )^{q}
   \left \vert\frac{ \zeta}
{ \zeta_{n+1}  }\right \vert^{\frac{2q}{p}}
 \left \vert\frac{ \eta}
{ \eta_{m+1}  }\right\vert^{\frac{2q}{p} }
       dV(\Theta) \right )^{\frac{1}{q}},
$$
where $q$  verifies  $\frac{1}{p}+\frac{1}{q}=1.$ Now applying    Theorem 7.1   yields
\begin{equation}
\left \vert (\text{grad}\ T_3 f)(Z)  \right \vert \leq C  \Vert f
\Vert_{\M_{N,p}      } \left (1  -\vert Z\vert \right )^{ -\frac{1}{2}-\frac{\vert N\vert +3}{p}}.
\end{equation}
so that by the classical  Hardy-Littlewood  lemma for the euclidean  ball  $\B_N$ we see  that
$$
\vert (T_3f)(Z)-  (T_3f)(Z^{'})      \vert \leq  C\Vert f\Vert_{\M_{N,p}}     \vert Z-Z^{'}\vert^{ \frac{1}{2}-\frac{\vert N\vert +3}{p}},\quad \text{for all}\ Z,Z^{'}\in\B_N.
$$
Therefore, choosing $Z^{'}=0,$  we obtain
\begin{equation}
\Vert T_3 f\Vert_{L^{\infty}(\B_N)} \leq \vert (T_3 f)(0)\vert + C\Vert f\Vert_{\M_{N,p}} \leq C\Vert f\Vert_{\M_{N,p}} .
\end{equation}
For every $u\in \mathcal{C}^{1}(\B_N),$   set
$$ \left(\text{grad}^{t}u\right) (Z):=\underset{\gamma\in \mathcal{C}^{2}_{1}(\B_N):\ \gamma(0)=Z}{\sup}{\vert (u\circ \gamma)^{'}(0)\vert }.
$$
By the proof  Lemma 4.8 in \cite{Kr1}, we see  that $\left  (\text{grad}^{t}\  \vert 1 -\overline{ \Theta} \bullet Z   \vert              \right )(Z)   \leq C  \vert \Theta -Z\vert$. Therefore,
a straightforward  calculation shows that
\begin{eqnarray*}
\left  (\text{grad}^{t}\ T_3 f\right )(Z)    & \lesssim & \sum\limits_{k=1}^{\vert N\vert +2}   \int\limits_{\M_N}
 \frac{\vert \Theta - Z\vert^{2} }{ \left\vert 1 - \Theta \bullet \overline{Z}  \right \vert^{  \vert N\vert +2}}
   \left (  1 +  \frac{\vert z\vert }{\vert \zeta\vert }
 \right ) \left (  1 +  \frac{\vert w\vert }{\vert \eta\vert }
 \right ) f_k(\Theta)
       dV(\Theta) \\
       &+& \sum\limits_{k=1}^{\vert N\vert +2}   \int\limits_{\M_N}
 \frac{1 }{ \left\vert 1 - \Theta \bullet \overline{Z}  \right \vert^{  \vert N\vert +1}}
   \left (  1 +  \frac{\vert z\vert }{\vert \zeta\vert }
 \right ) \left (  1 +  \frac{\vert w\vert }{\vert \eta\vert }
 \right ) f_k(\Theta)
       dV(\Theta).
\end{eqnarray*}
Hence, arguing as in the proof of  (8.6), we  see that
\begin{equation}
\left  (\text{grad}^{t}\ T_3 f\right )(Z)   \leq C  \Vert f
\Vert_{\M_{N,p}      } \left (1  -\vert Z\vert \right )^{ -\frac{\vert N\vert +3}{p}}.
\end{equation}
 Combining (8.6),(8.7) and (8.8), the  lemma  follows from Lemma 4.7 in \cite{Kr1}.
\end{proof}

To prove assertion (iii), we need the following
\begin{lem}
\begin{eqnarray*}
\vert J_1(Z)-J_1(Z^{'}) \vert \leq C\Vert f\Vert_{\M_{N,p}} \vert
Z-Z^{'}\vert^{\frac{1}{2} -\frac{\vert N\vert +3}{p}              }, \qquad \text{for all} \  Z,Z^{'}\in \M_N.
 \end{eqnarray*}
\end{lem}
\begin{proof}
Observe that  the polynomial  $A(\Theta,Z)$ satisfies
\begin{eqnarray*} \frac{1}{ \vert \zeta\vert^{2}  \vert \eta\vert^{2} }              \left\vert \text{grad}_Z \frac{A(\Theta,Z)}{ \left \vert \Theta -Z\right\vert^{2\vert N\vert }}  \right \vert
\lesssim  \frac{
    \left (   1+\frac{\vert z\vert}{\vert \zeta \vert} \right )
     \left (   1+\frac{\vert w\vert}{\vert \eta \vert} \right )}
    {   \left \vert \Theta -Z\right\vert^{2\vert N\vert }}
   +
 \frac{\frac{1}{\vert \zeta \vert}
                \left ( 1+  \frac{\vert w\vert}{\vert \eta \vert} \right )}
    {   \left \vert \Theta -Z\right\vert^{2\vert N\vert -1 }}
    +
\frac{\frac{1}{\vert \eta \vert}
                \left ( 1+  \frac{\vert z\vert}{\vert \zeta \vert} \right )}
    {   \left \vert \Theta -Z\right\vert^{2\vert N\vert -1 }}.
\end{eqnarray*}
In addition, if we  set $u\equiv 1$ in   Theorem 3.1,  then we see that
\begin{eqnarray*}
 \int\limits_{\partial\M_N}   \frac{A(\Theta,Z)}{ \vert  Z-\Theta \vert^{2\vert N\vert}}  \frac{  d\sigma(\Theta)
}{  \vert \zeta\vert^2   \vert \eta\vert^2}
=1.
\end{eqnarray*}
 Setting $Z:= rZ^{'},\ Z^{'}\in \partial\M_N,$  this implies that
\begin{eqnarray*}
 \left (\text{grad}_{\M_N} J_1\right ) (Z)& \lesssim &
 \int\limits_{\partial\M_N}      \left\vert \text{grad}_Z \frac{A(\Theta,Z)}{ \left \vert \Theta -Z\right\vert^{2\vert N\vert }}  \right \vert
  \left \vert
(T_1f)(\Theta) -(T_1f)(Z^{'}) \right \vert
 \frac{  d\sigma(\Theta)
}{  \vert \zeta\vert^2   \vert \eta\vert^2}.
\end{eqnarray*}
  Combining (8.4) and Lemma 8.3, we obtain
\begin{eqnarray*}
  \left \vert  (T_1f)(\Theta) -(T_1f)(Z^{'}) \right  \vert \leq C\Vert
f\Vert_{\M_{N,p}}  \vert \Theta - Z^{'}\vert^{\frac{1}{2} -\frac{\vert N\vert +3}{p}                     }.
\end{eqnarray*}
 Hence,
 \begin{eqnarray*}
  & & \left (\text{grad}_{\M_N} J_1\right ) (Z)              \leq  C\Vert
f\Vert_{\M_{N,p}} \left (
  \int\limits_{\partial\M_N}
\frac{
    \left (   1+\frac{\vert z^{'}\vert}{\vert \zeta \vert} \right )
     \left (   1+\frac{\vert w^{'}\vert}{\vert \eta \vert} \right )}
    {   \left \vert \Theta -rZ^{'}\right\vert^{2\vert N\vert }} \vert \Theta - Z^{'}\vert^{\frac{1}{2} -\frac{\vert N\vert +3}{p}                     }
 d\sigma(\Theta)  \right .\\
& +  & \left.
 \int\limits_{\partial\M_N}       \frac{\frac{1}{\vert \zeta \vert}
                \left ( 1+  \frac{\vert w^{'}\vert}{\vert \eta \vert} \right )}
    {   \left \vert \Theta -rZ^{'}\right\vert^{2\vert N\vert -1 }}\vert \Theta - Z^{'}\vert^{\frac{1}{2} -\frac{\vert N\vert +3}{p}                     }
 d\sigma(\Theta)
    +\int\limits_{\partial\M_N}
\frac{\frac{1}{\vert \eta \vert}
                \left ( 1+  \frac{\vert z^{'}\vert}{\vert \zeta \vert} \right )}
    {   \left \vert \Theta -rZ^{'}\right\vert^{2\vert N\vert -1 }}
 \vert \Theta - Z^{'}\vert^{\frac{1}{2} -\frac{\vert N\vert +3}{p}                     }
 d\sigma(\Theta)\right ).
 \end{eqnarray*}
 We shall establish in    Proposition 9.4 below  that the  latter three  integrals
 are dominated by   $C\left (  1 -\vert Z\vert \right )^{-\frac{1}{2} -\frac{\vert N\vert +3}{p}                     }.$
 Taking for granted  Proposition 9.4,  it follows that
 $$ \left (\text{grad}_{\M_N} J_1\right ) (Z)              \leq  C\Vert
f\Vert_{\M_{N,p}}\left (  1 -\vert Z\vert \right )^{-\frac{1}{2} -\frac{\vert N\vert +3}{p}                     }.$$
Finally, applying  Lemma 8.1  to this gives the desired conclusion.
\end{proof}

We now  complete the proof of   Theorem 8.2.
By   H\"{o}lder's inequality and        Theorem  7.2,  we have
\begin{eqnarray*}
  \left \vert  J_2 (Z) - J_2(Z^{'}) \right  \vert \leq C\Vert
f\Vert_{\M_{N,p}}  \vert Z - Z^{'}\vert^{1 -\frac{2\vert N\vert +4}{p}                     },\qquad
\text{for all}\  Z,Z^{'}\in\M_N  .
\end{eqnarray*}
This, combined with Lemma 8.4 gives that
\begin{eqnarray*}
  \left \vert  (T_1f)(Z) -(T_1f)(Z^{'}) \right  \vert \leq C\Vert
f\Vert_{\M_{N,p}}  \vert Z - Z^{'}\vert^{\frac{1}{2} -\frac{\vert N\vert +3}{p}                     },\qquad
\text{for all}\  Z,Z^{'}\in\M_N                  .
\end{eqnarray*}
Arguing as in the proof of   (8.7), one can show that
\begin{equation*}
\Vert T_1 f\Vert_{L^{\infty}(\M_N)}   \leq C\Vert f\Vert_{\M_{N,p}} .
\end{equation*}
This proves assertion (iii).
\end{proof}
\section{A   Stokes type theorem     on the manifold $\M_N$  and  applications}
The main result of this    section is the following   Stokes type theorem~:
\begin{thm}
Consider for every function $v\in \mathcal{C}^1\left (\overline{\M}_N\right)$ and every real numbers
$\lambda < 2n-1$ and $\mu< 2m -1,$  the function $u$ given by
$ u(\Theta):= \frac{v(\Theta)}{ \vert \zeta\vert^{\lambda}  \vert \eta\vert^{\mu}},$   for $\Theta \in \M_N.$
Then there is a constant  $C:=C(N)$ such that
\begin{eqnarray*}
\left \vert \int\limits_{\partial\M_N} ud\sigma  \right \vert \leq  C
\int\limits_{\M_N} \left (
\vert \xi \vert \vert (\text{grad}_{\xi} u )(\Theta) \vert +\vert \zeta \vert \vert (\text{grad}_{\zeta} u )(\Theta) \vert
+\vert \eta \vert \vert(\text{grad}_{\eta} u )(\Theta) \vert  +\vert u(\Theta)\vert  \right )  dV(\Theta).
\end{eqnarray*}
 \end{thm}
\begin{rem} {\rm  We do not know whether it is possible to establish  a theorem of reduction of estimates  from $\partial\M_N$
to  $\partial\B_{\vert N\vert},$ similar to Theorem 5.6.  To overcome this,  we use Theorem 9.1   to estimate  difficult  integrals     taken over    $\partial \M_N$ by  simpler ones  taken  over   $\M_N$ and then
apply Theorem  5.6.
    We have  already  encountered   this type of integral  estimates  in the proof of  Lemma 8.4.}
\end{rem}
\begin{proof}
 Set $d\xi  := d\xi_1\wedge \ldots  \wedge d\xi_l$ and
 \begin{eqnarray*}
\alpha_n(\zeta)        &:=  &  \frac{1}{n+1}\sum\limits_{j=1}^{n+1}  \frac{(-1)^{j-1}}{\zeta_{j}}
d\zeta_{1}\wedge
  \ldots \wedge  \widehat{d\zeta_{j}}\wedge
\ldots \wedge d\zeta_{n+1},\\
\alpha_m(\eta)        &:=  &  \frac{1}{m+1}\sum\limits_{k=1}^{m+1}  \frac{(-1)^{k-1}}{\eta_{k}}
d\eta_{1}\wedge
  \ldots \wedge  \widehat{d\eta_{k}}\wedge
\ldots \wedge d\eta_{m+1}.
 \end{eqnarray*}
 By   Proposition  2.1 in  \cite{VY1}  and Proposition 2.1 above, we  see that
 \begin{equation}
 \begin{split}
 dV_l(\xi)  &= C d\xi\wedge d\overline{\xi},\qquad\qquad\qquad\quad\quad \ \
 dV_n(\zeta)=  \left. C\vert \zeta\vert^2  \alpha_n(\zeta) \wedge \alpha_n(\overline{\zeta})\right |_{\H_n}, \\
 dV_m(\eta)&= \left. C\vert \eta\vert^2  \alpha_m(\eta) \wedge \alpha_m(\overline{\eta})\right|_{\H_m};
\qquad dV(\Theta)= dV_l(\xi)\wedge dV_n(\zeta) \wedge dV_m(\eta).
 \end{split}
 \end{equation}
Next, put
\begin{eqnarray*}
\omega_j(\xi) &:=& d\xi_{1}\wedge
  \ldots \wedge  \widehat{d\xi_{j}}\wedge
\ldots \wedge d\xi_{l},\qquad (1\leq j\leq l);\\
\omega_k(\zeta) &:=& d\zeta_{1}\wedge
  \ldots \wedge  \widehat{d\zeta_{k}}\wedge
\ldots \wedge d\zeta_{n+1},\qquad (1\leq k\leq n+1);\\
\omega_p(\eta) &:=& d\eta_{1}\wedge
  \ldots \wedge  \widehat{d\eta_{p}}\wedge
\ldots \wedge d\eta_{m+1},\qquad (1\leq p\leq m+1).\\
\end{eqnarray*}
Finally, we  define  $\omega_{jk}(\zeta),\widetilde{ \omega}_{jk}(\zeta)$ $(1\leq j,k\leq n+1)$  and
$\omega_{pq}(\eta),\widetilde{ \omega}_{pq}(\eta)$ $(1\leq p ,q\leq m+1)$  in just  the same  way as
$\omega_{pj}(z),\widetilde{ \omega}_{jk}(z)$  in \cite[p. 507--508]{MY}.

Consider the mapping $g: \rbrack 0,+\infty \lbrack \times \C^{\vert N\vert +2}
\longrightarrow \C^{\vert N\vert +2}$  given by  $$ g(t\Theta):= t\Theta.$$  Using (9.1)
and  proceeding as in the proof of     Lemma 2.1  in \cite{MY}, we   see  that
\begin{equation}\begin{split}
\left ( g^{\ast} dV \right )(t,\Theta)&=  t^{ 2\vert N\vert -1} dt\wedge \left\lbrack I_{\xi}\wedge dV_n(\zeta)\wedge dV_m(\eta) + I_{\zeta} \wedge dV_l(\xi)\wedge dV_m(\eta) \right.\\
& + \left . I_{\eta}\wedge dV_l(\xi)\wedge dV_n(\zeta)  \right \rbrack +  t^{ 2\vert N\vert } dV(\Theta),
\end{split}\end{equation}
where
\begin{eqnarray*}
I_{\xi}&:=& C \sum\limits_{p=1}^{l} (-1)^{p-1}\left (
  \overline{\xi}_p  d\xi \wedge \omega_p(\overline{\xi})+   \xi_p\omega_p(\xi)\wedge d\overline{\xi}  \right ),\\
I_{\zeta}&:= &C \vert \zeta\vert^2  \sum\limits_{j,k=1}^{n+1} \frac{ (-1)^{j+k}}{ \zeta_j \overline{\zeta}_k }
\widetilde{ \omega}_{jk}(\zeta),\qquad \text{and}\qquad
I_{\eta}:= C \vert \eta\vert^2  \sum\limits_{p,q=1}^{m+1} \frac{ (-1)^{p+q}}{ \eta_p \overline{\eta}_q }
\widetilde{ \omega}_{pq}(\eta).
\end{eqnarray*}
Now    set
\begin{equation}
\omega(\Theta):=I_{\xi}\wedge dV_n(\zeta)\wedge dV_m(\eta) + I_{\zeta} \wedge dV_l(\xi)\wedge dV_m(\eta)
 +  I_{\eta}\wedge dV_l(\xi)\wedge dV_n(\zeta).
\end{equation}
Since  $g$  is a diffeomorphism from  $\rbrack 0, +\infty\lbrack \times \partial\M_N\longrightarrow \H_N,$
it follows from (9.2) and  (9.3)  that
\begin{eqnarray*}
\int\limits_{\H_N}  u(\Theta) dV(\Theta)=\int\limits_{0}^{\infty}
t^{ 2\vert N\vert  -1} \int\limits_{\partial\M_N}  u(t\Theta) \omega(\Theta),\qquad \text{for all}\ u\in \mathcal{C}_{0}(\H_N),
 \end{eqnarray*}
 so that by   Lemma 2.3, we obtain  $ d\sigma = C\omega|_{\partial\M_N}.$ Therefore,
  since  by the  hypothesis  $\lambda < 2n-1,$  $\mu< 2m -1,$
$\vert u(\Theta)\vert \lesssim \frac{1}{ \vert \zeta\vert^{\lambda}  \vert \eta\vert^{\mu}}$   for  all $\Theta \in \M_N,$
   the    homogeneity properties  of the differential   form $\omega(\Theta)$ and the same
   arguments as in the  proof of (2.12), (2.13) and (2.17) of
 Proposition 2.5 imply that
 \begin{eqnarray}
 \int\limits_{\partial\M_N}  ud\sigma= \lim\limits_{ r \to 0} \int\limits_{\partial\M_{r    }     }  u\omega .
 \end{eqnarray}
  Stokes  theorem   gives that
 \begin{eqnarray}
 \int\limits_{\partial\M_{r   }     }  u\omega=\int\limits_{\M_{r   }     }  du \wedge\omega
 +\int\limits_{\M_{r    }     }  ud\omega.
 \end{eqnarray}
We shall estimate  $\left \vert \int\limits_{\M_r        }  du \wedge\omega \right \vert .$
Let  $Z$ be  a point of $\M_{r   }.$  Choose   $j$ and  $k$  with $1\leq j\leq n+1,\ 1\leq k\leq m+1$ so that
in a  sufficiently  small neighborhood  $\mathcal{U}:=\mathcal{U}(Z)$  in  $\M_{r    },$  we  have
\begin{equation}
\vert \zeta_j\vert \geq  \frac{1}{2} \max\limits_{p\not = j}{  \vert \zeta_p\vert }  \qquad
\text{and}\qquad  \vert \eta_k\vert \geq  \frac{1}{2} \max\limits_{q\not = k}{  \vert \eta_q\vert }.
\end{equation}
By (9.3) and  (9.5),  we  obtain
\begin{equation}\begin{split}
\left \vert \int\limits_{\mathcal{U} }  du \wedge\omega \right \vert  &\leq
\left \vert \int\limits_{\mathcal{U} }  du \wedge  I_{\xi}\wedge dV_n(\zeta)\wedge dV_m(\eta)            \right \vert
+\left \vert \int\limits_{\mathcal{U} }  du \wedge   I_{\zeta} \wedge dV_l(\xi)\wedge dV_m(\eta)                     \right \vert \\
& +\left \vert \int\limits_{\mathcal{U} }  du \wedge I_{\eta}\wedge dV_l(\xi)\wedge dV_n(\zeta)                       \right \vert.
\end{split}\end{equation}
We  shall estimate  for  example  $\left \vert \int\limits_{\mathcal{U} }  du \wedge   I_{\zeta} \wedge dV_l(\xi)\wedge dV_m(\eta)                     \right \vert.$    It should be noted that  the following   identity is  implicit
in the proof of  Lemma 2.1
of \cite{MY}~:
$$       \left . \frac{1}{ (n+1)^2 } I_{\zeta}\right|_{\H_n}=  \left .\frac{  (-1)^{j+k} \vert \zeta \vert^2 }{  \zeta_j  \overline{\zeta}_k}
\widetilde{\omega}_{jk}(\zeta)\right|_{\H_n},\qquad  \text{for}  \ 1\leq j,k\leq n+1.$$
 Therefore, $I_{\zeta}|_{\H_n}$ is  equal to
 \begin{eqnarray*}
  C\left .\frac{  \vert \zeta \vert^2 }{ \vert  \zeta_j  \vert^2}
\widetilde{\omega}_{jj}(\zeta)\right|_{\H_n}=
\vert \zeta\vert^2 \sum\limits_{p=1}^{j-1} \frac{  (-1)^{p-1} \zeta_p}{\zeta_j}  \omega_{pj}(\zeta)\wedge \alpha_n(\overline{\zeta})
 + \vert \zeta\vert^2 \sum\limits_{p=j+1}^{n+1} \frac{  (-1)^{p} \zeta_p}{\zeta_j}  \omega_{jp}(\zeta)\wedge \alpha_n(\overline{\zeta})\\
+\vert \zeta\vert^2 \sum\limits_{q=1}^{j-1} \frac{  (-1)^{n+q-1} \overline{\zeta}_q}{\overline{\zeta}_j}  \omega_{qj}
(\overline{\zeta})   \wedge \alpha_n(\zeta)  +
\vert \zeta\vert^2 \sum\limits_{q=j+1}^{n+1} \frac{  (-1)^{n+q} \overline{\zeta}_q}{\overline{\zeta}_j}  \omega_{jq}
(\overline{\zeta})   \wedge \alpha_n(\zeta).
\end{eqnarray*}
Combining the identity   $\alpha_n(\zeta)= (n+1)\frac{(-1)^{p-1}}{\zeta_p} \omega_p(\zeta), \ 1\leq p\leq n+1,$
(see formula (2.6) in \cite{MY}) and formula (9.1), a straightforward  calculation  gives that
 \begin{eqnarray*}
& &\int\limits_{\mathcal{U} }  du \wedge   I_{\zeta} \wedge dV_l(\xi)\wedge dV_m(\eta)
=
\int\limits_{\mathcal{U} }  d_{\zeta} u \wedge   I_{\zeta} \wedge dV_l(\xi)\wedge dV_m(\eta) \\
&=& C \int\limits_{\mathcal{U} } \sum\limits_{p=1,p\not = j}^{n+1} \left (
\frac{\partial u}{\partial \zeta_j} (-1)^p \zeta_p  +
\frac{\partial u}{\partial \zeta_p} (-1)^{p-1} \frac{\zeta_p^2}{\zeta_j} +
\frac{\partial u}{\partial \overline{\zeta}_j} (-1)^{p+n} \overline{\zeta}_p +
\frac{\partial u}{\partial \overline{\zeta}_ p} (-1)^{n+p-1} \frac{\overline{\zeta}_p^2}{\overline{\zeta}_j}  \right )
dV(\Theta).
\end{eqnarray*}
By virtue of  (9.6), we  majorize easily   the latter integral and   obtain
\begin{eqnarray*}
\left \vert \int\limits_{\mathcal{U} }  du \wedge   I_{\zeta} \wedge dV_l(\xi)\wedge dV_m(\eta)                     \right \vert \leq C(N)
\int\limits_{ \mathcal{U}       } \vert \zeta \vert \vert (\text{grad}_{\zeta} u )(\Theta) \vert
 dV(\Theta).
\end{eqnarray*}
Hence, in view of  (9.7), it follows that
\begin{equation*}
\left \vert \int\limits_{ \mathcal{U}    } du\wedge \omega  \right \vert \leq  C(N)
\int\limits_{ \mathcal{U}       } \left (
\vert \xi \vert \vert (\text{grad}_{\xi} u )(\Theta) \vert +\vert \zeta \vert \vert (\text{grad}_{\zeta} u )(\Theta) \vert
+\vert \eta \vert \vert(\text{grad}_{\eta} u )(\Theta) \vert   \right )  dV(\Theta).
\end{equation*}
On the other hand, we can  prove in just the same  way that
\begin{eqnarray*}
\left \vert \int\limits_{ \mathcal{U}    } ud \omega  \right \vert \leq  C(N)
 \int\limits_{ \mathcal{U}    } \vert u \vert dV.
\end{eqnarray*}
These two estimates, combined with  (9.4) and (9.5), complete  the proof.
\end{proof}

We now  present   two  applications of  Theorem
9.1.
\begin{prop}
Let   $\lambda,\alpha_1,\alpha_2$  be real numbers such that $ 0 < \lambda < 1,$    and
$0\leq \alpha_1  < 2n,$  $0\leq \alpha_2  < 2m.$  Then there   exists  a constant
$C:= C(N,\lambda,\alpha_1,\alpha_2) $   such that  for  every $ Z\in \B_{N},$
\begin{eqnarray*} \int\limits_{\partial\M_N}
\frac{1}{ \left\vert 1 - Z \bullet \overline{\Theta}   \right \vert^{  \vert N\vert +1-\lambda}}
\left (  1 +  \frac{\vert z\vert }{\vert \zeta\vert }
 \right )^{\alpha_1} \left (  1 +  \frac{\vert w\vert }{\vert \eta\vert }
 \right )^{\alpha_2}
       d\sigma(\Theta)\leq  C\left (1-\vert Z\vert^2\right )^{ \lambda -1}.
\end{eqnarray*}
\end{prop}
\begin{proof}
Applying  Theorem 9.1 gives that
\begin{eqnarray*} &  &\int\limits_{\partial\M_N}
\frac{1}{ \left\vert 1 - Z \bullet \overline{\Theta}   \right \vert^{  \vert N\vert +1-\lambda}}
\left (  1 +  \frac{\vert z\vert }{\vert \zeta\vert }
 \right )^{\alpha_1} \left (  1 +  \frac{\vert w\vert }{\vert \eta\vert }
 \right )^{\alpha_2}
       d\sigma(\Theta)\\
       &\leq&  C
\int\limits_{\M_N}
\frac{1}{ \left\vert 1 - Z \bullet \overline{\Theta}   \right \vert^{  \vert N\vert +2-\lambda}}
\left (  1 +  \frac{\vert z\vert }{\vert \zeta\vert }
 \right )^{\alpha_1} \left (  1 +  \frac{\vert w\vert }{\vert \eta\vert }
 \right )^{\alpha_2}
       dV(\Theta)\\
   & \leq  &    C\left (1-\vert Z\vert^2\right )^{ \lambda -1},
\end{eqnarray*}
where  the latter  inequality   follows from   Theorem 7.1.
\end{proof}
The following proposition completes the missing point  in the proof of Lemma 8.4 on page 40.
\begin{prop}
Suppose that   $ 0 < \lambda < 1.$  Then  there  is  a   constant
$C:= C(N,\lambda) $ such that  for every  $ 0< r< 1$ and $ Z\in \partial\M_{N},$
\begin{eqnarray*}
I_1 &:=& \int\limits_{\partial\M_N}
\frac{
    \left (   1+\frac{\vert z\vert}{\vert \zeta \vert} \right )
     \left (   1+\frac{\vert w\vert}{\vert \eta \vert} \right )}
    {   \left \vert \Theta -rZ\right\vert^{2\vert N\vert }} \vert \Theta - Z\vert^{\lambda                     }
 d\sigma(\Theta) \leq  C\left (1- r\right )^{ \lambda -1},\\
 I_2  &:= &
 \int\limits_{\partial\M_N}       \frac{\frac{1}{\vert \zeta \vert}
                \left ( 1+  \frac{\vert w\vert}{\vert \eta \vert} \right )}
    {   \left \vert \Theta -rZ\right\vert^{2\vert N\vert -1 }}\vert \Theta - Z\vert^{\lambda                     }
 d\sigma(\Theta)\qquad
    \leq C\left (1- r\right )^{ \lambda -1}.
\end{eqnarray*}
\end{prop}
\begin{proof}
We only  give  the proof of the  estimate for  $I_2.$   Starting from the elementary estimate
$ \vert \Theta - rZ \vert \approx (1-r) + \vert \Theta -Z \vert$ for all $ \Theta \in \partial \M_N,$
we  see that
\begin{eqnarray*}
I_2 &\lesssim  & \int\limits_{\partial\M_N}       \frac{\frac{1}{\vert \zeta \vert}
                \left ( 1+  \frac{\vert w\vert}{\vert \eta \vert} \right )}
    {   \left \lbrack        (1-r) +    \vert  \Theta -Z \vert \right\rbrack^{2\vert N\vert -1 }}\vert \Theta - Z\vert^{\lambda                     }
 d\sigma(\Theta)  \\
 &\lesssim  &
 \int\limits_{\M_N}       \frac{\frac{1}{\vert \zeta \vert}
                \left ( 1+  \frac{\vert w\vert}{\vert \eta \vert} \right )}
    {   \left \lbrack        (1-r) +    \vert  \Theta -Z \vert \right\rbrack^{2\vert N\vert -1 }}\vert \Theta - Z\vert^{\lambda   -1                   }
 dV(\Theta)  \\
&\lesssim  &
\int\limits_{\M_N}       \frac{  \left ( 1+  \frac{\vert z\vert}{\vert \zeta \vert} \right )
                \left ( 1+  \frac{\vert w\vert}{\vert \eta \vert} \right )}
    {   \left \lbrack        (1-r) +    \vert  \Theta -Z \vert \right\rbrack^{2\vert N\vert -1 }\vert \zeta-
    z \vert}\vert \Theta - Z\vert^{\lambda   -1                   }
 dV(\Theta),
\end{eqnarray*}
where the last two estimates  follow respectively from    Theorem 9.1  and
   the very elementary  inequality   $\frac{1}{\vert \zeta \vert} \leq
   \left ( 1+  \frac{\vert z\vert}{\vert \zeta \vert} \right )\frac{1}{\vert \zeta-
    z \vert       }.$  Hence, by part 3) of
  Theorem 5.6, the latter  integral is  dominated  by   $C\widetilde{I}_2,$ where
 \begin{eqnarray*}
 \widetilde{I}_2:= \int\limits_{\B_{2,\vert N\vert } }      \frac{\vert\widetilde{ \Theta} - \widetilde{Z}\vert^{\lambda   -1                   }  }
    {   \left \lbrack        (1-r) +    \vert  \widetilde{\Theta} -\widetilde{Z} \vert \right\rbrack^{2\vert N\vert -1 }\vert \widetilde{\zeta}-
    \widetilde{z} \vert}
 dV(\widetilde{\Theta}).
 \end{eqnarray*}
Dividing the domain of integration of  $\widetilde{I}_2$   into  two regions
$$\widetilde{E}_1:=
\left\lbrace \Theta \in \B_{2,\vert N\vert }:\ \vert  \widetilde{\Theta} -\widetilde{Z} \vert < 1-r \right\rbrace\quad
\text{and}\quad  \widetilde{E}_2:=
\left\lbrace \Theta \in \B_{2,\vert N\vert }:\ \vert  \widetilde{\Theta} -\widetilde{Z} \vert \geq 1-r \right\rbrace,$$
we thus  break $\widetilde{I}_2$ into  two corresponding  terms  $\widetilde{I}_{21}$ and   $\widetilde{I}_{22}.$
We then  apply  Lemma  7.4  to  estimate  each of  these  terms  and  obtain
\begin{eqnarray*}
 \widetilde{I}_{21}&\lesssim & \frac{ 1}{ (1-r)^{ 2\vert N\vert -1        }   }\int\limits_{ \widetilde{E}_1         }      \frac{\vert\widetilde{ \Theta} - \widetilde{Z}\vert^{\lambda   -1                   }  }
    {   \vert \widetilde{\zeta}-
    \widetilde{z} \vert}
 dV(\widetilde{\Theta}) \lesssim  (1-r)^{\lambda -1}, \\
  \widetilde{I}_{22}&\lesssim &\int\limits_{  \widetilde{E}_2      }      \frac{\vert\widetilde{ \Theta} - \widetilde{Z}\vert^{\lambda   -1                   }  }
    {      \vert  \widetilde{\Theta} -\widetilde{Z}  \vert^{2\vert N\vert -1 }\vert \widetilde{\zeta}-
    \widetilde{z} \vert}
 dV(\widetilde{\Theta})    \lesssim  (1-r)^{\lambda -1}.
 \end{eqnarray*}
 In summary, we  have
 \begin{eqnarray*}
 I_2  \lesssim \widetilde{I}_{2}= \widetilde{I}_{21}+ \widetilde{I}_{22}\leq  C  (1-r)^{\lambda -1},
  \end{eqnarray*}
which  completes the proof of  the proposition.
\end{proof}
\section{Proof of the main results}
In this  section  we prove    Theorems 1.1 and  1.2.  For this  purpose,
  we first     establish some    preparatory  results.

Consider  the holomorphic mapping  $F_N :\  \overline{\M}_N \longrightarrow
\overline{\Omega}_N \setminus \{0\} $ which maps every $Z\equiv (x,z,w)\equiv
  \left ( x_1,\ldots, x_l,z_1,\ldots,z_{n+1},w_1,\ldots,w_{m+1}  \right ),$ element of  $\overline{\M}_N$   to
 $$  F_N (Z):= \widetilde{Z}:=\left (\frac{x_1}{\sqrt{2}},\ldots, \frac{x_l}{\sqrt{2}}      ,z_1,\ldots,z_n,w_1,\ldots,w_m \right).$$
Recall that  $dV(\widetilde{\Theta} )$ is the canonical  volume form of  $\C^{\vert N\vert}.$   It follows from formula (5.2) in \cite{VY1} and formula
(9.1)  that
\begin{equation}
 \frac{  \vert \zeta_{n+1}\vert^2\vert \eta_{m+1}\vert^2}
 {\vert \zeta\vert^2\vert \eta\vert^2 }
dV(\Theta)= C F_N^{*} \left ( dV(\widetilde{\Theta} ) \right )
 ,\quad
\text{for}\  \Theta  \in \M_N\ \text{and}\ \widetilde{\Theta}=F_N(\Theta).
 \end{equation}
\begin{prop}
Consider a   $\overline\partial$-closed   $(0,1)$-form $f$
 of class  $\mathcal{C}^1$  defined  in a neighborhood   of $\overline{\Omega}_N.$ Then  the solution
 $T_1(F_N^{*}f)$    given by   Theorem 8.2 satisfies
 \begin{equation*}
   \left (  T_1(F_N^{*}f)      \right )(Z)  = \left (  T_1(F_N^{*}f)      \right )(Z^{'})       ,  \end{equation*}
   for  all
 $Z,Z^{'}   \in\M_N $ such that $F_N(Z)= F_N(Z^{'}).$
\end{prop}
\begin{proof}  Suppose that $f\in \mathcal{C}^{1}_{0,1}( r\Omega_N)$  for  some $r>1.$   Since
$r\Omega_N$ is pseudoconvex, there exists
a function  $u\in \mathcal{C}^{1}(\overline{\Omega}_N)$  such that $\overline\partial  u=f$ in
$\Omega_N.$   Therefore, it follows from   (8.3) and (8.5)  that
for every $  Z\in\M_N,$
 \begin{equation*}  \left (  T_1(F_N^{*}f)      \right )(Z) = (u\circ F_N)(Z)
 - \int\limits_{\partial\M_N} \frac{R(\Theta,Z)}{  \left (1-Z\bullet \overline{\Theta}\right )^{\vert N\vert}}
(u\circ F_N)(\Theta)\frac{ d\sigma(\Theta)}{  \vert \zeta\vert^2   \vert \eta\vert^2}
 .  \end{equation*}
 Using this  and   the explicit formula  of  $R(\Theta,Z),$ we see that the proof follows.
\end{proof}
\begin{thm}
For every $0<\lambda \leq  \frac{1}{2},$
 there is a constant  $ C:= C(N,\lambda)$ such that
\begin{equation*}
\Vert u\Vert_{ \Gamma_{\lambda, 2\lambda}(\M_N)}  \leq  C  \Vert u\Vert_{ \Lambda_{\lambda}(\partial\M_N)},
\end{equation*}
 for  all  functions  $ u$ in $ \mathcal{C}(\overline{\M}_N)$ which are holomorphic in $\M_N.$
\end{thm}
\begin{proof}
Consider the holomorphic function $U\in H(\B_N)$ defined by
\begin{equation}
U(Z):= \int\limits_{\partial\M_N} \frac{R(\Theta,Z)}{  \left (1-Z\bullet \overline{\Theta}\right )^{\vert N\vert}}
u(\Theta)\frac{ d\sigma(\Theta)}{  \vert \zeta\vert^2   \vert \eta\vert^2},\qquad \text{for all}\ Z\in \B_N.
  \end{equation}
Applying  Theorem 3.3  to the function $u$   yields
\begin{equation*}
U(Z)=u(Z),\qquad \text{for all}\ Z\in\M_N.
\end{equation*}
This  shows that  the theorem will follow  from the estimate
\begin{equation}
\Vert U\Vert_{ \Gamma_{\lambda, 2\lambda}(\B_N)}  \leq  C  \Vert u\Vert_{ \Lambda_{\lambda}(\partial\M_N)}.
\end{equation}
To prove this, observe by  (10.2) and formula (2) in \cite[Section 6.4.4]{Ru} that the    radial  derivative   $(\mathcal{R} U)$  of $U$  is given by
\begin{eqnarray*}
(\mathcal{R} U)(Z)             &=&  \int\limits_{\partial\M_N}   \left\lbrack    \frac{      \left( C\vert  \xi\vert^2  +  C\vert  \zeta\vert^2 +C\vert  \eta\vert^2 \right )
      \left (
(\vert \eta \vert^2  +w\bullet \overline{ \eta} ) z\bullet \overline{ \zeta}     +
  (\vert \zeta \vert^2  +z\bullet \overline{ \zeta} ) w\bullet \overline{ \eta}
   \right )     }
{ \left ( 1 - Z\bullet \overline{\Theta}\right )^{\vert N\vert } }  \right. \\
&+ & \left .\frac{
 N\left( C\vert  \xi\vert^2  +  C\vert  \zeta\vert^2 +C\vert  \eta\vert^2 \right )
      (\vert \zeta \vert^2  +z\bullet \overline{ \zeta} )
(\vert \eta \vert^2  +w\bullet \overline{ \eta} ) Z\bullet \overline{\Theta} }
{ \left ( 1 - Z\bullet \overline{\Theta}\right )^{\vert N\vert +1} }\right \rbrack u(\Theta) \frac{ d\sigma(\Theta)}{  \vert \zeta\vert^2   \vert \eta\vert^2}.
\end{eqnarray*}
Using this   and  arguing as in the proof of   Theorem 6.4.9 of  \cite{Ru},
it can be shown that
\begin{eqnarray*}
\left \vert (\mathcal{R} U)(Z)\right \vert
\leq  C  \Vert u\Vert_{ \Lambda_{\lambda}(\partial\M_N)}
\int\limits_{\partial\M_N}
\frac{1}{ \left\vert 1 - Z \bullet \overline{\Theta}   \right \vert^{  \vert N\vert +1-\lambda}}
\left (  1 +  \frac{\vert z\vert }{\vert \zeta\vert }
 \right )\left (  1 +  \frac{\vert w\vert }{\vert \eta\vert }
 \right ) d\sigma(\Theta).
       \end{eqnarray*}
       Therefore, by Proposition 9.3,
 $$ \left \vert (\mathcal{R} U)(Z)\right \vert                      \leq  C \Vert u\Vert_{ \Lambda_{\lambda}(\partial\M_N)}\left (1-\vert Z\vert^2\right )^{ \lambda -1},
 \qquad \text{for all}\  Z\in\B_N ;$$
  so that by Theorem 6.4.10 of \cite{Ru},  inequality (10.3) follows
  and the proof  is now complete.
\end{proof}
\begin{thm}
  For every $(0,1)$-form $f$ and   real numbers  $p,\alpha$ satisfying the  hypothesis of Theorem 8.2, we have
$$\left \lbrace
\begin{array}{l}
T_1f \in \Gamma_{\alpha,2\alpha}(\M_N)\ \text{and} \ \Vert
T_1f        \Vert_{ \Gamma_{\alpha,2\alpha}(\M_N)               }\leq  C\Vert f\Vert_{\M_{N,p}}\quad
\text{if}\ p<\infty, \\
T_1f \in \Gamma_{\frac{1}{2},\tilde{1}}(\M_N)\ \text{and} \ \Vert
T_1f      \Vert_{ \Gamma_{\frac{1}{2},\tilde{1}         }(\M_N)               }\leq  C\Vert f\Vert_{\M_{N,\infty}}\quad
\text{if}\ p=\infty,
\end{array}
\right.
$$
where the constant  $C:=C(N,p).$
\end{thm}
\begin{proof}
Let $\frac{1}{2} < r <1$   and set
\begin{eqnarray*}
r\M_N:=\left \lbrace Z\in\M_N:\  \vert Z\vert < r \right \rbrace \quad \text{and} \quad
r\partial\M_N:=\left \lbrace rZ,\ Z\in\partial\M_N \right \rbrace.
\end{eqnarray*}
Now  we define the  norm $\Vert f\Vert_{r\M_{N,p}}$  in  the same  way as
$\Vert f\Vert_{\M_{N,p}}$ given  in  (8.1)  by   substituting the  domain of integration
$\M_N$ by  $r\M_N.$  It is  obvious  that  $\Vert f\Vert_{r\M_{N,p}}\leq  \Vert f\Vert_{\M_{N,p}}.$

Applying  Theorem 8.2  to the complex manifold $r\M_N$  gives an integral operator
$T_r$  that  satisfies the following  properties~:
\begin{equation}
\overline{\partial}_{\H_N} (T_r f)= f\quad \text{on}\  r\M_N,
\end{equation}
and
\begin{equation}
(T_rf)|_{r\partial\M_N} \in \Gamma_{\alpha,2\alpha}(r\partial\M_N)\ \text{and} \ \Vert
(T_r f)|_{r\partial\M_N}        \Vert_{ \Gamma_{\alpha,2\alpha}(r\partial\M_N)               }\leq  C\Vert f\Vert_{\M_{N,p}} .\end{equation}
Setting
\begin{equation}
u(Z):= (T_1 f)(Z) - (T_r f)(Z),\qquad \text{for all}\ Z\in\M_N,
\end{equation}
  Theorem 8.2  and (10.4) imply that  $u$ is holomorphic  on $r\partial\M_N.$ On the other hand,  by
Theorem 8.2(ii) and  (10.5), (10.6), we get
\begin{eqnarray*}
\Vert
u|_{r\partial\M_N}        \Vert_{ \Lambda_{\alpha}(r\partial\M_N)               } \leq
\Vert
(T_r f)|_{r\partial\M_N}        \Vert_{ \Lambda_{\alpha}(r\partial\M_N)               }
+\Vert
(T_1 f)|_{r\partial\M_N}        \Vert_{ \Lambda_{\alpha}(r\partial\M_N)               }
\leq C
\Vert f\Vert_{\M_{N,p}}.
\end{eqnarray*}
Applying  Theorem 10.2 to this   estimate  yields
$$ \Vert
u|_{r\partial\M_N}        \Vert_{ \Gamma_{\alpha,2\alpha}(r\partial\M_N)               } \leq   C
\Vert f\Vert_{\M_{N,p}},$$
so that by   (10.5), we get   $T_1f =u+T_rf \in \Gamma_{\alpha,2\alpha}(r\partial\M_N)$ and
\begin{equation*}
 \Vert
(T_1 f)|_{r\partial\M_N}        \Vert_{ \Gamma_{\alpha,2\alpha}(r\partial\M_N)               }\leq  C\Vert f\Vert_{\M_{N,p}} .\end{equation*}
  Since  all   admissible curves $\gamma\in \mathcal{C}^{2}_{1}(\B_N)$ such that
$\gamma\subset \M_N$  lie  on some  manifold    $r\partial\M_N,$   the proof of the theorem is now
complete.
\end{proof}

\smallskip

\noindent {\bf Proof of  Theorem  1.1.}
Consider first  the case  where $f$  is a $\overline\partial$-closed   $(0,1)$-form
of class  $\mathcal{C}^1$   defined in a neighborhood    of  $\overline{\Omega}_N.$  The general case  will be treated later.

In view of  (8.1) and (10.1), it can be checked that
\begin{equation}
\Vert f\Vert_{L^p(\Omega_N)} = C(N,p) \Vert F^{*}_{N} f \Vert_{\M_{N,p}}.
\end{equation}
By Proposition 10.1, we can define the $\overline{\partial}$-solving operator  $T$  on $\Omega_N$  as
\begin{eqnarray}
  \left (Tf \right )(\widetilde{Z}):=\left ( T_1(F^{*}_N f) \right ) (Z),
             \end{eqnarray}
 for every  $\widetilde{Z}\in \Omega_N$ and $Z\in\M_N$  such that $F_N(Z)= \widetilde{Z}.$

 Combining   Proposition 10.1,  Theorem 10.3
and  equalities (10.7) and (10.8), we see that the operator   $ T          $ satisfies
\begin{eqnarray}
  \overline{\partial}(Tf)  &= & f\qquad  \text{on}\  \Omega_N; \\
  \Vert
F^{*}_N (T f)        \Vert_{ \Gamma_{\alpha,2\alpha}(\M_N)               }&\leq &  C\Vert f\Vert_{ L^p(\Omega_N)          }.
\end{eqnarray}
Let $\widetilde{Z}\equiv (\widetilde{x},\widetilde{z},\widetilde{w})$ and    $ \widetilde{Z^{'}}      \equiv (\widetilde{x^{'}},\widetilde{z^{'}},\widetilde{w^{'}}) $ be two  elements of      $\Omega_N.$  We shall show that there  exists
a constant $C:= C(N,p)$   such that
\begin{eqnarray}
\left\vert (T f) ( \widetilde{Z})-(T f) ( \widetilde{Z^{'}})  \right\vert \leq C \Vert f\Vert_{ L^p(\Omega_N)          }
\vert \widetilde{Z} - \widetilde{Z^{'}} \vert^{\alpha}.
\end{eqnarray}
Using the  remark made at the beginning   of the proof of   Lemma 8.1,  we  only need  prove  (10.11)  in
one of the  following three cases:\\
 1)  $\widetilde{x}=\widetilde{x^{'}},\widetilde{z}=\widetilde{z^{'}};$ \qquad\qquad 2) $\widetilde{x}=\widetilde{x^{'}},\widetilde{w}=\widetilde{w^{'}};$ \qquad \qquad  3) $\widetilde{z}=\widetilde{z^{'}},\widetilde{w}=\widetilde{w^{'}}.$

 Consider for example the case    $\widetilde{x}=\widetilde{x^{'}},\widetilde{w}=\widetilde{w^{'}}.$   In this case,  estimate  (10.11) becomes
\begin{eqnarray*}
\left\vert (T f) ( \widetilde{x},\widetilde{z},\widetilde{w}          )-(T f) (\widetilde{x},\widetilde{z^{'}},\widetilde{w} )  \right\vert \leq C \Vert f\Vert_{ L^p(\Omega_N)          }
\vert \widetilde{z} - \widetilde{z^{'} }\vert^{\alpha},
\end{eqnarray*}
 which can be proved by using (10.10) and  arguing as in the proof of  case 2 in
Section 5 of \cite{VY1}.

It remains to treat  the general case.  If merely $f\in
 L^p(\Omega_N), $   we can  regularize   $f$ by convolution
with a  $\mathcal{C}^{\infty}_{0}$ function of    sufficiently small  support. Then the same limiting argument  as
in     \cite[p. 361-362]{Ru}  shows that the
conclusion of the theorem   holds also for such  $f$. This completes the proof of the theorem.
\hfill
$\square$\\

\smallskip

\noindent{\bf  Proof of Theorem  1.2.}
First  suppose that $p<\infty.$ 
We break  the proof into  four  cases. In the course of the proof, we  shall  see that  the general case  can be reduced
to one  of  these  four cases.  In the sequel, we write  for every $Z\in\C^{\vert N\vert},$   $Z\equiv
(x, z,w)        \in\C^l\times\C^n\times \C^m.$

\smallskip

\noindent{\bf Case 1:} $n>2$ and $m\geq 2.$

For  every real number  $\lambda_0$   such that $\frac{1}{4}< \lambda_0  <\frac{1}{2},$   consider two  real numbers
 $\lambda,\mu >0 $  related by   $2\mu^{2}=\frac{1}{2}\left (  \frac{1}{2}  -\lambda \right )=  \frac{1}{2}  -\lambda_0.$
Let $ c\in \C$   such that $\vert c\vert \leq 1,$ and
consider the following elements of  $\Omega_N:$
 \begin{eqnarray*}
 Z_{\lambda,c} &:= &\left ( \underbrace{0,\ldots,0}_{l},\lambda,i\lambda,\mu c, \underbrace{0,\ldots,0}_{n-3},\frac{1}{2},\frac{i}{2},
   \underbrace{0,\ldots,0}_{m-2}  \right );\\
 Z_{\lambda_0,c}&:=& \left ( \underbrace{0,\ldots,0}_{l},\lambda_0,i\lambda_0,\mu c, \underbrace{0,\ldots,0}_{n-3},\frac{1}{2},\frac{i}{2},
   \underbrace{0,\ldots,0}_{m-2}  \right ).\\
   \end{eqnarray*}
Now  we put   $f:=\overline\partial u_0,$   where the function $u_0$ is given by
\begin{eqnarray*}
u_0(Z):=\frac{ \vert z_3 \vert^2}
{\left (  1  - \frac{  z_1}{2}+ \frac{ iz_2}{2}- \frac{ w_1}{2}+ \frac{ iw_2}{2}  \right )^{ \frac{1}{2} + \frac{\vert N\vert +3}{p}}},\qquad \text{for all}\ Z      \in\Omega_N.
\end{eqnarray*}
Then we have
\begin{eqnarray*}
f(Z):=\frac{  z_3   d\overline{z}_3  }
{\left (  1  - \frac{  z_1}{2}+ \frac{ iz_2}{2}- \frac{ w_1}{2}+ \frac{ iw_2}{2}  \right )^{ \frac{1}{2} + \frac{\vert N\vert +3}{p}}}.
\end{eqnarray*}
Suppose that $u$ is a  solution of the equation $ \overline\partial u=f$ on  $\Omega_N.$   Since $u-u_0$ is holomorphic  on $\Omega_N$ and $u_0( Z_{\lambda,0}) = u_0 (Z_{\lambda_0,0})=0,$  by Cauchy formula  we have
\begin{eqnarray*}
\frac{1}{2\pi }  \int\limits_{0}^{2\pi} u(Z_{\lambda,e^{i\theta}}) d\theta -u(Z_{\lambda,0}) &=&
\frac{ \vert \mu \vert^2}
{ (   \frac{  1}{2}  -\lambda )^{ \frac{1}{2} + \frac{\vert N\vert +3}{p}}} , \\
\frac{1}{2\pi }  \int\limits_{0}^{2\pi} u(Z_{\lambda_0,e^{i\theta}}) d\theta -u(Z_{\lambda_0,0}) &=&
\frac{ \vert \mu \vert^2}
{ (   \frac{  1}{2}  -\lambda_0 )^{ \frac{1}{2} + \frac{\vert N\vert +3}{p}}}.
\end{eqnarray*}
If $u\in\Lambda_{\alpha+\epsilon} (\Omega_N)$ for some  $\epsilon >0,$ then the  difference between
the two left hand sides  is $O(\vert  \lambda -\lambda_0 \vert^{\alpha +\epsilon} ).$  On the other hand,  the difference between
the  two  right sides  is   greater than  $ C\vert  \lambda -\lambda_0 \vert^{\alpha }.$
Letting $\lambda_0$ tend  to $\frac{1}{2},$  we reach a contradiction. Hence  $u\not\in\Lambda_{\alpha+\epsilon} (\Omega_N).$

It now remains to check that   $f\in L^{p-\epsilon}(\Omega_N)$ for all  $ \epsilon >0.$  Applying   (10.7)
and using the local coordinates
$\Phi^{z}$ and $\Phi^{w}$ of Theorem 4.1 with
$z:=(\frac{1}{2},\frac{i}{2},\underbrace{0,\ldots,0}_{n-1})\in\H_n$ and $w:=(\frac{1}{2},\frac{i}{2},\underbrace{0,\ldots,0}_{m-1})\in\H_m,$
 it follows that for every  $\epsilon\geq 0,$
$$
\Vert f\Vert^{p-\epsilon}_{L^{p-\epsilon} (\Omega_N)}=C\Vert F^{*}_{N}f\Vert^{p-\epsilon}_{\M_{N,p-\epsilon }}
\approx
\int\limits_{\mathcal{U}\cap  \B_{\vert N\vert} }\frac{ \vert z_2\vert^{p-\epsilon}
\vert z_3\vert^2  \vert z_4\vert^2}{  \vert 1 - z_1\vert^{ (p-\epsilon)\left (\frac{1}{2} + \frac{\vert N\vert +3}{p}\right )}}
dV(Z),$$
where $\mathcal{U}$  is a sufficiently small  neighborhood of the  point $(1,0,\ldots,0)\in \C^{\vert N\vert}$ and $dV(Z)$ is
  the Lebesgue measure of  $\C^{\vert N\vert}.$ We now  explain briefly
how the estimate $\approx$ in the latter line could be obtained. Indeed, using the local
 coordinates $\Phi^z$ and $\Phi^w,$ the function
$\vert \zeta_{n+1}\vert$ (resp. $\vert \eta_{m+1}\vert$) appearing in the
$\Vert\cdot\Vert_{\M_N,p}$ norm in (8.1) becomes the function $\vert z_3\vert$ (resp. $\vert
z_4\vert$) defined in $\C^{\vert N\vert}.$

By integration in    polar coordinates,  it is easy to  reduce the estimate of the latter    integral  to that of the following one
$$
\int\limits_{z\in\C:\vert  1- z\vert <1 }\frac{ dz\wedge d\overline{z}}{  \vert 1 - z\vert^{ 2- \frac{\epsilon(\vert N\vert +3)}{p}}}.
$$
From   this integral,  we see that    $f\in L^{p-\epsilon}(\Omega_N)$ for all $ \epsilon >0.$        This completes the proof in the first case. Furthermore,
we remark that the method presented here can  be applied to all domains  $\Omega_N$ where $N:=(n_1,\ldots,n_m)$
satisfies the condition $n_m > 2.$

\smallskip

\noindent{\bf  Case 2:} $l\geq 1$ and $n,m \geq 2.$

 Choose  $\lambda_0,\lambda$ and $\mu$ as in case  1.
Let $ c\in \C$ such that $\vert c\vert \leq 1$ and
consider   the following points of    $\Omega_N:$
 \begin{eqnarray*}
 Z_{\lambda,c} &:= &\left ( \mu c,\underbrace{0,\ldots,0}_{l-1},\lambda,i\lambda, \underbrace{0,\ldots,0}_{n-2},\frac{1}{2},\frac{i}{2},
   \underbrace{0,\ldots,0}_{m-2}  \right );\\
 Z_{\lambda_0,c}&:=& \left (\mu c, \underbrace{0,\ldots,0}_{l-1},\lambda_0,i\lambda_0, \underbrace{0,\ldots,0}_{n-2},\frac{1}{2},\frac{i}{2},
   \underbrace{0,\ldots,0}_{m-2}  \right ).\\
   \end{eqnarray*}
We set   $f:=\overline\partial u_0,$  where the function $u_0$ is given by
\begin{eqnarray*}
u_0(Z):=\frac{ \vert x_1 \vert^2}
{\left (  1  - \frac{  z_1}{2}+ \frac{ iz_2}{2}- \frac{ w_1}{2}+ \frac{ iw_2}{2}  \right )^{ \frac{1}{2} + \frac{\vert N\vert +3}{p}}},\qquad \text{for  all}\ Z      \in\Omega_N.
\end{eqnarray*}

The rest of the proof follows along the same lines as that of    case 1. Finally, we remark that
 the method  used in  this  second case   works also for    all  domains $\Omega_N$  where $N:=(n_1,\ldots,n_m)$
satisfies the condition $n_1=1$ and $ n_m >1 .$

\smallskip

\noindent{\bf Case 3:} $l=0$ and $n=m=2.$

For every  $\lambda_0$  such that   $\frac{1}{2\sqrt{2}} < \lambda_0
  <\frac{1}{\sqrt{2}},$ let $\lambda$ and $\mu$  be two positive  real numbers satisfying
   $\mu^{2}=\frac{1}{2}\left (  \frac{1}{\sqrt{2}}  -\lambda \right )=  \frac{1}{\sqrt{2}}  -\lambda_0.$
Let $ c\in \C$ such that $\vert c\vert \leq 1$ and
consider  the following     elements of  $\Omega_N$
 \begin{eqnarray*}
 Z_{\lambda,c} := \left (\mu c,0,\lambda,i\lambda \right ),\qquad\text{and}\qquad
 Z_{\lambda_0,c}:= \left ( \mu c,0,\lambda_0,i\lambda_0  \right ).
   \end{eqnarray*}
We set  $f:=\overline\partial u_0,$  where the function $u_0$ is defined by
\begin{eqnarray*}
u_0(Z):=\frac{ \vert z_1 \vert^2}
{\left (  1  - \frac{ w_1}{\sqrt{2}}+ \frac{ iw_2}{\sqrt{2} } \right )^{ \frac{1}{2} + \frac{6}{p}}},\qquad \text{for  all}\ Z\equiv (z_1,z_2,w_1,w_2)      \in\Omega_N.
\end{eqnarray*}
Proceeding as in  the proof of  case 1,  it can be checked that if a function $u$ satisfies $\overline\partial u=f$ then   $u\not\in\Lambda_{\alpha+\epsilon} (\Omega_N),\ \forall \epsilon >0.$                       It now remains
to establish that   $f\in L^{p-\epsilon}(\Omega_N)$ for all $\epsilon >0.$

 We first  apply   (10.7),
then use the local coordinates
$\Phi^{w}$ in Theorem 4.1 with $w:=\left(\frac{1}{\sqrt{2}},\frac{i}{\sqrt{2}},0\right)\in \H_2,$
and conclude that  for every $\epsilon\geq 0,$
$$
\Vert f\Vert^{p-\epsilon}_{L^{p-\epsilon} (\Omega_N)}
= C\Vert F^{*}_{N}f\Vert^{p-\epsilon}_{\M_{N,p-\epsilon }}
\approx
\int\limits_{\mathcal{U}\cap  \B_4 }\frac{ \vert z_2\vert^{p-\epsilon}
\vert z_3\vert^2  }{  \vert 1 - z_1\vert^{ (p-\epsilon)\left (\frac{1}{2} + \frac{6}{p}\right )}}
dV_4(Z).$$
Here $\mathcal{U}$ is a  sufficiently small neighborhood of the point $(1,0,0,0)$ in $ \C^4$ and $\B_4$ (resp. $dV_4(Z)$) is the      euclidean unit ball  (resp.
 the Lebesgue measure) of $\C^4.$

By integration in   polar  coordinates, the estimate of the latter   integral is reduced to that of  the integral
$$
\int\limits_{z\in\C:\vert  1- z\vert <1 }\frac{ dz\wedge d\overline{z}}{  \vert 1 - z\vert^{ 2- \frac{6\epsilon  }{p}}}.
$$
From  this integral  we conclude that    $f\in L^{p-\epsilon}(\Omega_N).$  The proof of the theorem is
complete
in this third case.
It should be noted that  this method is applicable to  all domains $\Omega_N$   where $N:=(n_1,\ldots,n_m)$
satisfies the condition $n_1=\cdots=n_m=2.$

\smallskip

\noindent{\bf  Case 4:}  $l=m=0$ and $n=2.$

In this  case  $\alpha(N,p)= \frac{1}{2} -\frac{3}{p}.$  Let $z$ be a  strongly convex point   of the boundary $\partial\Omega_N.$
It then follows from the work of Krantz in \cite[Section 6]{Kr3} that there exists  a   $(0,1)$-form  $f \in\mathcal{C}^{\infty}(\mathcal{U})
$ that  satisfies the  conclusion of the theorem  if    $\Omega_N$ is replaced by  $\mathcal{U}.$   Here $\mathcal{U}$
is an open strongly convex neighborhood of   $z$   in $\Omega_N\cup \{ z\}.$    In view of \cite{Kr3},  we  see easily that
the form  $f$ can be extended to  a form of  class $\mathcal{C}^{\infty}(\Omega_N)$ satisfying the conclusion
of the theorem.  The proof  is thus complete  in this last  case.

This  argument  also  shows that the Lipschitz  $\left (\frac{1}{2}+\epsilon\right )$-estimates  $(\epsilon >0)$  do not hold for  the case  $p=\infty.$   This completes the proof of  Theorem  1.2.\hfill
$\square$ \\

\smallskip

%
%
%
%
%
%
%
%
%
%
 Finally, we conclude  this paper  by some  remarks and  open problems.

1. It seems to be of some  interest  to  establish the $(L^p,L^q)$
type optimal  regularity for the   $\overline\partial$-equation on
$\Omega_N.$

2.   We conjecture that  the Lipschitz
$\widetilde{\frac{1}{2}}$-regularity   corresponding to the case
$p=\infty$
 in Theorem  1.1 is optimal.  More precisely,
this regularity can not be  improved to  Lipschitz  $\frac{1}{2}.$

3. Does there exist a natural way  to define  the Nevanlinna class
on the non-smooth domains $\Omega_N$ and  find  a related Blaschke
type condition that characterizes the zeroes of the functions  of
this class?

\end{document}